\documentclass[a4paper,11pt]{article}

\usepackage{mathrsfs} 

\usepackage{amsthm}
\usepackage{hyperref}
\usepackage{amssymb}
\usepackage{latexsym}
\usepackage{amsmath}
\usepackage{amsfonts}
\usepackage{mathabx}
\usepackage{subfig,graphicx}
\usepackage[utf8]{inputenc}

\usepackage[LGR,T1]{fontenc}%
\usepackage[french,english]{babel}
\usepackage{url}
\usepackage{enumerate}
\usepackage{hyperref}
\usepackage{xcolor}

\usepackage{fancyhdr}
\newtheorem{Theorem}{Theorem}[part]

\newtheorem{Proposition}{Proposition}[part]

\newtheorem{Lemma}{Lemma}[part]
\newtheorem{Corollary}{Corollary}[part]
\newtheorem{Remark}{Remark}[part]

\makeatletter \@addtoreset{equation}{section}

\@addtoreset{Definition}{section}

\@addtoreset{Theorem}{section}

\@addtoreset{Proposition}{section}

\@addtoreset{Property}{section}

\@addtoreset{Assumption}{section}

\@addtoreset{Corollary}{section}

\@addtoreset{Lemma}{section}

\@addtoreset{Remark}{section}

\@addtoreset{Example}{section}

\newcommand{\cA}{\mathcal{A}}

\newcommand{\cE}{\mathcal{E}}
\newcommand{\cF}{\mathcal{F}}

\newcommand{\cL}{\mathcal{L}}
\newcommand{\cM}{\mathcal{M}}
\newcommand{\cN}{\mathcal{N}}

\newcommand{\cP}{\mathcal{P}}

\newcommand{\cV}{\mathcal{V}}
\newcommand{\cW}{\mathcal{W}}
\newcommand{\cX}{\mathcal{X}}

\newcommand{\E}{\mathbb{E}}
\newcommand{\F}{\mathbb{F}}

\renewcommand{\P}{\mathbb{P}}

\newcommand{\R}{\mathbb{R}}

\def \proof{{\noindent \bf Proof. }}
\def \eproof{\hbox{ }\hfill$\Box$}

\newcommand{\ud}{\mathrm{d}}

\newcommand{\HYP}[1]
    {\ensuremath{({H#1} ) }}

\newcommand{\1}{{\bf 1}}
    
\newcommand{\set}[1]
    {\ensuremath{\{ #1 \}}}
    
\newcommand{\HP}[1] 
    {\ensuremath{\mathscr{H}^{#1}}}

\newcommand{\esp}[1]{\ensuremath{\mathbb{E} \! \left[#1\right] }}

\renewcommand{\Xi}[1]{X_{i #1}}

\def \a{\alpha}
\def \b{\beta}
\def \d{\delta}
\def \la{\lambda}
\def \ve{\varepsilon}
\def \g{\gamma}

\def \s{\sigma}

\def \t{\underline t}
\def \us{ \underline s}

\def \ets{\eta^\s} 
\def \etb{\eta^b}   
\def \tcX{\tilde{\cX}} 
\def \tnu{\tilde{\nu}} 
\def \eg{\mathfrak{e}} 

\def \vtstr{\vartheta^\star} %
\def \astr{\a}
\def \bstr{\b} 
\def \rstr{r^\star}
\def \pstr{p^\star}
\def \nstr{\nu^\star}

\newcommand{\NF}[1] 
    { \ensuremath{ \|#1\|_{_F} } } 

\def \hypl2{\HYP{^\star}}

\def \hypas{\HYP{^\star_{ \small as }}}

\newcommand{\ebig}[1]{\ensuremath{\mathbb{E}\bigg[#1 \bigg] }}

\definecolor{jcg}{HTML}{0aa344}

\title{Computing the invariant distribution of McKean-Vlasov SDEs by ergodic simulation
}

\author{
{\sc  
Jean-Francois Chassagneux \thanks{ENSAE-CREST and Institut Polytechnique de Paris, France. Email: {\tt chassagneux@ensae.fr}}~\thanks{This research benefited from the support of the ``Chaire Futures of Quantitative Finance''.} \quad
Gilles Pag\`es} \thanks{Laboratoire de Probabilit\'es, Statistique et Mod\'elisation, UMR~8001, Sorbonne Universit\'e, case 158, 4, pl. Jussieu, F-75252 Paris Cedex 5, France. E-mail: {\tt  gilles.pages@sorbonne-universite.fr}}~\thanks{This research
benefited from the support of the ``Chaire Risques Financiers'', Fondation du Risque.}  
}
\begin{document}
\maketitle

\begin{abstract}
    We design a fully implementable scheme to compute the invariant distribution of ergodic McKean-Vlasov SDE satisfying a uniform confluence property. Under natural conditions, we prove various convergence results notably we obtain rates for the Wasserstein distance in quadratic mean  and almost sure sense.
\end{abstract}



\section{Introduction}
Let $(\Omega,\cA,\P,\F:=(\cF_t)_{t\ge 0})$ be  a filtered probability space, with $\F$ satisfying the usual conditions. Let $q,d$ be two positive integers. We are  given a $q$-dimensional $\F$-Brownian motion $W$ (in particular independent of $\cF_0$) and 
we consider 
the following McKean-Vlasov Stochastic Differential Equation (MKV SDE), starting from $\xi \in L^2(\cF_0,\P)$, 
\begin{align}\label{eq de MKVSDE}
X_t = \xi + \int_{0}^tb(X_s,[X_s])\ud s + \int_0^t \sigma(X_s,[X_s]) \ud W_s, \quad t \ge 0\;,
\end{align}
where,  for a random variable $Y$, we denote by $[Y]$ its distribution. We assume that the function $(b,\sigma):\R^d\times\cP_2(\R^d) \rightarrow \R^d\times\cM_{d,q}(\R)$ satisfies
 the following assumption 
\begin{align}
    \HYP{L}\!:& \quad| b(x,\mu) - b(y,\nu) |  + \| \sigma(x,\mu) - \sigma(y,\nu)\|_{_F}
    \le L\left( |x-y| + \cW_2(\mu,\nu)\right) ,
\end{align}
for a positive constant $L$.
\textcolor{black}{Note that all starting values (denoted $X_0$, $X^\star_0$, $\cX_0$, $\bar\cX_0$, etc.) of stochastic differential equations driven by $W$  or any of their time discretizations schemes throughout the paper will be implicitely supposed to be $\cF_0$-measurable, except specific mention.
}

We know (see among others \cite{sznitman1991topics,bernou2022path}) that there exists a unique strong solution to \eqref{eq de MKVSDE}  and it satisfies for all $T\ge0$:
\begin{align}\label{eq basic int prop}
    \ebig{\sup_{t \in [0,T]}|X_t|^2} < +\infty.
\end{align}
An invariant distribution for the above dynamics is a probability measure $\nu$ such that
\begin{align}
    [X_t] = \nu, \; \forall \,t \ge 0.
\end{align}
In this work, we study the numerical approximation of the invariant distribution of \eqref{eq de MKVSDE}, if any (and when unique). 

\bigskip 

    The study and the computation of invariant measures for SDEs and more generally for continuous time Markov Processes is an important topic, which has motivated a lot of research: 
    Ranging from random mechanics \cite{kree1983random}, 
    stationary Fokker-Planck equation \cite{zbMATH00652291} to Langevin stochastic optimization \cite{zbMATH06775356,zbMATH07692273} see also \cite{MR4634732,MR4730248,MR1110084,MR4254553} for simulated annealing version. More recently, there has been some applications  to the modelling of stochastic volatility in finance \cite{zbMATH05816108}. 
    The case of McKean-Vlasov SDEs has also attracted a lot of attention as they have a large range of applications, from  natural sciences \cite{keller1971model, delarue2015particle} to economics \cite{jourdain2015capital,carmona2018systemic} as well as in physics or chemistry \cite{martzel2001mean, carrillo2003kinetic}, see also e.g. \cite{sharrock2021parameter,carmona2018probabilistic}. Invariant distributions appear also naturally in mean-field optimization problems with applications to neural network training \cite{chizat2018global,hu2021mean} or the learning of ergodic Mean-Field Games and Mean-Field Control \cite{lauriere2022learning}.  
    Our goal here is to present an approximation method of the invariant measure, which can be used in the previously mentionned examples.
   Precise applications to these fields are however outside of the scope of the present work. 
   
   \bigskip Before introducing the numerical method and stating our main results,
   we present the setting for our work, which guarantees among other things, the existence of a unique invariant distribution.

\bigskip
\noindent Our main structural assumption on $(b,\sigma)$ is a strong confluence assumption.
For every $p\ge 1$ and for $\a> \b\,\ge 0$, we consider:
 
\smallskip \noindent $\HYP{C}_{p, \a,\b}$: for every  $x$, $y\!\in \R^d$ and every $\mu$, $\nu\!\in {\cal P}_2(\R^d)$, $(b,\sigma)$ satisfies: 
\begin{align}\label{eq hyp C}
2(b(x,\mu)-b(y,\nu)\,|\, x-y\big) +(2p-1)\|\s(x,\mu)-\s(y,\nu)\|_{_F}^2 \le -\a|x-y|^2+\b {\cal W}_2^2(\mu,\nu).
\end{align}

\smallskip
\noindent When $p=1$, we shall write $\HYP{C}_{\a,\b}$ for $\HYP{C}_{1, \a,\b}$.

\smallskip
\noindent  We also introduce the $\HYP{MV}_{p, K,{\a},{\b}}$ mean-reverting assumption for  $p \ge 1, {\a}> 0$, $K$, ${\b}\ge 0$ as follows:

\smallskip \noindent $\HYP{MV}_{p, K, \a,\beta}$:  For every  $x$, $y\!\in \R^d$ and every $\mu \in {\cal P}_2(\R^d)$,
\begin{align}
    2(b(x,\mu)\,|\, x\big) +(2p-1)\|\s(x,\mu)\|_{_F}^2 \le K -{\a}|x|^2+{\b} {\cal W}_2^2(\mu,\d_0). \label{eq hyp mv}    
\end{align}

\medskip
We shall prove approximation results for the invariant distribution of \eqref{eq de MKVSDE} both in the $L^2$-sense and in the almost sure (a.s.) sense. We now formulate two standing assumptions which are required for our approach: 
\begin{align}
\hypl2\!:&  \text{ $\HYP{L}$ and $\HYP{C}_{\astr,\bstr}$ with $\astr>\bstr \ge 0$ hold and we set $\vtstr := 1-\frac{\bstr}{\astr} \in (0,1]$.} \label{eq Hstar for Gilles}
\\
        \hypas\!:&  \text{ \hypl2 and $\HYP{MV}_{2,K',\a',\b'}$ hold for $ \a'>\b' \ge 0$ and $K'\ge 0$.
        } \label{eq Hstaras for Gilles}
\end{align}

As stated in Corollary \ref{co main res invariant distrib mkv} of the next section, there exists a unique invariant distribution $\nstr$ to \eqref{eq de MKVSDE} under $\hypl2$. Moreover, we shall prove that there exists $p^\star>1$ such that  
\begin{align} \label{eq integrability inv meas first}
    \int_{\R^d}|y|^{2p^\star} \nu^\star(\ud y) < +\infty .    
\end{align}
For later use, we define $X^\star$ to be the solution of \eqref{eq de MKVSDE}, starting from some random vector $\xi\sim\nu^\star$. (Note that the notation makes no reference to $\xi$).

\noindent 
 To stress out the applicability of our main convergence results, we now formulate further assumptions on the $\sigma$ coefficients. These assumptions are used to obtain quantitative upper bound for the convergence error of our approximation method:

 \medskip
\noindent $\HYP{_\sigma}$: For every $(x,\mu)\in \R^d \times \cP_2(\R^d)$, $|\s(x,\mu)|\le  L_{\s}\big(1+|x|\eta_1(x) +{\cal W}_2(\mu, \d_0)\eta_2(\mu)\big)$ 
with 
\[ \lim_{|x|\rightarrow +\infty} \eta_1(x) = 0 \quad\text{ and }\quad \lim_{{\cal W}_2(\mu, \d_0)\to +\infty}\eta_2(\mu)=0,
    \]
where $L_\s$ is a positive constant,
\noindent  and

\medskip
\noindent $\HYP{_{\sigma_\infty}}$: There exists $\mu \in \cP_2(\R^d)$,  $\sup_{x \in \R^d}|\sigma(x,\mu)| < +\infty$.

\medskip
The connection between the assumptions $\HYP{L}$,$\HYP{C}$,$\HYP{MV}$ and $\HYP{_\sigma}$,$\HYP{_{\sigma_\infty}}$ is discussed in Lemma \ref{le structural assumptions3} below.

\noindent To approximate the invariant distribution of \eqref{eq de MKVSDE}, we consider the following numerical scheme, inspired by the decreasing step scheme introduced in \cite{lamberton2002recursive}. \\
Let $(\g_n)_{n\ge 1}$ be a non-increasing sequence of positive steps satisfying 
\begin{equation}\label{eq:stepdecrease}
\g_n \to 0 \quad \mbox{ and}\quad \Gamma_n:= \sum_{k=1}^n \g_k = +\infty \quad \mbox{ as}\quad n\to +\infty 
\end{equation}
(We set $\Gamma_0:=0$ by convention).

\smallskip \noindent
Denote $(\bar{\cX}_{\Gamma_n})_{n \ge 0}$ the scheme, which is defined by
\begin{align}
\label{eq:SDEdisc0}\bar\cX_0 &=\cX_0,\; \bar \nu_0= \d_{\cX_0} \text{ for some } \cX_0 \in \mathscr{L}^2(\cF_0,\P)\\
\label{eq:SDEdisc1}\bar \nu_{\Gamma_n} & = \frac{1}{\Gamma_n} \sum_{k=1}^n \g_k \d_{\bar \cX_{\Gamma_{k-1}}}, \quad n\ge 1,\\
\label{eq:SDEdisc2}\bar \cX_{\Gamma_n} & = \bar \cX_{\Gamma_{n-1}} +\g_n b(\bar \cX_{\Gamma_{n-1}}, \bar \nu_{\Gamma_{n-1}})+\sqrt{\g_n} \s(\bar \cX_{\Gamma_{n-1}}, \bar \nu_{\Gamma_{n-1}}) Z_n,\,\quad n\ge 1,\\
\label{eq:SDEdisc3}\mbox{with } \qquad Z_n   &= \frac{W_{\Gamma_{n}}-W_{\Gamma_{n-1}}}{\sqrt{\g_n}},\quad n\ge 1.
\end{align}
Other avatars of such decreasing schemes can be found in the literature:   See among others \cite{zbMATH02071135,zbMATH05196021} and more recently \cite{zbMATH06298329,zbMATH06525943,zbMATH06255773}.

\smallskip
The above algorithm is motivated by the following observation. It is known that the invariant measure (if any) of the SDE
\begin{align}
    \ud X_s = b(X_s,\mu)\ud s + \sigma(X_s,\mu) \ud W_s
\end{align}
can be obtained as the $a.s.$ weak limit of the random empirical measures $\nu^X_t(\omega) := \frac1t \int_0^t \delta_{X_u(\omega)} \ud u$, as $t\rightarrow \infty$. Inserting this estimation in the measure variable of the coefficients $(b,\sigma)$ leads to the following dynamics 
\begin{align}\label{eq de sid intro}
    \ud \cX _t =  b(\cX_t, \nu^\cX_t)\ud t + \sigma(\cX_t, \nu^\cX_t)\ud W_t\, \; \text{ with } \; \nu^{\cX}_t := \frac1t\int_0^t \delta_{\cX_s} \ud s.
\end{align}
The above process is known as a self-interacting diffusion (SID) and  $(\bar{\cX}_{\Gamma_n})_{n \ge 0}$ turns out to be its Euler scheme approximation with variable stepsize.

\medskip 
It is not the first time that self-interacting diffusions are used to approximate the invariant distribution of MKV SDEs. This idea has been suggested in \cite{alrachid2019new} and $L^2$-rate of convergence for the Wasserstein distance have been proved in \cite{du2023empirical}. The novelty of our paper is to focus on a fully implementable algorithm, compared to \cite{du2023empirical}. On top of $L^2$ and $a.s.$ convergence of our algorithm, we provide, for both cases, rates of convergence for the Wasserstein distance which is, to the best of our knowledge, completely new in this setting. On the road to our main results, which are given in the next section, we do study the convergence of the empirical measure of the self-interacting diffusion associated to an MKV SDE toward the invariant distribution of the latter. 
The ergodic behavior of self-interacting diffusion has already been proved under some specific conditions in \cite{benaim2002self, kleptsyn2012ergodicity} (see also the references therein) with physical applications in mind. In \cite{kleptsyn2012ergodicity}, the volatility is fixed and the drift is given by the gradient of the sum of an interaction potential and an external potential. Under convexity assumptions, an almost sure rate of convergence for the Wasserstein distance between the empirical measure of the SID and the invariant measure is provided. The methods used in this paper do rely on the structural convexity assumption: This is not the case in our approach where a confluence assumption is used instead. Moreover, the almost sure rates of convergence that we prove in Section \ref{subse conv emp meas}, though certainly suboptimal, are better than those proved in \cite{kleptsyn2012ergodicity}. Let us mention however that \cite{kleptsyn2012ergodicity} also treat (a specific example of) the case where $\astr=\bstr=1$, where uniqueness of the invariant distribution of the MKV SDE fails. The study of the convergence of our algorithm in this case seems out of reach by a direct application of the methods we develop below and is left for further research.
Regarding the convergence of the empirical measure of the SID to the invariant MKV SDEs, the closest paper to ours is certainly \cite{du2023empirical} (see also \cite{cao2024empirical} for the case of Langevin dynamics). In a similar setting, the authors prove the $L^2$-convergence of the Wasserstein distance with a rate. As a first step,  we recover here their results but we are also able to establish better rates of convergence under several stronger -- but reasonnable -- assumptions on the coefficients. 
\textcolor{black}{
It is important to notice that to prove the ergodic behavior of the SID, one must first control the rate of convergence in Wasserstein distance between the empirical measure of a standard SDE and its invariant distribution. 
}
\textcolor{black}{
    For this specific result, we rely on the companion paper \cite{chapag24}.
}

Let us also mention the work \cite{du2023self} that uses exponentially weighted empirical measure and combines this with an annealing method to obtain convergence but again no theoretical study of the implemented algorithm is provided.

\bigskip 

{
    \color{black}
The rest of the paper is organised as follows. In Section \ref{sub se main results}, we present ``ready to use'' convergence results for our numerical method. In section \ref{subse numerical illustration}, we present  some numerical experiments, which validate in practice our theoretical findings. In Section \ref{se stationnary solution MKVSDE}, we study the stationnary solution to MKV SDE.
In Section \ref{se SID}, we study the associated self-interacting diffusion (SID), which serves as an intermediary object to obtain the convergence of our numerical algorithm. In particular, we prove the stability of the SID with respect to some perturbations. In Section \ref{se numerical scheme}, we focus on the fully implementable algorithm and prove the main results announced above. 
Finally, the Appendix collects some useful technical results.
}

\subsection{Main results}
\label{sub se main results}
Our first main results is concerned with the convergence of the (discrete) empirical measure computed through our algorithm to the invariant distribution of the MKV~SDE.
\begin{Theorem}[Convergence of the empirical measure]\label{th main intro no rate}
    \begin{enumerate}[i)]
        \item Assume that $\hypl2$ hold. Then, 
        \begin{align*}
            \esp{\cW_2^2(\nu^\star,\bar{\nu}_{\Gamma_n})} \rightarrow 0 \text{ as } n \rightarrow +\infty.
        \end{align*}
        \item If moreover $\hypas$ hold and $\sum_{n\ge 1} \frac{\g_{n+1}^2}{\Gamma_n}<+\infty$, then 
         \[
            {\cal W}_2\big(\nu^{\star}   , \bar \nu_{\Gamma_n} \big)  \to 0  \quad \P\mbox{-}a.s.  \quad \mbox{ as }\quad n\to +\infty.
            \]
    \end{enumerate}
\end{Theorem}


\noindent We are able to provide upper bound for the rate of convergence in $L^2$ sense.
\begin{Theorem}[Quadratic Mean $\cW_2$-rates] \label{th main results intro L2 rate}
     Assume that $\hypl2$ hold. Let $\rstr := \frac{\vtstr}{1+\vtstr} \in (0,1)$ with $\vtstr$ defined in \eqref{eq Hstar for Gilles}. 
     For $n \ge 1$, we set $\gamma_n = \gamma_1 n^{-\rstr}$ with $\gamma_1>0$. Then, 
     \begin{enumerate}[i)]
        \item for any small $\eta>0$,
        \begin{align*}
            \esp{\cW_2^2(\nu^\star,\bar{\nu}_{\Gamma_n})}^\frac12 =
            O_\eta\left(n^{-(1-\rstr)\zeta_{\pstr}}\right) + \left(1+\cW_2(\nu^\star,[\cX_0])\right)o_\eta\left(n^{-\frac{\rstr}2+\eta}\right) 
        \end{align*}
        where $\zeta_{\pstr} := \frac{2 p^\star -1}{2(2(d+2)+(2p^\star-1)(d+3))} $.
        \item Assume moreover that $\HYP{\sigma}$ holds then  
        \begin{align*}
            \esp{\cW_2(\nu^\star,\bar{\nu}_{\Gamma_n})}^\frac12 
            = o_{\eta'}\left(n^{-\frac{1-\rstr}{2(d+3)}+\eta'}\right) + \left(1+\cW_2(\nu^\star,[\cX_0])\right)o_\eta\left(n^{-\frac{\rstr}2+\eta}\right) 
        \end{align*}
        for any small $\eta',\eta >0$.
        \item Assume moreover that $\HYP{\sigma_\infty}$ holds then
        \begin{align*}
            \esp{\cW_2^2(\nu^\star,\bar{\nu}_{\Gamma_n})}^\frac12 
            = O\left(n^{-\frac{1-\rstr}{2(d+3)}}\log(n)^{\frac{d+2}{2(d+3)}}\right) + \left(1+\cW_2(\nu^\star,[\cX_0])\right)o_\eta\left(n^{-\frac{\rstr}2+\eta}\right).
        \end{align*} 
     \end{enumerate}

\end{Theorem}

\noindent We now present  convergence results with a rate in the almost sure sense.

\begin{Theorem}[Almost sure $\cW_2$-rates] \label{th main results intro L2 as}
    Assume that $\hypas$ hold and that $\sigma$ is uniformly elliptic \textcolor{black}{(namely, $q=d$ and for some $\epsilon>0$, $\upsilon \sigma\sigma^\top(x,\mu)\upsilon \ge \epsilon $, for all $\upsilon\in \R^d$ with $|\upsilon|=1$ and $(x,\mu)\in \R^d \times \cP(\R^d)$)}. Let $\rstr := \frac{\vtstr}{1+\vtstr} \in (0,1)$.
    For $n \ge 1$, we set $\gamma_n = \gamma_1 n^{-(\rstr\wedge\frac13)}$ with $\gamma_1>0$. Then, 
    \begin{enumerate}[i)]
        \item for any small $\eta',\eta>0$, 
        \begin{align*}
            \cW_2(\nstr,\bar{\nu}_{\Gamma_n}) = o_{\eta'}\left(n^{-(1-(\rstr\wedge\frac13))\hat{\zeta}_{\pstr}}\log(n)^{\frac12+\eta'} \right) + 
            \left(1+|X^\star_0-\cX_0|\right)o_{\eta}\left(n^{-(\frac{\rstr}2\wedge\frac16)+\eta} \right)
        \end{align*}
        where $\hat{\zeta}_{\pstr} := \frac{(2\pstr-1)^2}{2(2\pstr+1) \{2(d+2) + (d+3)(d+2\pstr-1) \}}$.
        \item Assume moreover that $\HYP{\sigma}$ holds then, for any small $\eta',\eta>0$,   
         \begin{align*}
            \cW_2(\nstr,\bar{\nu}_{\Gamma_n}) =  o_{\eta'}\left(n^{- \frac{1-(\rstr\wedge\frac13)}{2(d+3)}+\eta'}  \right) + 
            \left(1+|X^\star_0-\cX_0|\right)o_{\eta}\left(n^{-(\frac{\rstr}2\wedge\frac16)+\eta} \right).
         \end{align*}
        \item Assume moreover that $\HYP{\sigma_\infty}$ holds then, for any small $\eta',\eta>0$,    
        \begin{align*}
            \cW_2(\nstr,\bar{\nu}_{\Gamma_n}) =  o_{\eta'}\left(n^{- \frac{1-(\rstr\wedge\frac13)}{2(d+3)}}\log(n)^{\frac12+\eta'}  \right) + 
            \left(1+|X^\star_0-\cX_0|\right)o_{\eta}\left(n^{-(\frac{\rstr}2\wedge\frac16)+\eta} \right).
         \end{align*}
    \end{enumerate}
\end{Theorem}

{  \paragraph{Comments:}
As developped in the following sections, there are mainly two sources of error appearing in the above upper bounds: the first one (in the r.h.s. of the previous inequalities) is due to the rate of approximation of the empirical measure of the self-interacting diffusion towards the invariant distribution of the MKV~SDE, the second is due to our approximation method \eqref{eq:SDEdisc0}-\eqref{eq:SDEdisc1}-\eqref{eq:SDEdisc2}-\eqref{eq:SDEdisc3}, namely the Euler scheme for the self-interacting diffusion. 

%

A natural question is: when does the algorithm slow down the convergence of the empirical measure toward the invariant distribution $\nstr$? Elementary computations show that, with our rates,  such is the case, roughly speaking if $r^{\star} = \frac{\vtstr}{1+\vtstr}$, see Theorem~\ref{th main results intro L2 rate}-(iii), satisfies $r^{\star} \le  \frac{1-r^{\star}}{d+3} $ i.e. 
$$
\frac{\b}{\a} \ge \frac{d+2}{d+5}.
$$
It is likely however that the rate order $\frac1{2(d+3)}$ is suboptimal and could be replaced by $\frac1d$ in line with results from \cite{fournier2015rate}.

The same phenomenon occurs in Theorem \ref{th main results intro L2 as}-(iii) involving $r^\star \wedge \frac13$ instead of $r^\star$. then, if $\frac{\beta}{\alpha} \ge \frac12$ then $r^\star \le \frac13$ leading to the same condition as in the $L^2$-case. If $\frac{\beta}{\alpha} <\frac12$ then this never occurs.
Let us notice that, for the $a.s.$ rate, there is a barrier at the order $\frac16$ coming from some stochastic integral terms.

\noindent 
Consequently  as $d$ increases,  this ``bad'' situation -- the algorithm slowing down the convergence -- occurs   less and less often. 

\vspace{2mm}
{
    \color{black}
\begin{Remark}
    The results stated in Theorem \ref{th main results intro L2 rate} and Theorem \ref{th main results intro L2 as} rely on the precise convergence rates obtained in our companion paper \cite{chapag24}, see Corollary \ref{co L2 conv rate emp meas to inv meas} and Corollary \ref{co as conv rate emp meas to inv meas} below. In particular, if one is able to prove that the empirical measure $\nu^{X^\star}_t = \frac1t \int_0^t \delta_{X^\star_s}\ud s$ converges towards the stationary probability distribution $\nu^\star$ with a rate $\zeta$, namely $\esp{\cW_2(\nu^{X^\star}_t,\nu^\star)^2}^\frac12=O(t^{-\zeta})$, then Theorem \ref{th main results intro L2 rate} i) would read 
        \begin{align}\label{eq with other rate to stationary nu}
            \esp{\cW_2^2(\nu^\star,\bar{\nu}_{\Gamma_n})}^\frac12 =
            O \left(n^{-(1-\rstr)\zeta}\right) + \left(1+\cW_2(\nu^\star,[\cX_0])\right)o_\eta\left(n^{-\frac{\rstr}2+\eta}\right) 
        \end{align}
        for any small $\eta>0$, under $\hypl2$.
        \\
Similarly, if one is able to obtain the following almost sure result: $ \cW_2(\nu^{X^\star}_t,\nu^\star) =O(t^{-\hat{\zeta}})$ then Theorem \ref{th main results intro L2 as} i) would read
\begin{align*}
    \cW_2(\nstr,\bar{\nu}_{\Gamma_n}) = O \left(n^{-(1-(\rstr\wedge\frac13))\hat{\zeta}} \right) + 
    \left(1+|X^\star_0-\cX_0|\right)o_{\eta}\left(n^{-(\frac{\rstr}2\wedge\frac16)+\eta} \right)
\end{align*}
for any small $\eta>0$, under $\hypas$.
\end{Remark}
}

}

{
    \color{black}
\subsection{Numerical illustration} 
\label{subse numerical illustration}

 We test the algorithm \eqref{eq:SDEdisc0}-\eqref{eq:SDEdisc1}-\eqref{eq:SDEdisc2}-\eqref{eq:SDEdisc3} on various one dimensional models.  In particular, we are interested in estimating the empirical convergence rate of the method for a given input step rate.  Thus for a step rate $r \in [0,1)$, namely $\gamma_n = \gamma_1 n^{-r}$, $1 \le n \le N$, for  a large $N$, we simulate  $M$ paths of the algorithm  to obtain $(\bar{\nu}_{\Gamma_n}^m)_{1 \le m \le M}$, $M$ independent estimations of the invariant probability measure. We then compute, for $1 \le n \le N$,
\begin{align}\label{eq de error}
    E_{n} := \left(\frac1M \sum_{m=1}^M \cW^2_2(\bar{\nu}_{\Gamma_n}^m,\tilde{\nu}^\star)\right)^\frac12 \text{ with } \tilde{\nu}^\star \text{ a discrete approximation of }\nu^\star.
\end{align}
When $M$ is large and $\tilde{\nu}^\star$ is close to $\nu^\star$, $E_n$ is a good approximation of $\esp{\cW^2_2(\bar{\nu}_{\Gamma_n} ,{\nu}^\star)}^\frac12$. We then observe on simulations that $E_n \sim n^{-\mathfrak{r}}$ for some $\mathfrak{r}\ge 0$ that we estimate empirically.

\noindent In all the experiments below, the initial condition of the algorithm $\bar{\cX}_0$ follows a uniform distribution on $[0,1]$ and we set $\gamma_1 = 1$.

\noindent Theorem \ref{th main results intro L2 rate} suggests to fix the step rate to $r^\star = \frac{\vartheta^\star}{1+\vartheta^\star}$ with $\vartheta^\star = 1-\frac{\beta}{\alpha}$. A natural question is then how to determine $\alpha$ and $\beta$? We will use the following simple answer.

\begin{Lemma}\label{le simple majo} Denote 
    $ a := - \sup_{x \in \R,\mu \in \cP_2(\R)} \partial_x b(x,\mu) > 0 $
    and assume that 
    \[
      a>0, \quad |b(x,\mu)-b(x,\nu)| \le L_b\cW_2(\mu,\nu) \,\text{ and }\, |\sigma(x,\mu)-\sigma(x,\nu)| \le L_\s\cW_2(\mu,\nu)
    \]
    for all $x \in \R,\mu,\nu \in \cP_2(\R)$.
    Then, $\HYP{C}_{\a ,\b}$ holds with 
    \begin{align*}
        \alpha   = 2a - \epsilon^\star L_b\text{ , } \beta =\frac{L_b}{\epsilon^\star} + L_\s^2 \text{ and }\epsilon^\star=\frac{2a}{L_b + \sqrt{L_b^2+2 a L_b L_\s^2}} .
    \end{align*}
    In particular,
    \begin{enumerate}
        \item if $L_\s = 0$ then $ \alpha   = a$, $\beta  = \frac{L_b}{a}$ and $\vartheta^\star = 1- \frac{L_b^2}{a^2}$.
        \item if $L_b = 0$, then $ \alpha    = 2a$, $\beta =L_\s^2$ and $\vartheta^\star = 1- \frac{L_\s^2}{2a}$.
    \end{enumerate}
\end{Lemma}


\paragraph{Mean-field Ornstein-Uhlenbeck model} 
For $\mathfrak{b} \ge 0$ and $\mathfrak{d} >0$, let us consider the following specification of \eqref{eq de MKVSDE}:
\begin{align}\label{eq de OU model}
    \ud X_t = (- X_t + \mathfrak{b}\esp{X_t})\ud t +  \sqrt{2} \mathfrak{d} \ud W_t.    
\end{align}
One proves easily that the behavior of $X_t$ is driven by the behavior of $m:t \mapsto \esp{X_t}$ that satisfies 
\begin{align}\label{eq evol mt ou like model}
    \ud m_t = (\mathfrak{b}-1)m_t \ud t.    
\end{align}
For $\mathfrak{b}<1$, the process has a unique stationary distribution given by $\cN(0,\mathfrak{d}^2 )$. For $\mathfrak{b}=1$, the process has  infinitely many stationary distributions parametrised by $\esp{X_0}$ (the mean of the initial condition).
If $\mathfrak{b}>1$, the stationary measure is still $\cN(0,\mathfrak{d}^2 )$: But if $X$ starts from a distribution with non zero mean, $(m_t)$ diverges. Concerning the associated self-interacting diffusion and for $\mathfrak{b}<1$, \cite{kleptsyn2012ergodicity} proves the convergence of $\nu^\cX_t$ to $\cN(0,\mathfrak{d}^2)$, see also our results below. The critical case $\mathfrak{b}=1$ is investigated in \cite{kleptsyn2012ergodicity}, which shows the convergence to a random measure. We first report  numerical results obtained for the model \eqref{eq de OU model} considering $\mathfrak{b}\in[0,1)$: A case covered by our theoretical results. 

\noindent As mentioned, we tested the algorithm \eqref{eq:SDEdisc0}-\eqref{eq:SDEdisc1}-\eqref{eq:SDEdisc2}-\eqref{eq:SDEdisc3} for various value of $\mathfrak{b}$ and various specifications of the algorithm step rate, with initial condition following an uniform law on $[0,1]$ and $\mathfrak{d}=1$ for the model \eqref{eq de MKVSDE}. 

\begin{figure}[h]
    \begin{center}
    \includegraphics[width=0.8\textwidth]{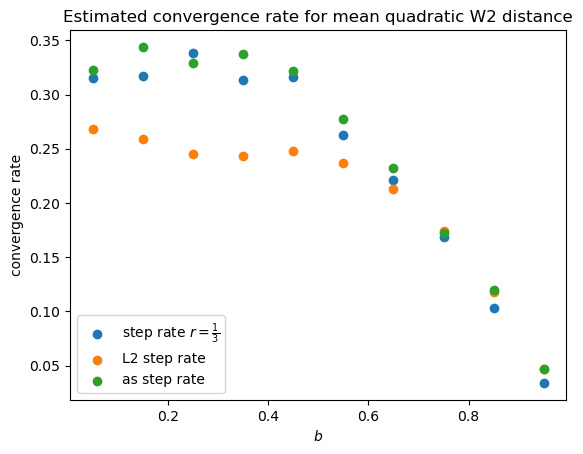}
    \end{center}
    \caption{Empirical estimation of the convergence rate  for Model \eqref{eq de OU model} as a function of $\mathfrak{b}$. Tests realised for three different specifications of the algorithm input step rate, $N=100000$ and $M=500$ in \eqref{eq de error}. $\bar{\cX}_0$ follows a uniform distribution on $[0,1].$} 
    \label{fig conv rate OU model various b}
\end{figure}

On Figure \ref{fig conv rate OU model various b}, we report the estimated convergence step rate of the numerical method for some values of $\mathfrak{b}\in (0,1)$. We use three specifications for the input step rate:
\\
- '\verb?L2 step rate?' is the rate obtained by using Lemma \ref{le simple majo}: here, we simply have that $\vartheta^\star = 1-\mathfrak{b}^2$ and $r^\star = \frac{1-\mathfrak{b}^2}{2-\mathfrak{b}^2}$ as given in Theorem \ref{th main results intro L2 rate}.\\
- '\verb?as step rate?' is the previous rate capped at $\frac13$ as suggested by Theorem \ref{th main results intro L2 as}.
\\
- '\verb?step rate =1/3?' corresponds to the optimal step rate obtained for classical Brownian SDEs when there is no mean field interaction, see e.g. \cite{lamberton2002recursive}.
\\
For small values of $\mathfrak{b}$, we observe that the rate $r=\frac13$ returns a convergence rate of $\frac13$ (approximately), from which we deduce that in \eqref{eq with other rate to stationary nu}, $\zeta$ must be equal to $\frac12$. So, for our experiments, we are in the case where the algorithm may slow down the convergence (see our Comments above). This appears on the graph when $\mathfrak{b}$ becomes large enough and corresponds to a strong mean field interaction (or self interacting feature).  Even though the observed rates are better than the theoretical ones, this slow down in the convergence rate is consistent with what is predicted by our theoretical results. Finally, we also observe that using the step rate $r=\frac13$ does not deteriorate dramatically the convergence. 

\vspace{2mm}

\begin{figure}
    \centering
    \subfloat{\includegraphics[width=.45\linewidth]{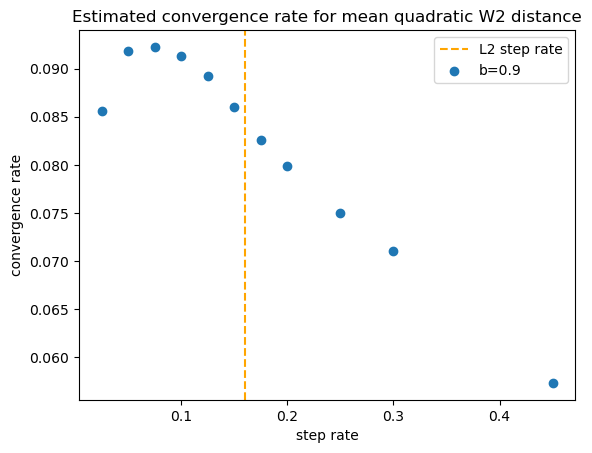}}\hfill
    \subfloat{\includegraphics[width=.45\linewidth]{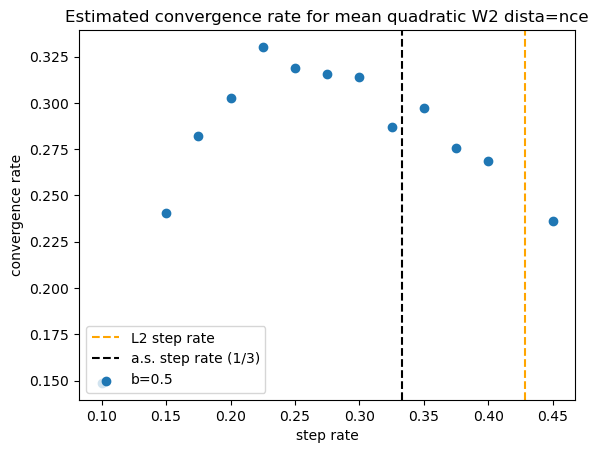}}\hfill 
     \caption{\footnotesize Comparison of optimal rates for $\mathfrak{b}=0.9$ (left) and $\mathfrak{b}=0.5$ (right) for model \eqref{eq de OU model}. $N=100000$ and $M=500$ in \eqref{eq de error}. $\bar{\cX}_0$ follows a uniform distribution on $[0,1].$
     }
     \label{fig optim}
\end{figure}

On Figure \ref{fig optim}, we plot the empirical convergence rate for various values of the step rate $r \in (0,0.5)$, again for two different values of $\mathfrak{b}$. In both cases, an optimal convergence rate seems to exist. In the case of a strong mean field interaction i.e. $\mathfrak{b}=0.9$, the optimum is localised around $0.1$ for a convergence rate slightly above $0.09$. In the case of a weaker mean field interaction, namely $\mathfrak{b}=0.5$, the optimal convergence rate is slightly below $1/3$ and seems to be realised for a range of input rate around $0.25$.   We observe that the realised optima are not those that we have predicted by the simple approach of Lemma \ref{le simple majo} and that we have used in Figure \ref{fig conv rate OU model various b}.
Because it is not so easy to predict the optimum algorithm step rate, we will use for our future experiments the fixed rate $r=\frac13$, which seems to be a good compromise.

\begin{figure}
    \centering
    \subfloat[Empirical mean of $(\bar{\cX}_{\Gamma_n})$ (M=$10000$) for various values of $\mathfrak{b}$]{\includegraphics[width=.5\linewidth]{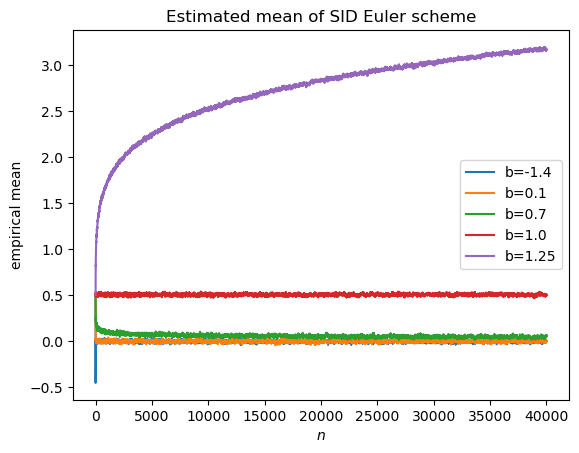}}\hfill  
    \subfloat[Convergence rate as a function of $\mathfrak{b}$. ($M=500$)]{\includegraphics[width=.45\linewidth]{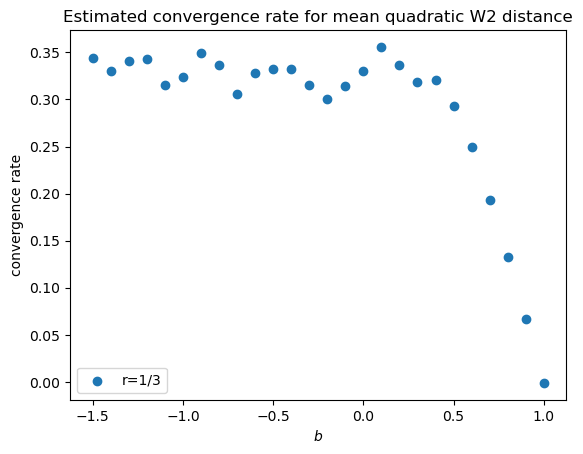}}\hfill  
     \caption{\footnotesize Behavior of $(\bar{\cX_n})$ and $(\bar{\nu}_{\Gamma_n})$ for model \eqref{eq de OU model}. $\bar{\cX}_0$ follows a uniform distribution on $[0,1].$}\label{fig final ou}
\end{figure}

With this choice of step rate ($r=\frac13$), we obtain Figure \ref{fig final ou}. On Figure \ref{fig final ou}(a), we plot the empirical mean of $(\bar{\cX}_n)_{0\le n \le N}$ (with $N=40000$ and the number of simulation $M=10000$). We observe a   behavior similar to the one of $(m_t)$ given in \eqref{eq evol mt ou like model}. In particular, for $\mathfrak{b}>1$, the empirical mean diverges and there is no hope for an ergodic behavior. The limit case $\mathfrak{b}=1$ exhibits a stationary behavior but Figure \ref{fig final ou}(b) shows that $(\bar{\nu}_{\Gamma_n})$ is not converging to $\cN(0,1)$, as expected from the theoretical SID behavior. For $\mathfrak{b}<1$, the empirical mean stabilises around $0$. On Figure \ref{fig final ou}(b), we observe convergence as predicted by our theoretical results for the range $\mathfrak{b}\in (-1,1)$ (see also the discussion above). But, interestingly, the algorithm is converging for $\mathfrak{b}\le -1$ with the optimal $\frac13$ rate. In this case, the mean-field feature does not perturb (possibly helps) the convergence: This case is not covered by our results below and its study is left for further research.

\paragraph{Non linear mean field Ornstein-Uhlenbeck model}
We now consider the model 
\begin{align}\label{eq fancy drift new formulation}
    \ud X_t = (-X_t + \mathfrak{b} \set{\|X_t\|_2-\mathfrak{d}})\ud t + \sqrt{2}\mathfrak{d} \ud W_t,
\end{align}
for parameters $\mathfrak{b} \in \R$, $\mathfrak{d}>0$.
\\
We make the following observation, see the appendix for computations.
\begin{Lemma}\label{le fancy drift stationary}
    The following holds 
    \begin{enumerate}
        \item  if $|\mathfrak{b}|\le 1$, there is a unique stationary measure to \eqref{eq fancy drift new formulation} given by $\cN(0,\mathfrak{d}^2)$.
        \item if $|\mathfrak{b}|>1$, there are two stationary measures to \eqref{eq fancy drift new formulation} given by $\cN(0,\mathfrak{d}^2)$ and $\cN(\frac{2\mathfrak{b}\mathfrak{d}}{\mathfrak{b}^2-1},\mathfrak{d}^2)$.
    \end{enumerate}
\end{Lemma}

Using Lemma \ref{le simple majo}, we have that $\vartheta^\star = 1-\mathfrak{b}^2$ for $|\mathfrak{b}| < 1$. For this range of parameter, Theorem \ref{th conv emp sid to emp mkv no rate} below, guarantees the convergence of the empirical measure of the SID to the unique stationary measure. Moreover, Theorem \ref{th main results intro L2 rate} shows that the algorithm converges with a rate and this is confirmed by our simulations.

\begin{figure}[h]
    \begin{center}
        \includegraphics[width=0.8\textwidth]{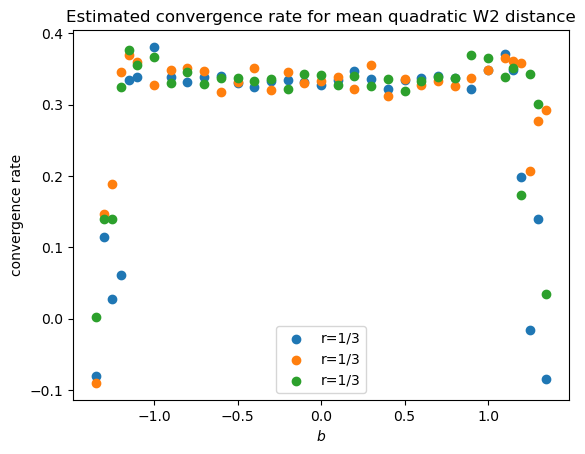}
    \end{center}
    \caption{Estimated convergence rate for the model \eqref{eq fancy drift new formulation} for a fixed algorithm step rate $r=1/3$. Three different runs are plotted, $N=100000$ and $M=500$ in \eqref{eq de error}. $\bar{\cX}_0$ follows a uniform distribution on $[0,1].$}  
    \label{fig conv rate fancydrift}
\end{figure}

On Figure \ref{fig conv rate fancydrift}, we plot the empirical convergence rate for different values of the parameter $\mathfrak{b}$ ($\mathfrak{d}=1$ for the simulations). The algorithm input step rate is fixed at $r=\frac13$. We observe that the convergence rate is better than the one predicted by our theoretical results and is in fact optimum (around $1/3$) for the value of $|\mathfrak{b}|\in (-1,1)$. A surprising fact is that the algorithm converges to the stationary measure $\cN(0,1)$ for value $|\mathfrak{b}|\in [1,1.15]$ with a rate $1/3$. Afterwards, the rate of convergence estimated is more erratic and possibly null. \textcolor{black}{On Figure \ref{fig second moment fancydrift}, we observe that the empirical second moment of $(\bar{\cX}_{\Gamma_n})_{0\le n \le N}$  seems to diverge for $\mathfrak{b}\in\set{-1.4,1.35,1.4}$ and to converge for $\mathfrak{b}\in\set{0.1,0.7,1.2}$.} One can show that, for values $\mathfrak{b}\ge 1$, the law of the MKV equation started with a law close to the stationary measure $\cN(0,\mathfrak{d}^2)$ will converge back to the stationary measure. This might explain why the algorithm work for some values $\mathfrak{b}\ge 1$. Note however that the behavior of the associated SID or the algorithm itself is not known for this parameter range and this study is left for further research. 

\begin{figure}[h]
    \begin{center}
        \includegraphics[width=0.8\textwidth]{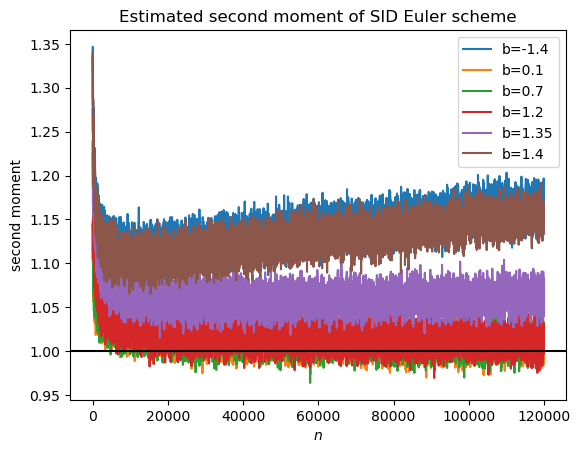}
    \end{center}
    \caption{Estimated second moment of the scheme for the model \eqref{eq fancy drift new formulation} for a fixed algorithm step rate $r=1/3$ (M=20000). $\bar{\cX}_0$ follows a uniform distribution on $[0,1].$ }  
    \label{fig second moment fancydrift}
\end{figure}

\paragraph{A model with mean field interaction in the volatility}
\begin{figure}[h]
    \begin{center}
    \includegraphics[width=0.8\textwidth]{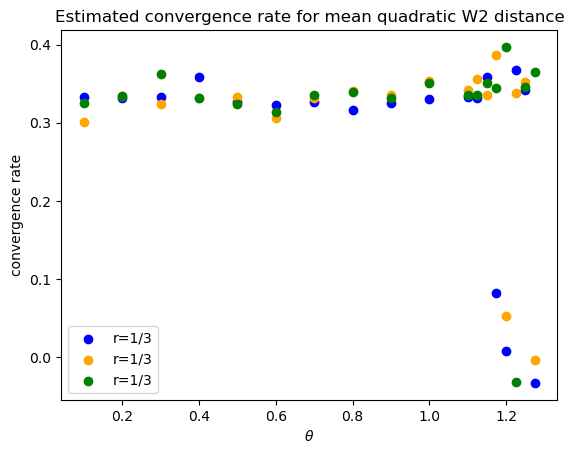}
    \end{center}
    \caption{Estimated convergence rate for the model \eqref{eq duetal light} for a fixed algorithm step rate $r=1/3$, $N=100000$ and $M=500$ in \eqref{eq de error}. $\bar{\cX}_0$ follows a uniform distribution on $[0,1].$}  
    \label{fig conv rate duetal}
\end{figure}
Lastly, we revisit a toy model introduced in \cite{du2023empirical}, given by, for $\theta>0$,
\begin{align}\label{eq duetal light}
    \ud X_t = - X_t \ud t + \sqrt{2}\theta\Big(1-\sqrt{\esp{|X_t|^2}}\Big)\ud W_t.
\end{align}
We write the stationary version as
\begin{align}
    \ud X^\star_t = - X^\star_t  \ud t + \sqrt{2}\mathfrak{d}^\star\ud W_t
\end{align}   
setting $\mathfrak{d}^\star:=   \| X^\star_t\|_2 $. From \eqref{eq duetal light}, we also have that $\mathfrak{d}^\star$ is solution to the following equation
$(\mathfrak{d}^\star)^2 = \theta^2\Big(1- \mathfrak{d}^\star\Big)^2$.
We easily deduce that the unique solution is given by $\mathfrak{d}^\star =\frac{\theta}{1+\theta}$. 

On Figure \ref{fig conv rate duetal}, we plot three runs of the algorithm, all using an input step rate of $r=\frac13$.
Using Lemma \ref{le simple majo}, we obtain that $\HYP{C}_{\alpha,\beta}$ is satisfied here when $\theta < 1$.  We observe  that the method does converge for this parameter range with the optimal rate $\mathfrak{r}=1/3$. The graph shows moreover that the convergence holds true up to a certain point (near $\theta=1.2$ ) where the empirical convergence rate becomes erratic and the algorithm unstable. Indeed, on Figure \ref{fig second moment mfvol}, we plot the empirical second moment of the Euler scheme $(\bar{\cX}_{\Gamma_n})$: One observes that for values $\theta \in \set{1.1,1.2}$, the estimated value are reasonnably close to the expected ones (the second moment of the stationary distribution). For $\theta=1.3$, the estimated values are highly oscillating. This is due to some paths of the algorithm taking extreme values.  Also, one can verify that started from a law with second moment close enough to $\mathfrak{d}^\star$, the law of $(X_t)$ will converge to the stationary limit. This might explain that the algorithm is still working for value of $\theta$ greater than one. Note however that the behavior of the associated SID or the algorithm itself is not known for this parameter range and this study is left for further research. 

\begin{figure}[h]
    \begin{center}
        \includegraphics[width=0.8\textwidth]{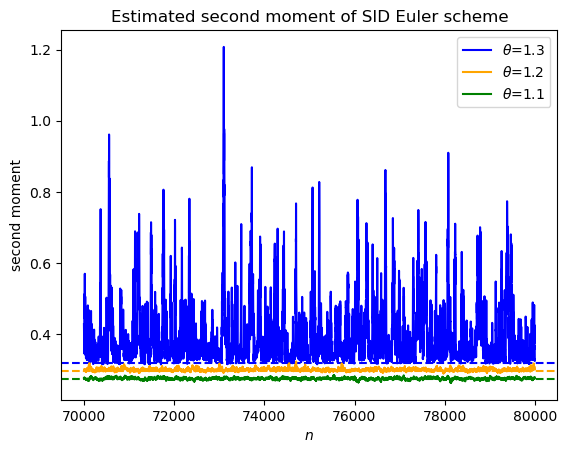}
    \end{center}
    \caption{Estimated second moment of the scheme for the model \eqref{eq duetal light} for a fixed algorithm step rate $r=1/3$. Dashlines give the value for the associated stationary distribution. $M=20000$ and $\bar{\cX}_0$ follows a uniform distribution on $[0,1].$  }  
    \label{fig second moment mfvol}
\end{figure}

\vspace{2mm}
The last two experiment tends to confirm that the "easy" choice of the input step rate $r=\frac13$ is a reasonnable one and could be used on more complex model.

}

\subsection{Notations}
$\bullet$ $(u\,|\,v) = \sum_{1\le i \le d}u^i v^i$ denotes the canonical inner product on $\R^d$ \textcolor{black}{of $u= (u^1,\ldots, u^d)$, $v=(v^1, \ldots,v^d)\!\in \R^d$, $|\cdot|$ the associated Euclidean norm, and $(A\,|\, B)_{_F}= \sum_{ij} a_{ij}b_{ij} = \mathrm{Tr}[A^\top B]$ denotes the ``Fr\" obenius'' inner product between  matrices $A=[a_j]$ and $B=[b_{ij}] \! \in  \cM_{d,q}(\R)$ (set of matrices with $d$ rows and $q$ columns), $\| \cdot\|_{F}$ the associated Fr\" obenius norm.}

\vspace{4pt}
\noindent $\bullet$ For $p\ge 1$, let 
\begin{align}\label{de eq a lyapunov function}
    \R^d \ni y\mapsto \cV_p(y)=|y|^{2p} \in \R. 
\end{align} 

We know that, for all $y,h \in \R^d$, 
\begin{align*}
    \cV_{p}(y+h)=\cV_p(y) + D\cV_p(y)(h) + \frac12 D^2\cV_p(y)(h,h) + o_y(|h|^2)
\end{align*}
where 
\begin{align*}
    &D\cV_p(y)(h) = 2p|y|^{2(p-1)}(y\,|\,h) \text{ and } 
    D^2\cV_p(y)(h,h) = \left( \partial^2_{xx}\cV_p(y)\,|\,hh^\top \right)_F , 
    \\ &\text{ with } 
    \partial^2_{xx}\cV_p(y)= 2p|y|^{2(p-1)} \left(\frac{2(p-1)}{|y|^2}yy^\top + I_d \right).
\end{align*}
We also observe that 
\begin{align}\label{eq useful majo hess}
    |D^2\cV_p(y)(h,h)| \le 2p(2p-1)|y|^{2(p-1)}|h|^2.
\end{align}

$\bullet$ In the following, we denote by $C$ a generic positive constant that may change from line to line. It depends on $d,L,\astr$ and $\bstr$. If it depends on an extra parameter $\theta$, it will be denoted $C_\theta$. This convention is also in force when using the Landau notation. 

\medskip
Moreover, we consider measurable positive functions $\mathfrak{C}:(\R_+)^\mathfrak{q} \rightarrow \R_+^*$ ($\mathfrak{q}$ is a positive integer) which satisfies 
\begin{align}\label{eq lim constant} 
    \lim_{\theta \rightarrow \partial (\R_+)^\mathfrak{q}} \mathfrak{C}(\theta) = +\infty.
\end{align}
This will allow us to clarify the dependence on critical parameters of the various constants appearing in the statement of our results below (on top of the dependence upon the fixed parameters $\astr$, $\bstr$ and $L$).
\\
We also extend this notation to Landau notation. Namely, for a parameter $\theta \in (\R_+)^\mathfrak{q}$, we shall write $f(x)=O_\theta(g(x))$ ($f,g$ are measurable positive functions) if 
\[ \lim_{x \rightarrow +\infty} \frac{f(x)}{g(x)} \le \mathfrak{C}(\theta).
    \]
Similarly, $f(x)=o_\theta(g(x))$ if $f(x)=\mathfrak{C}(\theta) o_1(g(x))$ when $x \rightarrow +\infty$.
(The variable $x$ can be a real or an integer.)

\section{Stationary solution to McKean-Vlasov SDEs}

\label{se stationnary solution MKVSDE}
\subsection{Existence and uniqueness of the invariant distribution}
In this section, we first recall that there exists a unique invariant distribution for the MKV dynamics \eqref{eq de MKVSDE} under $\hypl2$. This results is well known, see e.g. \cite{wang2018distribution}. We will follow a simple fixed point approach,  exposed here for sake of completeness. The proofs   given  below should be seen as an introduction to the proof for self-interacting diffusions and the associated numerical scheme.

\noindent Let us thus consider, for a given $\mu \in \cP_2(\R^d)$, a version of  
\eqref{eq de MKVSDE} where the dependence upon the law is frozen, namely
\begin{align}\label{eq SDE for fixed point}
X^{\mu,\xi}_t = \xi + \int_0^t b(X^{\mu,\xi}_s,\mu)\ud s + \int_0^t\sigma(X^{\mu,\xi}_s,\mu)\ud W_s,
\end{align}
for $\xi \in \mathscr{L}^2(\cF_0,\P)$ (note that we do not impose $[\xi]=\mu$).
We now introduce 
$$ \cP_2(\R^d) \ni \mu \mapsto \Pi(\mu) \subset \cP_2(\R^d)$$ 
the set of invariant measure of $\cP_2(\R^d)$ of \eqref{eq SDE for fixed point}.

\noindent We are thus looking for probability measure satisfying $\mu^\star \in \Pi(\mu^\star)$. Actually, we focus here on the case where $\Pi$ is single-valued as we impose the strong confluence assumption $\HYP{C}_{\a,\b}$ on $(b,\sigma)$. 

We first report some integrability results of \eqref{eq SDE for fixed point}. Note that since $\mu$ is fixed in the coefficients $b$ and $\sigma$, $X^{\mu,\xi}$ is a ``classical'' SDE. 

\color{black}
\begin{Lemma}\label{le unif int}
  Assume $\HYP{MV}_{\bar p, \bar K, \bar\a,\bar\b}$ holds for some $\bar p \ge 1$ and $\bar K \ge 0$, $\bar \b \ge 0$ and $\bar \a > 0$.
   Then, for all  $p \in [1,\bar{p}]$,  for every $\mu \in \cP_2(\R^d)$ and $\xi \in \mathscr{L}^{2p}(\cF_0,\P)$, 
    \begin{align*}
        \sup_{t \ge 0}\esp{|X^{\mu,\xi}_t|^{2p}} \le \bar{C} \left(1 + \esp{|\xi|^{2p} }+ \cW_2^{2p}(\mu,\delta_0)\right),
    \end{align*} 
    where $\bar{C}$ is a positive constant that depends on $\bar K,\bar \a,\bar \b$ and $p$.
\end{Lemma}

\color{black}
\noindent {\em Proof.} For the reader's convenience, we simply write $X$ for $X^{\mu,\xi}$ in this proof. 
   %
      %

\noindent    Setting $\la = p \bar \a$ and applying Ito's formula to $t \mapsto e^{\la t}\cV_p(X_t)$, recall \eqref{de eq a lyapunov function}, we obtain 
    \begin{align}
        &e^{\la t}|X_t|^{2p}  = |\xi|^{2p} 
        + 2p\int_0^t e^{\la s}|X_s|^{2(p-1)}\left( X_s \,|\, \s(X_s,\mu)\right) \ud W_s
        \label{eq ito to Vp}
        \\
        & + \int_0^t\!\!e^{\la s} \set{\la |X_s|^{2p}\!+\!2p|X_s|^{2(p-1)}\left( X_s\,|\, b(X_s,\mu)\right) 
         \!+\! \frac12 \left( \partial_{xx}^2 \cV_p(X_s) \,|\, \s(X_s,\mu)\s(X_s,\mu)^\top\right)_F}\ud s. \nonumber
    \end{align}
    Then it follows successively  from~\eqref{eq useful majo hess} and $\HYP{MV}_{{p},\bar{K},\bar{\alpha},\bar{\beta}}$ that, for every $s\ge 0$,  
    \begin{align*}
   \la   |X_s|^{2p}+2p |X_s|^{2(p-1)} \left( X_s\,|\, b(X_s,\mu)\right) 
         & + \frac12 \left( \partial_{xx}^2 \cV_p(X_s) \,|\, \s(X_s,\mu)\s(X_s,\mu)^\top \right)_F \\
        & \le p e^{\la s} |X_s|^{2(p-1)} \big(\tfrac{\la}{p} |X_s|^{2}+2\left( X_s| b(X_s,\mu)\right)\\
          &\qquad+  (2p-1)\|\s(X_s,\mu)\|_F^2 \big)\\
          &\le p  e^{\la s} |X_s|^{2(p-1)}\big( \bar K+ \bar \b \cW_2^2(\mu,\delta_0) \big).
      \end{align*}
      Hence
      \begin{align} 
        e^{\la t}|X_t|^{2p} +(p-1) \int_0^t   e^{\la s} |X_s|^{2p} \ud s
        &
        \le |\xi|^{2p} + 2p\int_0^t e^{\la s}|X_s|^{2(p-1)}\left( X_s \,|\, \s(X_s,\mu)\right) \ud W_s \nonumber
        \\
        &+p \int_0^t e^{\la s} |X_s|^{2(p-1)}\big( \bar K+ \bar \b \cW_2^2(\mu,\delta_0) \big)\ud s. \label{eq proof for next prop}
      \end{align}
    Let $(\tau_n)_{n\ge 1}$ be a sequence of stopping times increasing to $+\infty$ such that the stochastic integral in~\eqref{eq ito to Vp} is a martingale and $X^{\tau_n}$ is bounded by $n$. Taking expectation in the above inequality taken at the stopping time $\tau_n\wedge t$ yields
 
     \begin{align*}
       \E\,  e^{\la (t\wedge \tau_n)}|X_{t\wedge \tau_n}|^{2p} &+(p-1) \E\, \int_0^{t\wedge \tau_n}   e^{\la s} |X_s|^{2p} \ud s\\
        &\le \E\,|\xi|^{2p} +p\E\, \int_0^{t\wedge \tau_n} e^{\la s} |X_s|^{2(p-1)}\big( \bar K+ \bar\beta \cW_2^2(\mu,\delta_0) \big)\ud s.
      \end{align*}
 Letting $n\to+\infty$, it follow from Beppo Levi's monotone convergence theorem that 
   \begin{equation}\label{eq unif int mkv temp}
      e^{\la t}\E\, |X_t|^{2p} +(p-1) \int_0^t   e^{\la s} \E\, |X_s|^{2p}   \ud s\le \E\, |\xi|^{2p}+p  \int_0^t e^{\la s} |X_s|^{2(p-1)}\big( \bar K+ \bar \beta \cW_2^2(\mu,\delta_0) \big)\ud s. 
    \end{equation}
    Young's inequality then yields 
    \begin{align*}
        |X_s|^{2(p-1)}\left(\bar{K}+ \beta \cW_2^2(\mu,\delta_0)\right)
        \le \frac{p-1}{p}|X_s|^{2p} + \frac1p\left(\bar{K}+ \bar \beta \cW_2^2(\mu,\delta_0)\right)^{p}.
    \end{align*}
    Inserting the previous inequality into \eqref{eq unif int mkv temp} completes the proof since we obtain
    \[
     \hskip 3,25cm   \esp{|X_t|^{2p} } \le  \esp{|\xi|^{2p}} 
         + \frac1\la\left(\bar{K}+ \bar \beta \cW_2^2(\mu,\delta_0)\right)^{p}, \text{ for all } t \ge 0.\hskip 1,25 cm \Box
    \]
\color{black}

\begin{Proposition}\label{pr one invariant measure}
    Let 
    $\hypl2$ 
    hold.   
    Then, for some $p^\star:=p^\star(\astr,\bstr)> 1$,
    \begin{align}
    \cP_2(\R^d) \ni \mu \mapsto \Pi(\mu) \in \cP_{2p^\star}(\R^d) \label{eq the app Pi}
    \end{align}
    is well defined and satisfies, for all $(\mu,\nu) \in \cP_2(\R^d)^2$,
\begin{align}\label{eq:contraction}
\cW_2(\Pi(\mu),\Pi(\nu)) \le \sqrt{\frac{\bstr}{\astr}}\, \cW_2(\mu,\nu).
\end{align}
\end{Proposition}
{\color{black}
\proof Under $\hypl2$,  $\HYP{MV}_{\bar p, \bar K, \bar\a,\bar\b}$ holds for some $\bar p > 1$ and $\bar K \ge 0$, $\bar \b \ge 0$ and $\bar \a > 0$, see Lemma \ref{le structural assumptions3}. \\
1.a For $x_0 \in \R^d$, let $(P_t(x_0,\ud y))_{t \ge 0}$ be the transition measure of the semi-group associated to $X^{x_0,\mu}$ and set
\begin{align}
    \nu_t(\ud y) := \frac1t \int_0^tP_s(x_0,\ud y)\ud s.
\end{align}
From Lemma \ref{le unif int}, we know that the family $(\nu_t)$ satisfies
\[
\sup_{t\ge 0}\nu_t(|\cdot|^{2\bar p}) \le C_{\a,\b,\bar p}\left(1 +  |x_0|^{2\bar p}  + \cW_2^{2\bar p}(\mu,\delta_0)\right).
\]
Hence, if $p^{\star}\!\in (1, \bar p)$,  the family $(\nu_t)$ is $\cW_{2p^{\star}}$-relatively compact. Thus there exists a subsequence $(\nu_{t_n})_{n \ge 0}$,  $t_n \rightarrow +\infty$ and $\bar{\nu} \!\in \cP_{2p^\star }(\R^d)$ such that $\cW_{2p^{\star}}(\nu_{t_n},\bar{\nu})\to 0$ as $n\to +\infty$. 
One then computes that, for any  $f$ bounded and continuous,
\begin{align}
    \nu_{t_n}f - \nu_{t_n}P_tf =O\big(\frac{2|f|_\infty}{t_n}\big).
\end{align}
As $ \nu_{t_n}f \to \bar \nu f$ and $\nu_{t_n}P_tf\to \bar \nu P_tf$, $\bar \nu f = \bar \nu P_t f$ which concludes the proof for this step.

\medskip
\noindent 1.b  Let us show that any invariant distribution of~\eqref{eq SDE for fixed point} lies (at least) in $\cP_{2}(\R^d)$.  From \eqref{eq proof for next prop}  with $p=1$, we obtain for every $\Lambda > 0$, (using the notations of the proof of Lemma \ref{le unif int})
\begin{align*}
    |X^{x,\mu}_{t\wedge \tau_n}|^2\wedge \Lambda \le \big( |x|^2 e^{- \bar \a (t\wedge \tau_n)}\big)&\wedge \Lambda +\frac{\bar \b \cW^2_2(\mu,\delta_0) + \bar K }{\bar \a}\wedge \Lambda 
    \\
    &+2 \Big( e^{-\bar \a( t\wedge \tau_n)}\int_0^{t\wedge \tau_n}e^{-\bar \a s}(X^{x,\mu}_s\,|\, \s(X^x_s))\ud W_s)\Big)    
\end{align*}

where the stopping times still localize the stochastic integral into a martingale. Taking expectations and then letting $n\to+\infty$, we obtain
\[
\E\, |X_t^{x,\mu}|^2\wedge \Lambda \le (|x|^2e^{-\bar \a t})\wedge \Lambda +\frac{\bar  \b \cW^2_2(\mu,\delta_0) + \bar K }{\bar \a}\wedge \Lambda 
\]
Integrating w.r.t.  to an invariant distribution of~\eqref{eq SDE for fixed point} denoted $\nu$ yields for every $t\ge 0$, 
\[
\nu(|\cdot|^2\wedge \Lambda) = \int_{\R^d} \esp{ |X_t^{x,\mu}|^2\wedge \Lambda}\nu(dx) \le  \int_{\R^d} (|x|^2e^{-\bar \a t})\wedge \Lambda)\nu(dx) +\frac{\bar  \b \cW^2_2(\mu,\delta_0) + \bar K }{\bar \a}\wedge \Lambda  .
\]
Letting successively $t$ and $\Lambda$ go to $+\infty$ yields $
\nu(|\cdot|^2) \le \frac{\bar  \b \cW^2_2(\mu,\delta_0) + \bar K }{\bar \a}$.\\
1.c. Assume that there are two invariant probability measures $\nu^j$, $j \in \set{1,2}$ and denote $X^j=X^{\mu,\xi^j}$ with $\xi^j \sim \nu^j$, $j \in \set{1,2}$, the associated diffusion. Then one computes (still using a localization procedure of the stochastic integral)
\begin{align*}
    &e^{{\astr } t}\esp{|X^1_t-X^2_t|^2} = \esp{|X^1_0-X^2_0|^2}
    + {\astr} \int_0^t e^{{\astr} s}\esp{|X^1_s-X^2_s|^2} \ud s
    \\
    &+\int_0^t e^{{\astr} s}\esp{
        2 \left( b(X^1_s,\mu) - b(X^2_s,\mu) | X^1_s-X^2_s \right) + \| \sigma(X^1_s,\mu) - \sigma(X^2_s,\mu)\|_{_F}^2 }\ud s.
\end{align*}
From $\HYP{C}_{\astr,\bstr}$, we obtain 
\begin{align*}
    \cW_2(\nu^1,\nu^2) \le  e^{-\astr t}\esp{|X^1_t-X^2_t|^2}\le e^{-\astr t} \esp{|X^1_0-X^2_0|^2}
\end{align*}
which proves uniqueness of the invariant probability measure.
\\
2. For $i \in \set{1,2}$, consider $\mu^i \in \cP_2(\R^d)$ and $\xi^i \sim \Pi(\mu^i)$. Denoting $X^i = X^{\mu^i,\xi^i}$, we observe then that $[X^i_t] = \Pi(\mu^i)$, for all $t \ge 0$. Applying It\^o's formula to $(e^{\astr  t}|X^1_t - X^2_t|^2)_{t \ge 0}$, we compute
\begin{align*}
&e^{{\astr} t}|X^1_t - X^2_t|^2 = |X^1_0 - X^2_0|^2 + 2\int_0^t e^{{\astr} s}(X^1_s - X^2_s)\left( \sigma(X^1_s,\mu^1) - \sigma(X^2_s,\mu^2) \right) \ud W_s
\\
&\!+ \!\! \int_0^t\!\! e^{{\astr} s}\Big(2 \big( b(X^1_s,\mu^1) - b(X^2_s,\mu^2) | X^1_s - X^2_s\big) \!+\! 
\| \sigma(X^1_s,\mu^1) - \sigma(X^2_s,\mu^2)\|_{_F}^2 \!+\! \astr|X^1_s - X^2_s|^2\Big) \ud s\\
 & \qquad\qquad\qquad   \le |X^1_0 - X^2_0|^2+ 2\int_0^t e^{{\astr} s}(X^1_s - X^2_s)\left( \sigma(X^1_s,\mu^1) - \sigma(X^2_s,\mu^2) \right) \ud W_s\\
&\qquad \qquad \qquad \qquad  \qquad  \qquad \; +\bstr \cW^2_2(\mu^1,\mu^2) \int_0^t e^{\astr s}\ud s 
 \end{align*}
  using $\HYP{C}_{\astr,\bstr}$.
Localizing the stochastic integral by stopping times $\tau_n\uparrow +\infty$ taking expectation , and letting $n\to+ \infty$ owing to Beppo Levi theorem, yields on both sides, 
\begin{align*}
e^{{\astr}t}\esp{|X^1_t - X^2_t|^2} \le \esp{|X^1_0 - X^2_0|^2}
+\bstr \cW^2_2(\mu^1,\mu^2) \int_0^t e^{\astr s}\ud s .
\end{align*}
We then deduce
\begin{align}\label{eq main estimate}
\cW_2^2(\Pi(\mu^1),\Pi(\mu^2)) \le e^{-\astr t}\cW_2^2(\Pi(\mu^1),\Pi(\mu^2))
+ \frac{\bstr}{\astr }\left( 1 - e^{-\astr t}\right)\cW^2_2(\mu^1,\mu^2),
\end{align}
which proves the property~\eqref{eq:contraction}  by letting $t\o +\infty$. Finally, as $\Pi(\mu)$ is an invariant distribution of~\eqref{eq SDE for fixed point}  associated to $\mu$, it follows from~1.b that   $\Pi(\mu)\!\in \cP_{2p^{\star}}(\R^d)$.
\eproof
}

\bigskip


\noindent We also straightforwardly make the following observation.

\begin{Corollary} \label{co main res invariant distrib mkv}
\begin{enumerate}[(a)]
    \item Assume that $\hypl2$ is in force. Then, there exists  a unique probability measure $\nu^\star \in \cP_2(\R^d)$ such that $\nu^\star = \Pi(\nu^\star)$. Denoting $X^\star = X^{\nu^\star,\xi}$, for some $\xi\sim\nu^\star$ with $\xi \in \mathscr{L}^{2{p^\star}}(\cF_0,\P)$ (Note that the notation makes no reference to the starting random variable $\xi$) we have then 
     \begin{align}\label{eq integrability X star}
        \sup_{t \ge 0} \esp{|X^\star_t|^{2p^\star} } = \int_{\R^d}|y|^{2p^\star} \nu^\star(\ud y) < +\infty.
    \end{align}
    \item Moreover, if $\HYP{MV}_{\bar p, \bar K, \bar \a, \bar \beta}$ holds, then one has $p^\star \ge \bar p$.
\end{enumerate}
\end{Corollary}



\subsection{$A.s.$ and $L^2$ $\cW_2$-convergence of the empirical measure}
We conclude this section by stating convergence results for the empirical measure of the stationary McKean-Vlasov SDEs $X^\star$ toward its invariant distribution $\nu^\star$.

\medskip
\noindent The following proposition is classical in the framework of regular diffusions.
\begin{Proposition}\label{pr empirical measure no rate}
    Under $\hypl2$, it holds that     
    \begin{align} \label{eq le almost sure conv mvk}
        \lim_{t \rightarrow +\infty} \cW_2(\nu^{X^\star}_t,\nu^\star) = 0 \quad \P\mbox{-}a.s.
    \end{align}or, equivalently, with obvious notation,
    \[
    \nu^{\star}(dx)\mbox{-}a.s. \quad \lim_{t\rightarrow +\infty}\cW_2(\nu^{x,\nu^{\star}}_t,\nu^\star) = 0 \quad \P_x\mbox{-}a.s.
    \]
    and 
    \begin{align} \label{eq le L2 conv mvk}
        \lim_{t \rightarrow +\infty} \esp{\cW_2^2(\nu^{X^\star}_t,\nu^\star)} = 0.
    \end{align}
\end{Proposition}

{
    \color{black}
    \proof  1.  As $\nu^{\star}$ is (also) the  unique invariant measure  of the diffusion~SDE~\eqref{eq SDE for fixed point} (with $\mu=\nstr$)  which satisfies $\HYP{MV}_{1, \bar K, \bar \a,\bar \b}$ for some $\bar K \ge 0$, $\bar \b \ge 0$ and $\bar \a > 0$. Hence it is trivially an extremal point of the set of its invariant distribution(s). Then it is an ergodic process. Which implies~\eqref{eq le almost sure conv mvk} (see e.g.~\cite{pages2001quelques} among others for more details).\\ 
    2. Let us observe that 
    \begin{align*}
        \cW_2^2(\nu^{X^\star}_t,\nu^\star) \le 2 \cW_2^2(\nu^{X^\star}_t,\delta_0) +
        2 \cW_2^2(\delta_0,\nu^\star),
    \end{align*}
    and
    \begin{align*}
        \cW_2^{2\pstr}(\nu^{X^\star}_t,\delta_0) &\le \left(\frac1t\int_0^t |X^\star_s|^2 \ud s\right)^{ \pstr }
        \le \frac1t\int_0^t |X^\star_s|^{2 \pstr} \ud s
    \end{align*}
    where we used the definition of the $\cW_2$-distance and Jensen inequality. Taking expectation on both sides of the previous inequality and invoking Corollary \ref{co main res invariant distrib mkv}, we obtain 
    \begin{align}
        \esp{\cW_2^{ 2 \pstr}(\nu^{X^\star}_t,\delta_0)} \le C\;.
    \end{align}
    By de la Vall\'ee Poussin Theorem, the family $\left(\cW_2^{2}(\nu^{X^\star}_t,\nu^\star)\right)_{t \ge 0}$ is thus uniformly integrable, which, combined with \eqref{eq le almost sure conv mvk}, yields \eqref{eq le L2 conv mvk}.
    \eproof
}

\bigskip

We can improve the previous convergence results by providing (sub-optimal) rate of convergence in various situations.

\begin{Corollary} \label{co L2 conv rate emp meas to inv meas} Assume that $\hypl2$ is in force. Then
    \begin{align*} 
        \esp{ {\cal W}^2_2(\nu^{X^\star}_t, \nu^\star)}^\frac12 = O\left( t^{-\frac{2 \pstr-1}{2(2(d+2)+(2 \pstr-1)(d+3))}} \right),
    \end{align*}
    where $p^\star$ is defined in Proposition \ref{pr one invariant measure}.
    \\
    If, furthermore $\HYP{MV}_{\bar p,\bar K,\bar \a,\bar \b}$ with $\bar K \ge 0$, $\bar \b \ge 0$ and $\bar \alpha > 0$ is in force for any $\bar p \ge 1$ (e.g. if $\HYP{\sigma}$ holds according to Lemma \ref{le structural assumptions3}) then  
    \begin{align}
        \esp{ {\cal W}^2_2(\nu^{X^\star}_t, \nu^\star)}^\frac12 = o_\eta\left(t^{-\frac1{2(d+3)}+\eta}\right), \; \text{ for any small $\eta>0$.}
    \end{align}
Finally, if $\HYP{\sigma_\infty}$ holds  then 
    \begin{align*}
        \big(\E\, \cW_{2}^2(\nu^{X^\star}_t,\nu^\star)\big)^\frac12 = O \left(t^{-\frac1{2(d+3)}} 
        (\log t)^{\frac{d+2}{2(d+3)}}\right).
    \end{align*}
\end{Corollary}

\color{black}
\proof  We use the result of Theorem 1.1 in \cite{chapag24}. 
We can do so because assuming $\hypl2$ implies that $x \mapsto \Sigma(x):=\sigma(x,\mu^\star)$ and $x \mapsto B(x):= b(x,\mu^\star)$ do satisfy a strong confluence assumption as required in \cite{chapag24}. 
For the second statement, it follows from  Lemma \ref{le unif int} that 
$\int_{\R^d} |y|^{2p} \nstr(\ud y) < +\infty$ for all $p \ge 1$ and the result is then a direct application of the first upper bound.
For the third statement, we observe that $\HYP{\sigma_\infty}$ combined with $\HYP{L}$ implies that $x \mapsto \Sigma(x)=\sigma(x,\nstr)$ is bounded and then the stationary distribution $\nu^\star$ has exponential moments, and then we can invoke directly Theorem 1.2 in \cite{chapag24}.
\eproof 
\color{black}

\bigskip
\noindent Let us now turn to the case of the almost sure convergence.

\begin{Corollary} \label{co as conv rate emp meas to inv meas}
    Assume that $\hypl2$ is in force and that $x \mapsto \sigma(x,\nu^\star)$ is uniformly elliptic. Then, 
    \begin{align*}
        \cW_2(\nu^{X^\star}_t,\nu^\star) = o_\eta\big(t^{-\hat{\zeta}_{\pstr}} \log(t)^{\frac12+\eta}\big)\,, \text{ for any small }\, \eta >0,
    \end{align*}
    where $\hat{\zeta}_{\pstr} = \frac{(2 \pstr - 1)^2}{2(2 \pstr -1) \{2(d+2) + (d+3)(d + 2 \pstr - 1) \}}$ and $\pstr$ is defined in Proposition \ref{pr one invariant measure}.  \\
    If, furthermore $\HYP{MV}_{\bar p,\bar K,\bar \a,\bar \b}$ with $\bar K \ge 0$, $\bar \b \ge 0$ and $\bar \alpha > 0$ holds for any $\bar p \ge 1$ (e.g. if $\HYP{\sigma}$ holds according to Lemma \ref{le structural assumptions3}) then  
    \begin{align}
         {\cal W}_2(\nu^{X^\star}_t, \nu^\star) = o_\eta\left(t^{-\frac1{2(d+3)}+\eta}\right)\,, \text{ for any small }\, \eta >0.
    \end{align}
Finally, if $\HYP{\sigma_\infty}$ holds  then 
    \begin{align*}
       \cW_{2}^2(\nu^{X^\star}_t,\nu^\star)\big) = o_\eta \left(t^{-\frac1{2(d+3)}} 
        (\log t)^{\frac12+\eta}\right) \,, \text{ for any small }\, \eta >0.
    \end{align*}
\end{Corollary}   

\color{black}
\noindent The proof of the previous result is a straightforward adaptation of the proof of Corollary \ref{co L2 conv rate emp meas to inv meas} based again on Theorem 1.1 and Theorem 1.2 in the companion paper \cite{chapag24}: It is left to the reader.  
\color{black}

\section{A class of self-interacting diffusions}

\label{se SID}
Our aim is to compute the invariant distribution of \eqref{eq de MKVSDE} i.e. the fixed point of the application $\Pi$, given in \eqref{eq the app Pi}. However, as the function $\Pi$ is not directly computable, in the following, we rely on the following fact 
\begin{align*}
\Pi(\mu) = \lim_{t \rightarrow +\infty} \nu^\mu_t:= \frac1t\int_0^t \delta_{X^{\mu,x}_s} \ud s , 
\end{align*}
see e.g. the first step of the proof of Proposition \ref{pr one invariant measure}.

\noindent One could then use a Picard iteration based on successive approximation of $\nu^\mu_T$ for given $T>0$ or as a limit case of the previous procedure introduce the self-interacting diffusion:
\begin{align}\label{eq de sid}
\cX^{\xi}_t = \xi + \int_0^t b(\cX^\xi_s, \nu^\cX_s)\ud s +  \int_0^t \sigma(\cX^\xi_s, \nu^\cX_s)\ud W_s\, \text{ with } \nu^\cX_s := \frac1s\int_0^s \delta_{\cX^\xi_u} \ud u,
\end{align}
and $\xi \in \mathscr{L}^2(\cF_0,\P)$. As announced in the introduction, this is the road we follow here. 

\smallskip \noindent
According to \cite[section 7]{du2023empirical}, there exists a unique solution to \eqref{eq de sid}, that satisfies, for all $t \ge 0$, 
\begin{align}
    \sup_{s\le t} \esp{|\cX^\xi_s|^2} < +\infty.
\end{align}

\noindent The main goal of this section is to prove that $\nu^\cX$ converges to $\nu^\star$ in $L^2$ and almost surely. This is the content of Theorem  \ref{th conv emp sid to emp mkv no rate}. Then, we improve the previous results in Theorem \ref{th conv emp sid to emp mkv with a rate} by proving convergence with a rate.

\subsection{Some key properties}

\subsubsection{Uniform boundedness}

\color{black}
\begin{Lemma}\label{le unif int sid}
    Assume that $\HYP{MV}_{\bar{p},\bar{K},\bar\alpha,\bar\beta}$ for $\bar p \ge 1$ and with $\bar \a > \bar \b\ge 0$ is in force.
     Then for all $p \in [1,\bar{p}]$, as soon as $\xi \in \mathscr{L}^{2p}(\cF_0,\P)$,
    \begin{align*}
        \sup_{t \ge 0}\esp{|\cX^{\xi}_t|^{2p}} \le \bar{C}\left(1 + \esp{|\xi|^{2p} }\right),
    \end{align*} 
    where $\bar{C}$ is a positive constant that depends on $\bar K,\bar \a,\bar \b$ and $p$.
\end{Lemma}
\color{black}

   \color{black}
    \proof 
Set $\la = p \bar \a$. It follows from It\^o's formula to $ e^{\la \cdot}\,\cV_p(\cX)$ (have in mind  \eqref{de eq a lyapunov function}) applied between $0$ and $t\wedge \tau$ where  $\tau= \tau_n$ is  a stopping time localizing the Brownian stochastic integral and such that  the stopped process $\cX^{\tau}$ is  bounded by $n$: 
    \normalsize
    \begin{align}
        e^{\la (t\wedge \tau)}&|\cX_{t \wedge \tau}|^{2p}  = |\xi|^{2p} 
        + 2p\int_0^{t\wedge \tau} e^{\la s}|\cX_s|^{2(p-1)}\left( \cX_s \,|\, \s(\cX_s,\nu_s)\right) \ud W_s
        \label{eq It\^o to Vp bis}
        \\
        & + \int_0^{t\wedge \tau} e^{\la s} \Big[\la |\cX_s|^{2p}+ 2p|\cX_s|^{2(p-1)}\left( \cX_s\,|\, b(\cX_s,\nu_s)\right) 
         + \frac12 \left( \partial_{xx}^2 \cV_p(\cX_s) \,|\, \s(\cX_s,\nu_s)\s(\cX_s,\nu_s)^\top\right)_F\Big]\ud s. \nonumber
    \end{align}

    \noindent Taking  expectation in \eqref{eq It\^o to Vp bis} and using \eqref{eq useful majo hess}, we get 
    \begin{align*}
        \E\, e^{\la (t\wedge \tau)} |\cX^{\tau}_t|^{2p}  &\le \E\,|\xi|^{2p} 
        + p \E\, \int_0^{t\wedge \tau} e^{\la s} |\cX_s|^{2(p-1)}
        \Big[\tfrac{\lambda}{p} |\cX|^2+
            2\left( \cX_s| b(\cX_s,\nu_s)\right) 
            +
            (2p-1)\|\s(\cX_s,\nu_s)\|_F^2   \Big]\ud s.
    \end{align*}
    Since $\HYP{MV}_{{p},\bar{K},\bar{\alpha},\bar{\beta}}$ holds, we deduce from the previous inequality,
    \begin{equation}    \label{eq unif int mkv temp bis}
       \E\, e^{\la (t\wedge \tau)} |\cX^{\tau}_t|^{2p}   \le \E\, |\xi|^{2p} 
       + p \E\,\int_0^{t\wedge \tau} e^{\la s} |\cX_s|^{2(p-1)}
        \Big[
             \bar{K} -\big(\bar\alpha-\tfrac{\la}{p}\big)|\cX_s|^{2} + \bar\beta \cW_2^2(\nu_s,\delta_0)\Big] \ud s. 
    \end{equation}
Since $ \cW_2^2(\nu_s,\delta_0) \le \frac1s\int_0^s |\cX_u|^2\ud u$, we get, using Young's inequality, 
\begin{align*}
    p|\cX_s|^{2(p-1)}\cW_2^2(\nu_s,\delta_0) &\le \frac{p}s\int_0^s |\cX_s|^{2(p-1)}|\cX_u|^2\ud u  \le (p-1)|\cX_s|^{2p} + \frac1s\int_0^s |\cX_u|^{2p}\ud u. 
\end{align*}
Using again Young's inequality, we also have 
\begin{align}
    p|\cX_s|^{2(p-1)}\bar{K} \le (p-1)|\cX_s|^{2p} + \bar{K}^p.
\end{align}
Inserting these two inequalities into \eqref{eq unif int mkv temp bis} with  in mind the above  the definition of $\la$, one has
\begin{align*}
     \E\, e^{\la (t\wedge \tau)} |\cX^{\tau}_t|^{2p}  &\le \E\, |\xi|^{2p} +
  \bar{K}^p \E\,  \frac{e^{\la (\tau \wedge t)}-1}{\la} + \bar\beta \E \int_0^{t\wedge \tau} \frac{e^{\lambda s}}{s}\int_0^s \E\, |\cX_u|^{2p}\ud u \ud s \\
  &\le \E\, |\xi|^{2p} +
 \bar{K}^p \frac{e^{\la t}-1}{\la}+ \bar \beta \int_0^{t} \frac{e^{\lambda s}}{s}\int_0^s\E\,  |\cX_u|^{2p}\mbox{\bf 1}_{\{u\le \tau\}}\ud u \ud s\\
 &\le \E\, |\xi|^{2p} + \frac{e^{\la t}-1}{\la}\Big(\bar{K}^p + \bar \beta \sup_{u\in [0,t]} 
\E\, |\cX_u|^{2p}\mbox{\bf 1}_{\{u\le \tau\}}\Big).
\end{align*}
Noting that $  \E\, e^{\la (t\wedge \tau)} |\cX^{\tau}_t|^{2p}  \ge e^{\la t}  \E\, |\cX_t|^{2p}\mbox{\bf 1}_{\{t\le \tau\}}$ finally yields (having in mind that $\bar \rho = \frac{\bar \beta}{\la})$, 
\[
 \E\, |\cX_t|^{2p}\mbox{\bf 1}_{\{t\le \tau\}} \le e^{-\la t} \E\, |\xi|^{2p} + \frac{\bar K^p}{\la}  +\bar \rho \sup_{s\in [0,t]} 
\E\, |\cX_s|^{2p}\mbox{\bf 1}_{\{s\le \tau\}}
\]
which in turn implies, since $\bar \rho \!\in (0,1)$, 
\[
\sup_{s\in [0,t]} \E\, |\cX_s|^{2p}\mbox{\bf 1}_{\{s\le \tau\}} \le \frac{1}{1-\bar \rho} \Big(e^{-\la t} \E\, |\xi|^{2p} + \frac{\bar K^p}{\la}\Big).
 \]
Letting $\tau =\tau_n\uparrow +\infty$  completes the proof  owing to Fatou's lemma.
    \eproof

\color{black}

\subsubsection{Stability}

Let  $(\etb,\ets)$ be a square integrable, continuous and adapted process.
We consider a continuous and adapted square integrable process $\tilde{\cX}$  with dynamics
\begin{align}\label{eq de tXi}
\tcX_t = \tcX_0 + 
    \int_0^t \left(b (\tcX_s,\nu^{\tcX}_s) + \etb_s \right) \ud s 
    + \int_0^t \left(\sigma(\tcX_s,\nu^{\tcX}_s) + \ets_s \right) \ud W_s .
\end{align}

\begin{Proposition}\label{pr stability for sid} 
    Assume that $\hypl2$ hold. We define, for $t\ge 0$,  
    \begin{align}\label{eq pr def delta}
        \eg_t :=  |\etb_t|^2 + |\ets_t|^2,
    \end{align}
    which is a continuous adapted and integrable process. \\
$(a)$ For any $\vartheta \in (0,\vtstr)$, the following holds 
    \begin{align}\label{eq pr stab L2 upper bound}
        \frac1t\int_0^t \esp{|\tcX_s -\cX_s|^2} \ud s 
        \le   \,t^{-\vartheta } \left( C\, \esp{|\tcX_0-\cX_0|^2} + \mathfrak{C}(\vtstr-\vartheta)\int_0^t s^{{\vartheta}-1} \mathbb{E}[\eg_s] \ud s \right).
    \end{align}
$(b)$ Moreover, set, for $t\ge0$,
    \begin{align}\label{eq de H}
        H_t := 2(\tcX_t -\cX_t|\sigma(\tcX_t,\nu^{\tcX}_t) - \sigma(\cX_t,\nu^{\cX}_t)+ \ets_t)   
    \end{align}
    and assume that,
    \begin{align}\label{eq ass H}
        \sup_{t \ge 0} \E\,|H_t|^{2} <+\infty.
    \end{align}
    Then the following holds, for any $\vartheta \in (0,\vtstr)$,  $\eta > 0$,
    \begin{align}\label{eq pr stab as upper bound}
        \frac1t\int_0^t |\tcX_s - \cX_s|^2 \ud s \le 
         t^{-\vartheta} \left( C|\tcX_0-\cX_0|^2 + \mathfrak{C}(\vtstr-\vartheta) \int_0^t s^{\vartheta-1} \eg_s  \ud s \right)+o_{\eta}\left(\frac{(\log t)^{\frac12+\eta}}{t^{\frac12\wedge \vartheta}}\right).  
    \end{align}
    
\end{Proposition}

{
    \color{black}
    \proof In this proof, we denote  $(\nu,\tnu) := (\nu^{\cX},\nu^{\tcX}) $, for the reader's convenience.  We set $\epsilon>0$ s.t. $\varrho = \frac{\bstr+\epsilon}{\astr-\epsilon}\!\in (0,1)$ and $\la := \astr-\epsilon>0$. Setting $\vartheta = 1-\varrho$, we observe that $\vartheta \in (0,\vtstr)$. 
\\
{\sc Step 1}({\em Proof of $(a)$}). 
Let us define the local martingale, for $t \ge 0$,
\begin{align}\label{eq de M}
M^{(\la)}_t := \int_0^te^{\la s} H_s \ud W_s 
\end{align}
and denote $(\tau_n)$ its localizing sequence.
By It\^o's formula applied to $(x,y)\mapsto |x-y|^2$ and $(\cX, \tcX)$ one has, for $t \ge 0$,
\begin{align}\label{eq:pr:It\^o1}
&e^{\la t\wedge\tau_n} |\tcX_{t\wedge \tau_n}-\cX_{t\wedge \tau_n}|^2 
= |\tcX_0-\cX_0|^2 +  M^{(\la)}_{t\wedge \tau_n}+ \la \int_0^{t\wedge\tau_n} e^{\la s}|\tcX_s-  \cX_s|^2 \ud s 
\\
\nonumber  &\int_0^{t\wedge \tau_n} e^{\la s} \bigg[ 2(b(\tcX_{ s}  , \tnu_{s}) - b(\cX_{s}, \nu_{s}) + \etb_s\,| \tcX_{ s} - \cX_s)+ \| \s(\tcX_{s}, \tnu_{s})-\s(\cX_{s}, \nu_{s}) + \ets_s\|^2_{_F}\bigg] \ud s.
\end{align}
%
We now compute 
\begin{align}
    \| \s(\tcX_{s}, \tnu_{s})-\s(\cX_{s}, \nu_{s}) &+ \eta^\sigma_s\|^2_{_F} 
    \le (1+\tfrac{\epsilon}{4L^2})\NF{\s(\tcX_{s}, \tnu_{s}) -\s(\cX_{s}, \nu_{s})}^2 +
    (1+\tfrac{4L^2}{\epsilon})\NF{\ets_s}^2 \nonumber
    \\
    \nonumber &\le \NF{\s(\tcX_{s}, \tnu_{s}) -\s(\cX_{s}, \nu_{s})}^2 +\frac{\epsilon}2 \left(|\tcX_s - \cX_s|^2 + \cW_2^2(\tnu_s,\nu_s) \right)\\
    &\hskip 4,6cm +(1+\tfrac{4L^2}{\epsilon})\NF{\ets_s}^2 \label{eq useful cont sigma}
\end{align}
where we used $\HYP{L}$ for the second inequality.
Observing that 
\[
    2(\etb_s\,| \tcX_{ s} - \cX_s) \le \frac{2}\epsilon|\etb_s|^2 + \frac{\epsilon}{2}|\tcX_s - \cX_s|^2
\]
and combining \eqref{eq useful cont sigma} with $\HYP{C}_{\astr,\bstr}$, we get 
\begin{align}
     2(b(\tcX_{ s}  , \tnu_{s}) &- b(\cX_{s}, \nu_{s}) + \etb_s\,| \tcX_{ s} - \cX_s)+ \| \s(\tcX_{s}, \tnu_{s})-\s(\cX_{s}, \nu_{s}) + \ets_s\|^2_{_F}
    \nonumber
    \\
    &\le
    - (\astr-\epsilon)|\tcX_s - \cX_s|^2 + (\bstr+\epsilon) \cW_2^2(\tnu_s,\nu_s) 
    +C_\epsilon\left(|\etb_s|^2 + |\ets_s|^2\right) . \label{eq cont ds term}
\end{align}
Inserting the previous estimate into \eqref{eq:pr:It\^o1}, we obtain 
\begin{align}\label{eq main majo stab}
    e^{\la t\wedge \tau_n} |\tcX_{t\wedge \tau_n}-\cX_{t\wedge \tau_n}|^2 
    &\le  |\tcX_0-\cX_0|^2 +   M^{(\la)}_{t \wedge \tau_n}
    \\
    &+ 
    \int_0^{t \wedge \tau_n} e^{\la s} \left[ (\bstr+\epsilon) \cW_2^2(\tnu_s,\nu_s) 
    +C_\epsilon \eg_s \right] \ud s. \nonumber
    \end{align}
Taking expectation in the above inequality, we obtain 
\begin{align}\label{eq main majo stab esp}
    \esp{e^{\la t\wedge \tau_n}|\tcX_t-\cX_t|^2}
    &\le \esp{ |\tcX_0-\cX_0|^2}
    \\
    &+  \esp{\int_0^{t \wedge \tau_n} e^{\la s} \left[ (\bstr+\epsilon) \cW_2^2(\tnu_s,\nu_s) 
    +C_\epsilon \eg_s \right] \ud s}.
    \end{align}
    Using Fatou's Lemma for the term in the l.h.s. of the above inequality and the monotone convergence theorem for the second term in the r.h.s., we let $n \rightarrow +\infty$ to get 
    \begin{align}\label{eq main majo stab esp final}
        \esp{e^{\la t}|\tcX_t-\cX_t|^2}
        &\le \esp{ |\tcX_0-\cX_0|^2}
        +  \esp{\int_0^{t } e^{\la s} \left[ (\bstr+\epsilon) \cW_2^2(\tnu_s,\nu_s) 
        +C_\epsilon \eg_s \right] \ud s}.
        \end{align}

    Set $g(t) := \frac1t\int_0^t \esp{|\tcX_s-\cX_s|^2}\ud s$, for $t>0$ and $g(0)=\esp{|\tcX_0-\cX_0|^2}$ (observe that $g$ is continuous). Integrating the previous inequality (after multiplication by $e^{-\lambda t}$), we compute, for $t>0$,
    \begin{align}
        g(t) &\le \frac{1-e^{\la t }}{\la t }\esp{|\tcX_0-\cX_0|^2}  
        + (\bstr+\epsilon) \frac1t \int_0^t e^{-\la s }\int_0^s e^{\la u}\esp{\cW_2^2(\tnu_u,\nu_u)} \ud u\, \ud s
        \nonumber
        \\
        & + \frac1t \int_0^t C_{\epsilon}\int_0^{s } e^{\la (u-{s })} \esp{\eg_u} \ud u\, \ud s.
    \end{align}
    Applying Fubini's Theorem to the last two terms in the r.h.s$.$ of this inequality yields
        \begin{align} 
            g(t)  \le &\frac{1-e^{-\la t}}{\la t}\esp{|\tcX_0-\cX_0|^2}+  
            \varrho \frac1t \int_0^t (1-e^{\la(s-t)})\esp{\cW_2^2(\tnu_s,\nu_s)} \ud s   \nonumber
            \\
        &   + C_\epsilon \frac1t \int_0^t (1-e^{\la(s-t)})\esp{\eg_s}\ud s. 
            \label{eq conclu 1}
        \end{align}
    We observe that 
        \begin{align*}
            \cW_2^2(\tnu_s,\nu_s) &\le  \frac1s \int_0^s |\tcX_u- \cX_u|^2 \ud u,  \label{eq bound w2}
        \end{align*}
        which leads to
        \begin{align}
            g(t) \le \frac{\varrho}{ t} \int_0^t g(s) \ud s + \frac{1-e^{-\la t}}{\la t}\esp{|\tcX_0-\cX_0|^2} + \frac{C_\epsilon}t \int_0^t \esp{\eg_s} \ud s .
        \end{align}
        Applying \eqref{eq prop gronwall bis} in Lemma \ref{le pseudo gronwall main} (observe that $\frac{1-e^{-\la t}}{\la t}=\frac1t\int_0^te^{-\la s} \ud s$), we get 
        \begin{align} \label{eq end step 1}
            \hskip -0,5cm g(t)=\frac1t\int_0^t \esp{|\tcX_s - \cX_s|^2} \ud s \le 
             t^{{\varrho}-1} \int_0^t s^{-{\varrho}} \left(e^{-\la s}
             \esp{|\tcX_0-\cX_0|^2} + C_\epsilon \esp{\eg_s} \right) \ud s.
        \end{align}
        Observing then that $\int_0^t s^{-\rho}e^{-\la s}\ud s = O(1)$ proves \eqref{eq pr stab L2 upper bound} in~$(a)$.
        
        \smallskip
        \noindent {\sc Step 2} ({\em Proof of $(b)$}).  For this part, we do not need localisation, in particular, owing to \eqref{eq ass H}, $M^{(\lambda)}$ is a square integrable martingale. So, using the same computations as in the previous step, we arrive at \eqref{eq main majo stab} without localisation namely
        \begin{align}\label{eq main majo stab as}
            e^{\la t} |\tcX_{t}-\cX_{t}|^2 
            &\le  |\tcX_0-\cX_0|^2 +   M^{(\la)}_{t} + 
            \int_0^{t} e^{\la s} \left[ (\bstr+\epsilon) \cW_2^2(\tnu_s,\nu_s) 
            +C_\epsilon \eg_s \right] \ud s. 
            \end{align}
        Set $g(t) := \frac1t\int_0^t |\tcX_s-\cX_s|^2\ud s$, for $t>0$ and $g(0)=|\tcX_0-\cX_0|^2$ (observe that $g$ is continuous). Integrating the previous inequality, we compute, for $t>0$ ,
    \begin{align}
        g(t) &\le \frac{1-e^{\la t}}{\la t}|\tcX_0-\cX_0|^2+
        \frac1t \int_0^t e^{-\la s} M^{(\la)}_{s } \ud s    \nonumber
        \\
        &\qquad + (\beta+\epsilon) \frac1t \int_0^t e^{-\la {s }}\int_0^{s} e^{\la u}\cW_2^2(\tnu_u,\nu_u) \ud u\, \ud s + \frac1t \int_0^t C_{\epsilon}\int_0^{s } e^{\la (u-{s })} \eg_u \ud u\, \ud s.
    \end{align}
    Applying Fubini's Theorem to the last two terms in the r.h.s$.$ of this inequality yields
        \begin{align} 
            g(t)  &\le \frac{1-e^{-\la t}}{\la t}|\tcX_0-\cX_0|^2+ \frac1t \int_0^t e^{-\la s} M^{(\la)}_s \ud s   \nonumber
            \\
        &     \qquad+
            \varrho \frac1t \int_0^t (1-e^{\la(s-t)})\cW_2^2(\tnu_s,\nu_s) \ud s + C_\epsilon \frac1t \int_0^t (1-e^{\la(s-t)})\eg_s\ud s. 
            \label{eq conclu 1}
        \end{align}
    We observe that 
        \begin{align*}
            \cW_2^2(\tnu_s,\nu_s) &\le  \frac1s \int_0^s |\tcX_u- \cX_u|^2 \ud u,  \label{eq bound w2}
        \end{align*}
        which leads to
        \begin{align}
            g(t) \le \frac{\varrho}{ t} \int_0^t g(s) \ud s + \frac{1-e^{-\la t}}{\la t}|\tcX_0-\cX_0|^2 + \frac{C_\epsilon}t \int_0^t \eg_s \ud s + \frac1t \int_0^t e^{-\la s} M^{(\la)}_s \ud s.
        \end{align}
        Applying \eqref{eq prop gronwall bis} in Lemma \ref{le pseudo gronwall main} (observe that $\frac{1-e^{-\la t}}{\la t}=\frac1t\int_0^te^{-\la s} \ud s$), we get 
        \begin{align} \label{eq end step 1}
            \hskip -0,5cm g(t)=\frac1t\int_0^t |\tcX_s - \cX_s|^2 \ud s \le 
             t^{{\varrho}-1} \int_0^t s^{-{\varrho}} \left(e^{-\la s}|\tcX_0-\cX_0|^2 + C_\epsilon \eg_s + e^{-\la s} M^{(\la)}_s\right) \ud s.
        \end{align}
    Let us denote $t\mapsto f(t):=\int_t^{+\infty} s^{-\varrho} e^{-\la s} \ud s$. Integration by parts formula applied between $0$ and $t$ to  $s \mapsto \left(f(s)-f(t)\right) M_s^{(\la)}$ yields
        \begin{align}
            t^{{\varrho}-1}\int_0^t s^{-{\varrho}} e^{-\la s} M^{(\la)}_s \ud s
        &= t^{{\varrho}-1}\int_0^t \set{f(s)-f(t)}\ud M^{(\la)}_s \nonumber 
        \\
       &= t^{{\varrho}-1} \int_0^t f(s) e^{\la s} H_s \ud W_s 
      - t^{{\varrho}-1} f(t) \int_0^t e^{\la s} H_s \ud W_s. \label{eq job to do}
        \end{align}
        We now study the two terms in the r.h.s$.$ of the previous equation.
       
        $\blacktriangleright$ Let us focus  on the first one. For some $\eta > 0$ we define 
        \begin{align*}
    t\mapsto h_{\eta}(t) =  t^{(\frac12-\varrho)^+}(\log t)^{\frac12+\eta} \; \text{ and } \; N^\eta_t := \int_1^t \frac{f(s)e^{\la s}H_s}{h_{\eta}(s)}\ud W_s . 
\end{align*} 
One has by Burkholder-Davis-Gundy inequality for every  $T>1$,
\begin{align}
\esp{\sup_{t\in[1,T]}|\tilde N^{\eta}_t|^{2}} 
 &\le C \esp{\int_1^T |H_s|^{2} \left[\frac{f(s)e^{\lambda s}}{h_{\eta}(s)}\right]^{2} \ud s } \nonumber
 \\
 &\le C  \sup_{t \ge 1}\esp{|H_s|^{2}} \left(\int_1^T \left[\frac{f(s)e^{\lambda s}}{h_{\eta}(s)}\right]^{2} \ud s\right). \label{eq almost there}
\end{align}
We also compute by integration by parts  that
\begin{align}
\label{eq useful domination}
    f(t) = \frac1{\la}e^{-\la t}t^{-\varrho} - \frac{\varrho}{\la}\int_t^{+\infty} e^{-\la s}s^{-(\varrho+1)}\ud  s\sim  \frac{e^{-\la t}t^{-\rho}}{\la} = O(e^{-\la t}t^{-\varrho}),
\end{align}
since $t^{-(\rho+1)}= o(t^{-\rho})$ so that  
\[
    \left[\frac{f(s)e^{\lambda s}}{h_{\eta}(s)}\right]^{2} 
    = O\left(\left[\frac{s^{-\varrho}}{s^{(\frac12-\varrho)^+}(\log s)^{\frac12+\eta}}\right]^{2} \right) = O\left(s^{-(1\vee (2\varrho))}(\log s)^{-(1+\eta)}\right)
\]
which leads, using \eqref{eq almost there}, to
\begin{align*}
   \esp{ \sup_{t\ge 1}|N^{\eta}_t|^2 }=\lim_{T\to +\infty} \esp{\sup_{t\in[1,T]}| N^{\eta}_t|^{2}} 
    \le C_{\eta} \sup_{t \ge 1}\esp{|H_s|^{2}} 
    <+\infty
\end{align*}
owing to Beppo Levi's Theorem. As consequence $(N^{\eta}_t)_{t\ge 1}$ is a true ${L}^{2}$-bounded martingale, hence $\P$-$a.s$ converging toward  an ${L}^{2}$-integrable ($a.s.$ finite) variable $N^{\eta}_\infty$. 
\color{black}
The Stochastic Kronecker lemma~\ref{le kronecker gil}$(b)$   implies that 
\begin{align*}
    \frac1{h_{\eta}(t)}\int_0^t f(s)e^{\la s} H_s \ud W_s \rightarrow 0 \quad\P-a.s.
\end{align*}
We then obtain that, for any $\eta > 0$,
\begin{align*}
    t^{{\varrho}-1} \int_0^t f(s) e^{\la s} H_s \ud W_s = o_{\eta}\left(t^{-(\frac12\wedge (1-\varrho))}(\log t)^{\frac12+\eta}\right).
\end{align*}

 $\blacktriangleright$ The second term can be treated likewise by introducing
for  any $\eta > 0$:
        \begin{align*}
    t\mapsto \tilde h_{\eta}(t) =e^{\la t} t^{\frac12}(\log t)^{\frac12+\eta} \; \text{ and } \; \tilde N^\varepsilon_t := \int_1^t \frac{e^{\la s}H_s}{\tilde h_{\eta}(s)}\ud W_s 
\end{align*} 
and we obtain using again~\eqref{eq useful domination},
 that, for any $\eta > 0$, 
\begin{align*}
    t^{{\varrho}-1} f(t) \int_0^t  e^{\la s} H_s \ud W_s = o_\eta\left(t^{-(\frac12\wedge (1-\varrho))}(\log t)^{\frac12+\eta}\right).
\end{align*}
 We conclude from this step that 
\begin{align*} 
    \frac1t\int_0^t |\tcX_s - \cX_s|^2 \ud s \le 
      t^{{\varrho}-1}\left(C|\tcX_0-\cX_0|^2 + C_\epsilon\int_0^t s^{-{\varrho}} \eg_s  \ud s\right) +o_\eta\left(t^{-(\frac12\wedge (1-\varrho))}(\log t)^{\frac12+\eta}\right). \hspace{4mm} \Box
\end{align*}
}

\subsection{Convergence toward  the stationary measure}

In this section, we apply the stability result to obtain the ergodic behavior of the SID. In particular, we prove $L^2$ and $a.s.$ rate of convergence for the Wasserstein distance between the empirical emasure of the SID and the invariant distribution of the MKV SDE.

\medskip 
\noindent Let us start however with a generic convergence result.

\label{subse conv emp meas}

\begin{Theorem} \label{th conv emp sid to emp mkv no rate}
  Assume that $\hypl2$ is in force.
    \begin{enumerate}[(a)] 
        \item Then, for $\cX_0 \in \mathscr{L}^2(\cF_0,\P)$,
    \begin{align*}
    \lim_{t \rightarrow + \infty}  \esp{|\cX_t - X_t^\star|^2} = 0.
    \end{align*}
In particular,
     \begin{align*}
         \lim_{t \rightarrow + \infty} \esp{\cW^2_2(\nu^{\cX}_t,\nu^\star)} =\lim_{t \rightarrow + \infty} \frac1t\int_0^t \esp{|X^\star_s - \cX_s|^2} \ud s  = 0.
    \end{align*}
    \item Assume that $\hypas$ is in force. \textcolor{black}{Then, for $\cX_0 \in \mathscr{L}^2(\cF_0,\P)$, the following holds} 
    \begin{align}\label{eq res conv no speed continous time new}
         \lim_{t \rightarrow + \infty}\cW_2(\nu^{\cX}_t,\nu^\star) =\lim_{t \rightarrow + \infty} \frac1t\int_0^t |X^\star_s - \cX_s|^2 \ud s=0 \;\;\P-a.s. 
    \end{align}
\end{enumerate}
\end{Theorem}
{
    \color{black}
    \proof
    {\sc Step 1.}({\em Proof of $(a)$}) We first observe that $X^\star$ satisfy \eqref{eq de tXi} with 
    \begin{align}
        \etb_t = b(X^\star_t,\nu^\star) - b(X^\star_t,\nu_t^{X^\star})
            \;\text{ and }\;
        \ets_t = \s(X^\star_t,\nu^\star) - \s(X^\star_t,\nu_t^{X^\star}). 
    \end{align}
    From $\HYP{L}$, we get, according to \eqref{eq pr def delta}, 
    \begin{align} \label{eq starting observation as no rate}
        \eg_t \le C \cW^2_2(\nu^{X^\star}_t,\nu^\star).
    \end{align}
    Then, combining \eqref{eq le L2 conv mvk} with Lemma \ref{le Cesaro gil} with $g(t)=t^{-\rho}$, we get
    \begin{align}\label{eq useful interm}
        \esp{\eg_t} \rightarrow 0 \quad \text{ and }\quad t^{{\varrho}-1} \int_0^t s^{-{\varrho}} \esp{\eg_s}  \ud s \rightarrow 0.
    \end{align}
    Inserting the previous limit into \eqref{eq pr stab L2 upper bound}, we obtain
    \begin{align}\label{eq L2 done}
        \esp{\cW^2_2(\nu^{\cX}_t,\nu^{X^\star}_t)} \le \frac1t\int_0^t \esp{|X^\star_s - \cX_s|^2} \ud s = o(1).
     \end{align}
     Observing that
 \begin{align}
    \esp{\cW^2_2(\nu^{\cX}_t,\nu^{\star})} \le 2\esp{\cW^2_2(\nu^{\cX}_t,\nu^{X^\star}_t)} + 2 \esp{\cW^2_2(\nu^{\star},\nu^{X^\star}_t)},
 \end{align}
 we obtain 
 $$\lim_{t \rightarrow + \infty} \esp{\cW^2_2(\nu^{\cX}_t,\nu^\star)} = 0,$$
 by invoking \eqref{eq le L2 conv mvk} again.
 \\
     1.b Going back to the proof of Proposition \ref{pr stability for sid} , \eqref{eq main majo stab esp final} reads, with $\tcX=X^\star$, $\nu = \nu^\cX$ and $\tilde{\nu}=\nu^{X^\star}$, 
     \begin{align*}
        \esp{|X^\star_t-\cX_t|^2}
        \le e^{-\la t}  \esp{|X^\star_0-\cX_0|^2} 
        + \int_0^t e^{\la (s-t)} \esp{(\bstr+\epsilon) \cW^2_2(\nu^{\cX}_s,\nu^{X^\star}_s) 
        +C_\epsilon \eg_s} \ud s,
        \end{align*}
        for some $\epsilon > 0$.
 Combining the previous inequality with \eqref{eq useful interm}, \eqref{eq L2 done} and Lemma \ref{le Cesaro gil}, we obtain the second claim in (a). 
 \\
 {\sc Step 2.}({\em Proof of $(b)$}) 
 i) Assume first that $\cX_0 \in \mathscr{L}^4(\cF_0,\P)$. 
 We first observe that applying Lemma \ref{le unif int sid} and \eqref{eq integrability X star}, under $\hypas$, yields that 
 \begin{align}
    \sup_{t \ge 0}\esp{|\cX_t|^{4} + |X^{\star}_t|^{4}} \le C.
 \end{align}
We thus deduce easily that \eqref{eq ass H} holds true with $H$ simply given here by 
$$H_t = 2(X^\star_t -\cX_t | \sigma(X^\star_t,\nu^{X^\star}_t) - \sigma(\cX_t,\nu^{\cX}_t)).$$ 
Combining \eqref{eq starting observation as no rate} and \eqref{eq le almost sure conv mvk} with Lemma \ref{le Cesaro gil} , we have 
    \begin{align}
        t^{{\varrho}-1} \int_0^t s^{-{\varrho}} \eg_s  \ud s \rightarrow 0.
    \end{align} 
We can then invoke Proposition \ref{pr stability for sid} and in particular insert the previous limit into \eqref{eq pr stab as upper bound} to obtain.
\begin{align}\label{eq almost done}
    \cW^2_2(\nu^{\cX}_t,\nu^{X^\star}_t) \le \frac1t\int_0^t |X^\star_s - \cX_s|^2 \ud s = o(1).
 \end{align}
 Observing that
 \begin{align}
    \cW^2_2(\nu^{\cX}_t,\nu^{\star}) \le 2\cW^2_2(\nu^{\cX}_t,\nu^{X^\star}_t) + 2\cW^2_2(\nu^{\star},\nu^{X^\star}_t)
 \end{align}
 and combining \eqref{eq almost done} with \eqref{eq le almost sure conv mvk}, we conclude for this step that 
 \begin{align*}
    \lim_{t \rightarrow + \infty}\cW_2(\nu^{\cX}_t,\nu^\star) =\lim_{t \rightarrow + \infty} \frac1t\int_0^t |X^\star_s - \cX_s|^2 \ud s=0 \;\;\P-a.s. 
\end{align*}
ii) Assume now that $\cX_0 \in \cL^2(\cF_0,\P)$ only. For any positive integer $A$, we denote $(\cX_t^A)_{t \ge 0}$ the self-interacting diffusion \eqref{eq de sid} started from $\cX_0^A:=\cX_0\1_{\set{|\cX_0|\le A}}$. One observes that $(\cX_t)_{t\ge 0}$ and  $(\cX_t^A)_{t\ge0}$ coincides on the set $\set{|\cX_0|\le A}$. 
Applying the result of the previous step, we obtain  
\begin{align*}
    \lim_{t \rightarrow + \infty}\cW_2(\nu^{\cX}_t,\nu^\star) =\lim_{t \rightarrow + \infty} \frac1t\int_0^t |X^\star_s - \cX_s|^2 \ud s=0 \;\;\text{ on } \;\set{|\cX_0|\le A}\setminus \cN_A,
\end{align*}
where $\cN_A$ is a negligible set. We observe then that $\Omega = \cup_{A \in \mathbb{N}^*}\set{|\cX_0|\le A}$ and that $\cup_{A \in \mathbb{N}^*} \cN_A$ is a negligible set as well, so that \eqref{eq res conv no speed continous time new} is proved.
\eproof
}


\bigskip 

\noindent Both for the $L^2$ and $a.s.$ case, it is possible to identify convergence
rate and thus improve substantially the result of the previous section.

\paragraph{$L^2$ rate}

\begin{Theorem} \label{th l2 conv with a rate sid to nu star}
    Assume that  $\hypl2$ is in force. Let $\cX_0 \in \mathscr{L}^2(\cF_0,\P)$.
    For any $\vartheta \in (0,\vtstr)\setminus \set{\zeta}$, the following holds for $\theta=(\vtstr-\vartheta,|\vartheta-\zeta|)$
  \begin{align*}
        \esp{\cW^2_2(\nu^{\cX}_t,\nu^\star)} \le  \frac1t\int_0^t \esp{|X^\star_s - \cX_s|^2} \ud s  = 
        O \left(t^{-\vartheta}\esp{|X^\star_0-\cX_0|^2} + \mathfrak{C}(\theta) t^{-\zeta}(\log t)^{\kappa}\right),
   \end{align*}
   where $(\zeta,\kappa)=\left(\frac{2 \pstr - 1}{2(2(d+2)+(2 \pstr - 1)(d+3))},0\right)$ . \\
   If furthermore $\HYP{\sigma_\infty}$ holds, then one can set $(\zeta,\kappa)=\left(\frac1{2(d+3)},\frac{d+2}{2(d+3)}\right)$.\\
\smallskip
    Moreover, the same upper bounds hold true for \esp{|\cX_t - X_t^\star|^2}.
\end{Theorem}

\color{black}
\proof 1. We first observe that $X^\star$ satisfy \eqref{eq de tXi} with 
\begin{align}
    \etb_t = b(X^\star_t,\nu^\star) - b(X^\star_t,\nu_t^{X^\star})
        \;\text{ and }\;
    \ets_t = \s(X^\star_t,\nu^\star) - \s(X^\star_t,\nu_t^{X^\star}). 
\end{align}
From $\HYP{L}$, we get, combining \eqref{eq pr def delta} and Corollary \ref{co L2 conv rate emp meas to inv meas}, 
    \begin{align} \label{eq starting observation L2 conv rate sid to mkv} 
        \esp{\eg_t} \le C \esp{\cW^2_2(\nu^{X^\star}_t,\nu^\star)} =O(t^{-\zeta}(\log t)^{\kappa}).
    \end{align}
For any $\vartheta \in (0,\vtstr)$, it follows moreover from \eqref{eq pr stab L2 upper bound}
\begin{align*}
    \frac1t\int_0^t \esp{|X^\star_s -\cX_s|^2} \ud s 
    \le C  t^{-\vartheta} \left( \esp{|X^\star_0-\cX_0|^2} + \mathfrak{C}(\vtstr-\vartheta) \int_0^t s^{\vartheta-1} \mathbb{E}[\eg_s] \ud s \right).
\end{align*}
The proof is concluded for this step by inserting \eqref{eq starting observation L2 conv rate sid to mkv} into the previous inequality.
\\
2. Going back to the proof of Proposition \ref{pr stability for sid} , \eqref{eq main majo stab esp final} reads, with $\tcX=X^\star$, $\nu = \nu^\cX$ and $\tilde{\nu}=\nu^{X^\star}$, 
\begin{align*}
   \esp{|X^\star_t-\cX_t|^2}
   \le e^{-\la t}  \esp{|X^\star_0-\cX_0|^2} 
   + \int_0^t e^{\la (s-t)} \esp{(\bstr+\epsilon) \cW^2_2(\nu^{\cX}_s,\nu^{X^\star}_s) 
   +C_\epsilon \eg_s} \ud s,
   \end{align*}
   for some $\epsilon > 0$. 
   Since $s \mapsto e^{\la s} s^{-a}$ is increasing for $s$ large enough for any $a$ and $e^{-\la t}=o(t^{-a})$ for any $a>0$, we obtain, using the result of the previous step,  
   \begin{align*}
    \esp{|X^\star_t-\cX_t|^2} = O \left(t^{-\vartheta}\esp{|X^\star_0-\cX_0|^2} + \mathfrak{C}(\theta) t^{-\zeta}(\log t)^{\kappa}\right) .
   \end{align*}
\eproof 
\color{black}

\begin{Remark}
    The convergence with a rate under $\hypl2$ is already obtained in \cite[Theorem 2.2]{du2023empirical}. The rate obtained by assuming $\HYP{\sigma_\infty}$ is new, to the best of our knowledge.
\end{Remark}


\paragraph{Almost sure rate}

To the best of our knowledge, the following results are new in this setting.

\begin{Theorem} \label{th conv emp sid to emp mkv with a rate}
    Assume that $\hypas$ is in force and that $\sigma(\cdot,\nstr)$ is uniformly elliptic. Let $\cX_0 \in \mathscr{L}^2(\cF_0,\P)$.
    Then, the following holds, for any $\vartheta \in (0,\vtstr)\setminus\set{\zeta}$, any small $\eta, \eta' >0$,
    \begin{align*}
         &\cW_2(\nu^{\cX}_t,\nu^\star) \le  \frac1t\int_0^t |X^\star_s - \cX_s|^2 \ud s
         \\
         &= 
         O \left(t^{-\vartheta}   |\tcX_0-\cX_0|^2 \right)
         + o_{|\vartheta-\zeta|,\eta'}
        \left(t^{-\zeta}(\log t)^{\frac12+\eta'}\right)
        +o_{\eta }\left(\frac{(\log t)^{\frac12+\eta}}{t^{\frac12\wedge \vartheta}}\right)
         \;\;\P-a.s. 
    \end{align*}
    where $\zeta = \frac{{a^\star}^2}{2(2+{a^\star}) \{2(d+2) + (d+3)({a^\star}+d) \}}$. \\ If moreover $\HYP{\sigma_\infty}$ holds,  then one can set $\zeta = \frac1{2(d+3)}$.
\end{Theorem}

\color{black}
\proof We assume that $\cX_0 \in \mathscr{L}^4(\cF_0,\P)$; the proof for the case $\cX_0 \in \mathscr{L}^2(\cF_0,\P)$ only is obtained following the same arguments as in step 2.(ii) of the proof of Theorem \ref{th conv emp sid to emp mkv no rate}.
\\
We now observe that $X^\star$ satisfy \eqref{eq de tXi} with 
\begin{align}
    \etb_t = b(X^\star_t,\nu^\star) - b(X^\star_t,\nu_t^{X^\star})
        \;\text{ and }\;
    \ets_t = \s(X^\star_t,\nu^\star) - \s(X^\star_t,\nu_t^{X^\star}). 
\end{align}
From $\HYP{L}$, we get, combining \eqref{eq pr def delta} and Corollary \ref{co as conv rate emp meas to inv meas}, 
    \begin{align} \label{eq starting observation as conv rate sid to mkv} 
        \eg_t \le C \cW^2_2(\nu^{X^\star}_t,\nu^\star) =o_{\eta'}(t^{-\zeta}(\log t)^{\frac12+\eta'})\;,\; \text{ for any small } \eta'>0.
    \end{align}
We compute, for $\vartheta \neq \zeta$,
\[
   t^{-\vartheta}\int_0^t s^{\vartheta-1}\eg_s \ud s = o_{|\vartheta-\zeta|,\eta'}
   \left(t^{-\zeta}(\log t)^{\frac12+\eta'}\right)
    \]
and obtain from \eqref{eq pr stab as upper bound}
    \begin{align*}
        \frac1t\int_0^t |\tcX_s - \cX_s|^2 \ud s \le &
        C  t^{-\vartheta}   |\tcX_0-\cX_0|^2 
         + o_{|\vartheta-\zeta|,\eta'}
        \left(t^{-\zeta}(\log t)^{\frac12+\eta'}\right)
        +o_{\eta }\left(\frac{(\log t)^{\frac12+\eta}}{t^{\frac12\wedge \vartheta}}\right).  
    \end{align*}
    
\eproof
\color{black}

\section{Numerical scheme}


\label{se numerical scheme}
\noindent For the sequence $(\Gamma_n)_{n \ge 0}$ introduced in~\eqref{eq:stepdecrease} and every instant $t\ge 0$, we define 
\[
\t= \Gamma_n \;  \mbox{ for} \; t\!\in [\Gamma_n , \Gamma_{n+1}).
\]

\noindent We will also introduce the right continuous inverse of $s\mapsto \us$ defined by
\[
N(s) = \min\{k: \Gamma_{k+1}>s\}
\]
which satisfies 
\[
\Gamma_{N(s)}= \us \mbox{ for every } s>0.
\]

The numerical scheme we consider has been defined in~\eqref{eq:SDEdisc0}-\eqref{eq:SDEdisc3}. We introduce now its associated continuous time version:
\begin{align}
\label{eq:Eulerstar}
\bar \cX_t & = \cX_0 +\int_0^t b(\bar \cX_{\us}, \bar \nu_{\us})\ud s +\int_0^t \sigma(\bar \cX_{\us}, \bar \nu_{\us})\ud W_s,\\
\bar \nu_{t}  &:=\frac{1}{t} \int_0^{t} \delta_{\bar \cX_{\us}} \ud s,
\end{align}
with the convention that, when $\t =0$ (i.e. $t\!\in [0, \g_1)$), $\bar \nu_{t} = \delta_{\cX_0}$ and, more generally, $\frac{1}{\t}\int_0^{\t} x(s)\ud s = x(0)$ when $\t=0$). 

Note that this definition is  consistent with~\eqref{eq:SDEdisc1} at times $t=  \Gamma_n$,  $n\ge1$. The fact that we have inserted $ \bar{\nu}_{\t}$ in the sheme ensures its tractability although it adds some additional technicality when computing {\em a priori} estimates. Moreover for technical matter we also introduce the empirical measure $\nu^{\bar \cX}$ of the continuous time Euler scheme, naturally defined at time $t$ by $\nu^{\bar \cX}_t= \frac1t \int_0^t \delta_{\bar \cX_s}\ud s$.

\smallskip To shorten some statements,  we introduce the notation $\hat t$: when appearing in  a statement, it means that this statement is valid for $t$ and $\t$.

\medskip
\noindent The main result of this section are the theorem and its corollary below.

\begin{Theorem}\label{thm:tpsdiscretsmulable} 
  Assume  that the step sequence $(\g_n)_{n\ge1}$ satisfies~\eqref{eq:stepdecrease}. \\
\smallskip 
\noindent $(a)$ {\em $L^2$-convergence}. Let $\hypl2$ hold.
Then
\begin{align}\label{eq th conv esp}
\esp{|{\cal X}_t-\bar {\cal X}_{\hat t}|^2 }\to 0 \quad \mbox{ as }\quad t\to +\infty,
\end{align}
and, a fortiori, $\displaystyle   \frac{1}{t}\int_0^t \esp{|{\cal X}_s-\bar {\cal X}_{\hat s}|^2}\ud s \to 0 , \; $
\begin{align}\label{eq th conv esp2}
\esp{{\cal W}_2^2\big(\nu^{\cX}_t  ,  \nu^{\bar\cX}_{t} \big)} \to 0  \quad\mbox{ and }\quad \esp{{\cal W}_2^2\big(\nu^{\cX}_t  , \bar \nu_{t} \big)} \to 0 \quad\mbox{as}\quad t\to +\infty.
\end{align}

\noindent $(b)$ {\em $A.s.$-convergence}.  Let $\hypas$ hold. If 
\begin{align}\label{eq th conv ps1}
 \sum_{n\ge 1} \frac{\g_{n+1}^2}{\Gamma_n}<+\infty,
\end{align} then 
\[
\frac{1}{t}\int_0^t\big({\cal X}_s-\bar {\cal X}_{\hat s}\big)^2\ud s \to 0  \quad \P\mbox{-}a.s.  \quad \mbox{ as }\quad t\to +\infty
\]
and, a fortiori,
\begin{align}\label{eq th conv ps2}
{\cal W}_2^2\big(\nu^{\cX}_t  , \bar \nu_{t} \big)  \to 0 \quad\mbox{ and }\quad {\cal W}_2\big(\nu^{\cX}_t  , \bar \nu_{t} \big)  \to 0  \quad \P\mbox{-}a.s.  \quad \mbox{ as }\quad t\to +\infty.
\end{align}
\end{Theorem}

\noindent As a straightforward consequence of Theorem~\ref{thm:tpsdiscretsmulable} and Proposition~\ref{pr empirical measure no rate}, we have the following results.

\begin{Corollary} \label{co conv scheme no rate}
  Assume that the step sequence $(\g_n)_{n\ge1}$ satisfies~\eqref{eq:stepdecrease}.

\smallskip
\noindent $(a)$ Let $\hypl2$ hold. Then, the following holds
\[
 \esp{{\cal W}_2^2\big(\bar \nu_{t} ,\nu^\star\big) }\to 0 \quad \mbox{ as }\quad t\to +\infty.
\]

\noindent $(b)$  Let $\hypas$ hold. Assume that  $\displaystyle \sum_{n\ge 1} \frac{\g_{n+1}^2}{\Gamma_n}<+\infty$, then 
\[
{\cal W}_2^2\big(\bar \nu_{t} ,\nu^{\star}\big)  \to 0 \quad \P\mbox{-}a.s. \quad \mbox{ as }\quad t\to +\infty.
\]
\end{Corollary}

\bigskip
\noindent We first need the two following technical lemmas.

\begin{Lemma}\label{prop:cXLpbounded} 
  Assume that $\HYP{MV}_{\bar{p},\bar{K},\bar\alpha,\bar\beta}$ for $\bar p \ge 1$ and with $\bar \a > \bar \b\ge 0$ is in force.
     Then for all $p \in [1,\bar{p}]$, as soon as $\cX_0 \in \mathscr{L}^{2p}(\cF_0,\P)$,
we have that 
\[
  \sup_{n\ge 0} \|\bar \cX_{\Gamma_n}\|_{2p} \le \bar{C} (1+ \|\cX_0\|_{2p}).
\]
where $\bar{C}$ is a positive constant that depends on $\bar K,\bar \a,\bar \b$ and $p$.


\end{Lemma}

\noindent The proof of this key Lemma is postponed to Section \ref{se proof of boundedness marginal scheme}

\begin{Lemma}[Moment and regularity control]\label{lem:MRcont} 
  Assume that $\HYP{MV}_{\bar{p},\bar{K},\bar\alpha,\bar\beta}$ for $\bar p \ge 1$ and with $\bar \a > \bar \b\ge 0$ is in force.
     Then for all $p \in [1,\bar{p}]$, as soon as $\cX_0 \in \mathscr{L}^{2p}(\cF_0,\P)$,

  \begin{align} \label{eq:L^2boundXstar}
 (i)  &\;\sup_{t\ge 0} \|\bar \cX_t\|_{2p} \le \bar{C}' (1+\|\cX_0\|_{2p})<+\infty,\\
 (ii) &\; \label{eq:L^2ctrlXstar} \sup_{n\ge 1} \big \|\sup_{t\in (\Gamma_n, \Gamma_{n+1})}|\bar \cX_t- \bar \cX_{\Gamma_n}|\,\big\|_{2p}  \le  \bar{C}'(1+\|\cX_0\|_{2p})\g_{n+1} \\
\nonumber \mbox{and} & \sup_{t\ge 0, t\neq \t} \frac{\|\bar \cX_t- \bar \cX_{\t}\|_{2p}}{(t-\t)^{\frac 12}} \le  \bar{C}''(1+\|\cX_0\|_{2p})
 \end{align}
where $\bar{C}',\bar{C}''$ are positive constants depending upon $\bar K,\bar \a,\bar \b$ and $p$.
  \end{Lemma}
  %
  {
    \color{black}
  \proof 
  
  \noindent $(i)$ follows from  $(ii)$ combined with Lemma~\ref{prop:cXLpbounded}. 
  
  \smallskip
  \noindent As for $(ii)$, we note that 
  \[
  \bar \cX-\bar \cX_{\t} = b(\bar \cX_{\t}, \bar \nu_{\t} )(t-\t) +\sigma(\bar \cX_{\t}, \bar \nu_{\t} )(W_t-W_{\t})
  \]
  where $W_t-W_{\t}$ is clearly independent of ${\cal F}_{\t}$. Hence, if $\t=\Gamma_n$ and $t\!\in (\Gamma_n, \Gamma_{n+1})$,
  \[
  \esp{\sup_{t\in (\Gamma_n, \Gamma_{n+1})} |\bar \cX-\bar \cX_{\t} |^p} \le 2^{p-1}\Big( \esp{\big| b(\bar \cX_{\t}, \bar \nu_{\t} ) \big|^p} (t-\t)^p +  \,\esp{\| \sigma(\bar \cX_{\t}, \bar \nu_{\t} )\|^p_{_F}} \esp{\sup_{t\in (\Gamma_n, \Gamma_{n+1})} |W_t-W_{\t}|^p}\Big).
  \]
  Now $ \E\,\sup_{t\in (\Gamma_n, \Gamma_{n+1})} |W_t-W_{\t}|^p=  \E\,\sup_{t\in (\g_{n+1})} |W_t|^p= 2 \E\, |W_{\g_{n+1}}|^p= c^W_p \g_{n+1}^{p/2}$ since $W\stackrel{d}{=}-W$ and $\sup_{0\le s \le t} W_s \stackrel{d}{=} |W_t|$    so that, 
  \[
  \E\, \sup_{t\in (\Gamma_n, \Gamma_{n+1})} |\bar \cX-\bar \cX_{\t} |^p \le 2^{p-1}\Big( \E\, \big| b(\bar \cX_{\t}, \bar \nu_{\t} ) \big|^p (\g_{n+1}^pp +  c^W_p\E\, \| \sigma(\bar \cX_{\t}, \bar \nu_{\t} )\|^p_{_F}\g_{n+1}^{p/2}  \Big).
  \]
  Hence
  \begin{align*}
  \frac{  \E\, \sup_{t\in (\Gamma_n, \Gamma_{n+1})} |\bar \cX-\bar \cX_{\t} |^p }{\g_{n+1}^{p/2}} & \le C  \big( 1+ \E\, |\bar \cX_{\t}|^p + \E {\cal W}_2(\bar \nu_{\t},\d_0)^p\big)\\
  & \le  C  \Big(1+  C^p_{p,\bar \cX} + \frac{1}{\t} \int_0^{\t} \E\, |\bar \cX_{\us}|^p \ud s\Big) \\
  & \le C  (1+  2\,C^p_{p,\bar \cX} )  (1+\|\cX_0\|_p)^p <+\infty,
  \end{align*}
  where for the last inequality we used Proposition~\ref{prop:cXLpbounded}.
  \eproof
  }

\bigskip
 
  \color{black}
  \noindent {\bf Proof of Theorem~\ref{thm:tpsdiscretsmulable}}
\noindent {Step 1.} {\em $L^2$-convergence}. \\$\blacktriangleright$ We first observe that $\bar{\cX}$ satisfies~\eqref{eq de tXi} with 
  \begin{align}
      \etb_t = b(\bar{\cX}_{\t},\bar{\nu}_{\t}) - b(\bar{\cX}_{t},{\nu}^{\bar{\cX}}_{t})
          \;\text{ and }\;
      \ets_t = \s(\bar{\cX}_{\t},\bar{\nu}_{\t}) - \s(\bar{\cX}_{t},{\nu}^{\bar{\cX}}_{t}). 
  \end{align}
  From $\HYP{L}$, we get for $\eg_t$ according to its definition~\eqref{eq pr def delta}, 
    \begin{align} \label{eq starting observation}
        \eg_t \le C\left(|\bar{\cX}_{\t}-\bar{\cX}_{t}|^2+\cW^2_2(\bar{\nu}_{\t},{\nu}^{\bar{\cX}}_{t})\right).
    \end{align}
  We now study the second term in the r.h.s. of the above inequality. We first note that
\[
\frac{1}{t}\int_0^t \delta_{\bar \cX_s}\ud s = \frac{\t}{t}\cdot\frac{1}{\t} \int_0^{\t} \delta_{\bar \cX_s}\ud s   + \Big(1-\frac{\t}{t}\Big) \cdot \frac{1}{t-\t}\int_{\t}^t \delta_{\bar \cX_s}\ud s.
\]
By (a standard argument about the) convexity of $\mu \mapsto {\cal W}_2^2(\mu, \nu)$ on ${\cal P}_2(\R^d)$, one has 
\begin{align*}
{\cal W}^2_2 \bigg(\frac{1}{t}\int_0^t \delta_{\bar \cX_s}\ud s ,\frac{1}{\t} \int_0^{\t} \delta_{\bar \cX_{\us}}\ud s\bigg)& \le \tfrac{\t}{t} {\cal W}_2^2\bigg(\frac{1}{\t} \int_0^{\t} \delta_{\bar \cX_s}\ud s , \frac{1}{\t}\int_0^{\t} \delta_{\bar \cX_{\us}}\bigg) \\
\\&\qquad+ \Big(1-\frac{\t}{t}\Big) {\cal W}^2_2\bigg(\frac{1}{t-\t}\int_{\t}^t \delta_{\bar \cX_s}\ud s, \frac{1}{\t}\int_0^{\t} \delta_{\bar \cX_{\us}}\ud s\bigg)
\end{align*}
and 
\[
 {\cal W}^2_2\bigg(\frac{1}{t-\t} \int_{\t}^t \delta_{\bar \cX_s}\ud s, \frac{1}{\t} \int_0^{\t} \delta_{\bar \cX_{\us}}\ud s\bigg)\le 2\bigg( {\cal W}^2_2\bigg(\frac{1}{t-\t} \int_{\t}^t \delta_{\bar \cX_s}\ud s,   \delta_{0}\bigg)+ {\cal W}^2_2\bigg(\delta_0, \frac{1}{\t} \int_0^{\t} \delta_{\bar \cX_{\us}}\ud s\bigg) \bigg).
\]

As ${\cal W}_2(\mu, \d_0)^2 = \int |\xi|^2\nu(d\xi)$ and ${\cal W}_2^2 \Big(\frac {1}{t}\int_0^t\delta_{x(s)}\ud s, \frac{1}{t}\int_0^t\delta_{y(s)}\ud s\Big)\le \frac{1}{t}\int_0^t |x(s)-y(s)|^2\ud s$,  we finally get 
\begin{align}\label{eq:W2}
\nonumber {\cal W}^2_2 \bigg(\frac{1}{t}\int_0^t \delta_{\bar \cX_s}\ud s ,\frac{1}{\t} \int_0^{\t} \delta_{\bar \cX_{\us}}\ud s\bigg) &\le \frac{1}{t}\int_0^t |\bar \cX_s-\bar \cX_{\us}|^2 \ud s + \frac{2}{t} \int_{\t}^t |\bar \cX_s|^2\ud s \\
&\hskip3.5cm + 2 \frac{t-\t}{t} \cdot \frac{1}{\t} \int_0^{\t} |\bar \cX_{\us}|^2\ud s
\end{align}
(with our convention  $\frac{1}{\t}\int_0^{\t} x(s) \ud s=x(0)$ when $\t =0$ and $x$ is continuous at $t=0$).
So that finally, we obtain 
\begin{align}\label{eq majo conv without a rate}
  \eg_t &\le C\, \mathfrak{E}_t 
  \end{align}
  with 
  \begin{align}\label{eq de main disc error}
    \mathfrak{E}_t := |\bar \cX_t-\bar \cX_{\t}|^2 +\frac{1}{t}\int_0^t |\bar \cX_s- \bar \cX_{\us}|^2\ud s +\frac{1}{t}\int_{\t}^t |\bar \cX_s|^2\ud s  
   + \frac{t-\t}{t} \cdot \frac{1}{\t}\int_0^{\t} |\bar \cX_{\us}|^2 \ud s.
  \end{align}
$\blacktriangleright$  From~\eqref{eq pr stab L2 upper bound} (see Proposition~\ref{pr stability for sid}) combined with~\eqref{eq majo conv without a rate}, we obtain 
\begin{align}\label{eq temp control disc error}
  \frac1t\int_0^t \esp{|\bar{\cX}_s -\cX_s|^2} \ud s 
        \le C_{\!\varrho} t^{{\varrho}-1} \left( \int_0^t s^{-{\varrho}} \mathbb{E}[\mathfrak{E}_s] \ud s \right).
\end{align}
By~\eqref{eq:L^2boundXstar} and~\eqref{eq:L^2ctrlXstar} from Lemma~\ref{lem:MRcont}$(ii)$-$(iii)$ applied with $p=1$, we get, taking expectation in~\eqref{eq de main disc error}, for every $t\ge \g_1$,
\[
 \esp{\mathfrak{E}_t}  \le (1+\|\cX_0\|_2)^2 \bigg( (C'_{2,Y})^2\big((t-\t) +\frac{1}{t}\int_0^t (s-\us)\ud s\big) +2(C_{2,Y})^2\Big(\frac{t-\t}{t}\Big) \bigg).
\]

One checks that, for $s\ge \g_1$,
\[
\frac{1}{t}\int_0^t (s-\us)\ud s \le \frac{1}{\t} \sum_{k=1}^{N(t)+1} \g_k^2 =  \frac{1}{\Gamma_{N(t)}} \sum_{k=1}^{N(t) } \g_k^2+ \frac{\g_{N(t)+1}^2}{\Gamma_{N(t)}}\to 0 \quad \mbox{as}\quad t\to +\infty
\]
since $\g_n \to 0$ as $n\to +\infty$. Moreover $\frac{t-\t}{t}\le \frac{\g_{N(t)+1}}{\g_1}\to 0$ as $t\to +\infty$. Finally, 
\begin{align}
  \esp{\mathfrak{E}_t}\to 0 \quad \text{as} \quad t\to +\infty, \label{eq error goes to 0}  
\end{align}
and using {Lemma~\ref{le Cesaro gil}} with $g(t)=t^{-\rho}$, we get
\begin{align}\label{eq useful interm disc error}
     t^{{\varrho}-1} \int_0^t s^{-{\varrho}} \esp{\mathfrak{E}_s}  \ud s \rightarrow 0.
\end{align}
Inserting the previous estimate into~\eqref{eq temp control disc error} yields 
\begin{align}\label{eq L2 disc error done} 
  \esp{\cW^2_2(\nu^{\cX}_s,\nu^{\bar{\cX}}_s)} \rightarrow 0.
\end{align}
$\blacktriangleright$  Going back to the proof of Proposition~\ref{pr stability for sid}, we take expectation in~\eqref{eq main majo stab} with $\nu = \nu^\cX$ and $\tilde{\nu}=\nu^{\bar{\cX}}$ to obtain 
\begin{align*}
   \esp{|\bar{\cX}_t-\cX_t|^2}
   \le 
    \int_0^t e^{\la (s-t)} \esp{(\beta+\epsilon) \cW^2_2(\nu^{\cX}_s,\nu^{\bar{\cX}}_s) 
   +C_\epsilon \mathfrak{E}_s} \ud s.
   \end{align*}
Combining the previous inequality with~\eqref{eq useful interm disc error},~\eqref{eq L2 done} and {Lemma~\ref{le Cesaro gil}}, we obtain the second claim in~\eqref{eq th conv esp}.
This concludes the proof of step~$(a)$  in Theorem~\ref{thm:tpsdiscretsmulable} when $\hat t =t$. One concludes by using Lemma~\ref{lem:MRcont}$(ii)$ with $p=1$ when $\hat t= \t$.

\medskip
\noindent {Step 2.}  {\em $A.s.$-convergence}. We assume that $\cX_0 \in \mathscr{L}^4(\cF_0,\P)$; the proof for the case $\cX_0 \in \mathscr{L}^2(\cF_0,\P)$ only is obtained following the same arguments as in step 2.(ii) of the proof of Theorem \ref{th conv emp sid to emp mkv no rate}. 
\\
We first check that 
$$
H_t = 2(\bar{\cX}_t -\cX_t | \s(\bar{\cX}_{\t},\bar{\nu}_{\t}) - \s(\cX_t,\nu^{\cX}_t))
$$
satisfies~\eqref{eq ass H}.
Indeed, from this, we derive  that 
 \begin{align*}
  |H_t|&\le 2L|\cX_t-\bar \cX_t|(|\cX_t-\bar \cX_{\t}|+ {\cal W}_2(\nu^{\cX}_t,  \bar \nu_{\t})\big)\\
  & \le  C\Big(  |\cX_t-\bar \cX_t|^2 +|\bar \cX_t-\bar \cX_{\t}|^2 + {\cal W}^2_2(\nu^{\cX}_t,  \bar \nu_{\t})^2\Big)\\
  &\le C\Big( |\cX_t|^2+|\bar \cX_t|^2 +|\bar \cX_t|^2+|\bar \cX_{\t}|^2 + {\cal W}_2(\nu^{\cX}_t,\d_0)+ {\cal W}_2^2(\bar \nu_{\t},\d_0)^2\Big)\\
  &\le C\Big(  |\cX_t|^2+|\bar \cX_t|^2 +|\bar \cX_t|^2+|\bar \cX_{\t}|^2 + \frac{1}{t}\int_0^t |\cX_s|^2\ud s+
  \frac{1}{\Gamma_{N(t)}}\sum_{k=1}^{N(t)}\g_\ell |\bar \cX_{\Gamma_{\ell-1}}|^2\Big).
\end{align*}
Consequently, under $\hypas$,
\[
\esp{|H_t|^{2}}\le C(1+\E|\cX_0|^{4})<+\infty.
\]
It follows from~\eqref{eq pr stab as upper bound} from Proposition~\ref{pr stability for sid} that
\begin{align}\label{eq appli pr stab}
  \frac1t\int_0^t |\tcX_s - \cX_s|^2 \ud s \le 
  C_{\!\varrho}  t^{{\varrho}-1}  \int_0^t s^{-{\varrho}} \mathfrak{E}_s \ud s +o_{\eta}\big (t^{-(\frac12\wedge(1-\varrho))}(\log t)^{\frac 12+\eta}\big) .  
\end{align}
First note that 
\[
\int_{\g_1}^{+\infty} \frac{\E\, |\bar \cX_s-\bar \cX_{\us}|^2}{s} \ud s\le C(1+ \|\cX_0\|_2)^2 \int_1^{+\infty} \frac{s-\us}{s}\ud s\le \sum_{n\ge 1}\frac{\g^2_{n+1}}{\Gamma_{n}}<+\infty.
\]
Consequently, owing to Fubini's Theorem the non-negative  random variable (r.v.)   \linebreak $\int_{\g_1}^{+\infty} \frac{ |\bar \cX_s-\bar \cX_{\us}|^2}{s} \ud s<+\infty$ $\P$-$a.s.$ since it is integrable. Kronecker's Lemma, see Lemma~\ref{le kronecker gil}$(a)$ in the Appendix, then implies both 
\begin{align} 
& \frac{1}{t} \int_0^{t}   |\bar \cX_s-\bar \cX_{\us}|^2  \ud s\to 0 \quad \P\mbox{-}a.s. \quad \mbox{ as }\quad t\to +\infty, \label{eq error a.s. conv to use}
\\
 \text{and }\hskip 2cm&
{t}^{\varrho-1} \int_0^{t}   s^{-\varrho} |\bar \cX_s-\bar \cX_{\us}|^2  \ud s\to 0 \quad \P\mbox{-}a.s. \quad \mbox{ as }\quad t\to +\infty.\hskip1,5cm  \label{eq first term error a.s. conv no speed}
\end{align}

Then, owing to~\eqref{eq error a.s. conv to use} and continuous time C\'esaro principle, see Lemma~\ref{le Cesaro gil} in the Appendix, we obtain  
\begin{align}\label{eq conv 2nd term}
  t^{{\varrho}-1}\int_0^t s^{-(\varrho+1)} \int_0^{s}   |\bar \cX_u-\bar \cX_{\underline u}|^2  du\,\ud s\to 0 \quad \P\mbox{-}a.s. \quad \mbox{ as }\quad t\to +\infty.
\end{align}
Now, we compute 
\begin{align*}
\int_{\g_1}^{+\infty}\frac{1}{s^2}\int_{\us}^s \E\, |\bar \cX_u|^2du\, \ud s& \le C (1+\|\cX_0\|_2)^2\int_{\g_1}^{+\infty}\frac{s-\us}{s^2}\ud s\\
&\le C (1+\|\cX_0\|_2)^2\sum_{n\ge 1}\frac{\g^2_{n+1}}{\Gamma^2_{n}}<+\infty.
\end{align*}
Consequently the r.v. $\displaystyle \int_1^{+\infty}\frac{1}{s^2}\int_{\us}^s |\bar \cX_u|^2du\, \ud s<+\infty$  $\P$-$a.s.$ so that, still by Lemma~\ref{le kronecker gil},
\begin{align}\label{eq conv 3rd term}
  t^{{\varrho}-1}\int_1^{t}s^{-(\varrho+1)} \int_{\us}^s |\bar \cX_u|^2\ud s \to 0   \quad \P\mbox{-}a.s. \quad \mbox{ as }\quad t\to +\infty.
\end{align}
As for the last term, one proceeds likewise 
\begin{align*}
\int_{\g_1}^{+\infty}\frac {1}{s} \frac{s-\us}{s}\frac{1}{\us}\int_0^{\us}\E\, |\bar \cX_{u}|^2du\,\ud s &\le C (1+\|\cX_0\|_2)^2 \int^{+\infty}_{\g_1} \frac{s-\us}{s^2} \,\ud s<+\infty
\end{align*}
so that
\[
\int_{\g_1}^{t}\frac{1}{s}\frac{s-\us}{s}\frac{1}{\us}\int_0^{\us} |\bar \cX_{\us}|^2du\ud s < + \infty 
\]
and one concludes by the Kronecker Lemma~\ref{le kronecker gil}. Combined with~\eqref{eq first term error a.s. conv no speed}-\eqref{eq conv 2nd term}-\eqref{eq conv 3rd term}, this shows that 
\[
  t^{{\varrho}-1}  \int_0^t s^{-{\varrho}} \mathfrak{E}_s \ud s \to 0 \quad \P\mbox{-}a.s.  \quad \mbox{ as }\quad t\to +\infty,
  \]
and, according to~\eqref{eq appli pr stab}, concludes the proof of statement~$(b)$ with $\hat t=t$.  When $\hat t= \t$, the  conclusion sraightforwardly follows using again~\eqref{eq error a.s. conv to use}. 
  \eproof
\color{black}

\subsection{$L^2$-Rate of convergence for the Euler approximation}

\begin{Proposition}[Numerical complexity]\label{pr num comp}
Assume $\hypl2$. Then, for any $\vartheta\in (0,\vtstr)$,
 
  \begin{align} \label{eq control num comp}
   \frac1{\Gamma_n}\!\int_0^{\Gamma_n}  \hskip-0,25cm\esp{|\cX_t-\bar{\cX}_{\hat t}|^2} \!\ud t
     = O_{\vtstr-\vartheta}\left((\Gamma_n)^{-\vartheta}\sum_{k=1}^{n-1} \frac{\gamma_k^2}{\Gamma_k^{1-\vartheta}}\right) \longrightarrow 0 \;\mbox{ as }\; n\to +\infty 
  \end{align}
 as soon as $\sum_{k\ge 1}\frac{\g^2_{k}}{\Gamma_k}<+\infty$ (with still in mind that $\hat t= t$ or $t$).
%
%
%
%
\end{Proposition}

{
  \color{black}
  \proof We first observe that the convergence in \eqref{eq control num comp} is a direct application of Kronecker's Lemma. We now prove the upper bound in \eqref{eq control num comp}.\\
  Step 1. Going back to the proof of Theorem~\ref{thm:tpsdiscretsmulable}, Equation~\eqref{eq temp control disc error} shows that the convergence is controlled via
  \begin{align}\label{eq temp control disc error bis}
    \frac1t\int_0^t \esp{|\bar{\cX}_s -\cX_s|^2} \ud s 
          \le C_{\!\varrho} t^{{\varrho}-1}  \int_0^t s^{-{\varrho}} \mathbb{E}[\mathfrak{E}_s] \ud s ,
  \end{align}
  where 
  \begin{align}\label{eq de main disc error bis}
    \mathfrak{E}_t := |\bar \cX_t-\bar \cX_{\t}|^2 +\frac{1}{t}\int_0^t |\bar \cX_s- \bar \cX_{\us}|^2\ud s +\frac{1}{t}\int_{\t}^t |\bar \cX_s|^2\ud s  
   + \frac{t-\t}{t} \cdot \frac{1}{\t}\int_0^{\t} |\bar \cX_{\us}|^2 \ud s.
  \end{align}
  We then compute, owing to Fubini's Theorem to integrate the second term in $\mathfrak{E}$, that 
  \begin{align*}
    \frac1t\int_0^t \esp{|\bar{\cX}_s -\cX_s|^2} \ud s  \le C t^{{\varrho}-1} 
     \left(\int_0^ts^{-{\varrho}} \esp{|\bar \cX_s-\bar \cX_{\us}|^2} \ud s + \sup_{u \ge 0} \esp{|\bar \cX_{\us}|^2} \int_0^ts^{-(\varrho+1)}(s-\us)\ud s \right).
  \end{align*}
  Invoking Lemma~\ref{lem:MRcont}, we deduce 
  \begin{align}
    \frac1t\int_0^t \esp{|\bar{\cX}_s -\cX_s|^2} \ud s  \le C \cE_t \;\text{ with }\; \cE_t := t^{{\varrho}-1} 
     \int_0^ts^{-\varrho}(s-\us)\ud s .
  \end{align}
 Elementary computations show that
  \begin{align}
    \cE_{\Gamma_n} = O_{\varrho}\left((\Gamma_n)^{\varrho-1}\sum_{k=1}^{n-1} \frac{\gamma_{k+1}^2}{\Gamma_k^\varrho}\right).
  \end{align}
\noindent Step 2. We first observe that 
\begin{align*}
  \esp{\cW^2_2(\nu^{\cX}_{\Gamma_n},\bar{\nu}_{\Gamma_n})} &\le 
  \frac1{\Gamma_n}\int_0^{\Gamma_n}\esp{|\cX_t - \bar{\cX}_{\t}|^2} \ud t
  \\
  & \le 
  \frac 2 {\Gamma_n}\int_0^{\Gamma_n}\esp{|\bar{\cX_t} - \bar{\cX}_{\t}|^2} \ud t
  +\frac 2{\Gamma_n}\int_0^{\Gamma_n}\esp{|\bar{\cX_t} - {\cX}_t|^2} \ud t.
\end{align*} 
Using Lemma~\ref{lem:MRcont}, we then obtain 
\begin{align*}
  \esp{\cW^2_2(\nu^{\cX}_{\Gamma_n},\bar{\nu}_{\Gamma_n})} &\le 
  \frac1{\Gamma_n}\int_0^{\Gamma_n}\esp{|\cX_t - \bar{\cX}_{\t}|^2} \ud t
  \\
  & \le 
  C\frac1{\Gamma_n}\int_0^{\Gamma_n}(t - {\t})\ud t
  +\frac 2{\Gamma_n}\int_0^{\Gamma_n}\esp{|\bar{\cX_t} - {\cX}_t|^2} \ud t
  \\
  &\le C\frac1{\Gamma_n}\int_0^{\Gamma_n}(t - {\t})\ud t + O_\varrho\left((\Gamma_n)^{\varrho-1}\sum_{k=1}^{n-1} \frac{\gamma_{k+1}^2}{\Gamma_k^\varrho}\right),
\end{align*} 
using the result obtained in the first step. Now, similarly compute that 
\begin{align*}
  \frac1{\Gamma_n}\int_0^{\Gamma_n}(t - {\t})\ud t = O\left((\Gamma_n)^{-1}\sum_{k=1}^{n-1}\gamma_{k+1}^2\right).
\end{align*}
Since $(\Gamma_n)_n$ is increasing, we have
\begin{align}
  (\Gamma_n)^{-1}\sum_{k=1}^{n-1} \gamma_{k}^2 \le  (\Gamma_n)^{\varrho-1}\sum_{k=1}^{n-1} \frac{\gamma_{k+1}^2}{\Gamma_k^\varrho} ,
\end{align}
and  we conclude 
\begin{align}\label{eq:majorhofinale}
  \frac1{\Gamma_n}\int_0^{\Gamma_n} \esp{|\bar{\cX}_{\t} -\cX_t|^2} \ud t  = O_\varrho\left((\Gamma_n)^{\varrho-1}\sum_{k=1}^{n-1} \frac{\gamma^2_{k+1}}{\Gamma_k^\varrho}\right).
\end{align}
  \eproof
}

\smallskip
\noindent As a first application of the above proposition, we give the following general bound of convergence.

\begin{Corollary}\label{co non param speed}Let $\hypl2$ hold. 

\noindent  $(a)$  Then 
  \begin{align}\label{eq the condition}
    \left(\varlimsup_n \bigg(\frac{\gamma_n-\gamma_{n+1}}{\gamma_{n+1}^2}\Gamma_{n+1}\bigg)   < \vtstr \right)\Longrightarrow 
    \frac1{\Gamma_n}\int_0^{\Gamma_n} \esp{|\cX_t-\bar{\cX}_{\hat t} |^2} \ud t 
    = O \left(\gamma_n\right).
  \end{align}
$(b)$ More precisely, let $\g_n = \frac{\g_1}{n^r}$, $r\!\in (0,1)$ and  $\vartheta = \frac{r}{1-r}$. Set also $r^{\star}= \frac{\vtstr}{1+\vtstr}$.
\[
 \frac1{\Gamma_n}\int_0^{\Gamma_n} \esp{|\cX_t-\bar{\cX}_{\hat t} |^2} \ud t  \le  \left\{\begin{array}{ll}
 O_{\vtstr-\vartheta}(n^{-r})& \mbox{if } r< r^{\star}\\
  o_{\eta}\big(n^{-(1-r)(1-\frac{\b}{\a})+\eta} \big)& \mbox{if } r\ge r^{\star}
 \end{array}\right.
\]
for every (small) $\eta>0$.
\color{black}
\end{Corollary}

  \color{black}
  \proof $(a)$  Let $\varrho$ be such that   $  \varlimsup_n[\ldots] < 1-\varrho< \vtstr$ in~\eqref{eq the condition}. We will study the upper bound in~\eqref{eq:majorhofinale} for such a $\rho$.
 We set $\tilde{\gamma}_n := \frac{\gamma_n}{\Gamma_n}$ and $\tilde{\Gamma}_n =\sum_{k=1}^n \tilde{\gamma}_k$, for $n \ge 1$ and 
  \begin{align}
    u_n := (\Gamma_n)^{\varrho-1}\sum_{k=1}^{n-1} \frac{\gamma_{k}^2}{\Gamma_k^\varrho}
    \;\text{ and }\; v_n := \frac{u_n}{\gamma_n}.
  \end{align}
  We then compute 
  \begin{align*}
    v_{n+1} = \frac{u_{n+1}}{\gamma_{n+1}} = \lambda_n v_n + \tilde{\gamma}_{n+1} \;\text{ with }\;
    \lambda_n = \left(\frac{\Gamma_n}{\Gamma_{n+1}}\right)^{1-\varrho} \!\!\frac{\gamma_n}{\gamma_{n+1}}
    = (1-\tilde{\gamma}_{n+1})^{1-\varrho}\frac{\gamma_n}{\gamma_{n+1}}.
  \end{align*}
Owing to~\eqref{eq the condition}, we have that for $n \ge n_0$,
\begin{align*}
  \frac{\gamma_n}{\gamma_{n+1}} \le 1 + (1-\varrho - \theta) \tilde{\gamma}_{n+1} \le e^{(1-\varrho - \theta) \tilde{\gamma}_{n+1}}
\end{align*}
for some $\theta > 0$. Hence,
\begin{align}
  v_{n+1} \le e^{-\theta}v_n + \tilde{\gamma}_n \;, n \ge n_0\,.
\end{align}
Iterating the previous inequality, we obtain 
\begin{align}\label{eq main step here}
  v_{n+1}\le e^{-\theta(\tilde{\Gamma}_{n+1} - \tilde{\Gamma}_{n_0})}v_{n_0}
  +
  \sum_{k = n_0 + 1}^{n+1} e^{-\theta(\tilde{\Gamma}_{n+1} - \tilde{\Gamma}_{k})}\tilde{\gamma}_k.
\end{align}
We also observe that 
\begin{align*}
  \sum_{k = n_0 + 1}^{n+1}e^{\theta  \tilde{\Gamma}_{k} } \tilde{\gamma}_k 
  &\le 
  e^{\theta \tilde{\gamma_1}}\sum_{k = n_0 + 1}^{n+1}e^{\theta  \tilde{\Gamma}_{k-1} } \tilde{\gamma}_k 
  \\
  &\le 
  e^{\theta \tilde{\gamma_1}} \sum_{k = n_0 + 1}^{n+1}\int_{k-1}^{k} e^{\theta t} \ud t  \le \frac{e^{\theta \tilde{\gamma_1}}}{\theta}e^{\tilde{\Gamma}_{n+1}\theta}.
\end{align*}
  Combining the previous inequality with~\eqref{eq main step here}, we finally get $ v_{n+1} \le C_{\theta,n_0}$, which yields the announced result. 

\smallskip
\noindent $(b)$ If $\g_n = \frac{\g_1}{n^r}$, one checks that condition~\eqref{eq the condition} reads $r < r^\star$, namely, 
\[
  \varlimsup_n \bigg(\frac{\gamma_n-\gamma_{n+1}}{\gamma_{n+1}^2}\Gamma_{n+1}\bigg) = \vartheta < \vtstr.
  \]
If one sets $r\ge r^\star$. Then, for any $\bar\vartheta < \vtstr$, we have according to \eqref{eq:majorhofinale}
\begin{align*}
  \frac1{\Gamma_n}\int_0^{\Gamma_n} \esp{|\bar{\cX}_{\t} -\cX_t|^2} \ud t & = O_{\bar\vartheta}\left((\Gamma_n)^{-\bar\vartheta}\sum_{k=1}^{n-1} \frac{\gamma^2_{k+1}}{\Gamma_k^{(1-\bar\vartheta)}}\right)\\
  &= O_{\bar \vartheta} \left(n^{-\bar \vartheta (1-r)}  \right)
  \\
  &= O_{\bar \vartheta}  \left(n^{- \vtstr(1-r)+(\vtstr- \bar\vartheta)(1-r)}  \right).
\end{align*}
The fact that $\vtstr- \bar\vartheta$ can be as closed as desired to $0$ leads to the statement. \eproof

\color{black}

\begin{Remark}\label{re best L2 upper bound}
  It follows from the previous result that setting $\gamma_n = \gamma_1 n^{-r^\star}$ allows to obtain the upper bound 
  \begin{align}
    \frac1{\Gamma_n}\int_0^{\Gamma_n} \esp{|\cX_t-\bar{\cX}_{\hat t} |^2} \ud t 
    = o_\eta(n^{-r^\star + \eta})
  \end{align}
  for every (small) $\eta > 0$. This is the best upper bound as for $r>r^\star$, $\vtstr(1-r)< \vtstr(1-r^\star)=r^\star$.
\end{Remark}

\color{black}

\subsection{$A.s.$-Rate of convergence for the Euler approximation}

\begin{Proposition}[$A.s.$ convergence rate]\label{prop: as-vitesse}
$(a)$ Assume that $\hypas$ is in force.
 Then,  for every $\vartheta \!\in \big(0, \vtstr\big)$  such that 
 \begin{equation}\label{eq:serie-mu}
\sum_{k\ge1} \frac{\g^2_{k}}{\Gamma_k^{1-\vartheta\wedge \frac 12}(\log \Gamma_k)^{1+\epsilon}}<+\infty\quad \mbox{ for some $\epsilon>0$}, 
\end{equation}
one has, for every (small) $\eta>0$, 
\begin{align*}
{\cal W}_2^2(\bar \nu_t, \nu^{\cX}_t \big) \le \frac1t \int_0^t |\bar \cX_{\underline s}-\cX_s|^2 \ud s& 
=  o_{\vartheta,\eta} \big(t^{-(\vartheta\wedge\frac12)}(\log t)^{\frac 12 +\eta}\big)
\;\; \P\mbox{-}a.s.
\end{align*}

\noindent $(b)$ Assume  $\g_n = \frac{\g_1}{n^r}$ for some $\g_1>0$, $r\!\in (0,1)$.  
Let $r^{\star}:=\frac{\vtstr}{1+\vtstr}\!\in (0,\frac 12] $. For every $\eta>0$, 
\[
{\cal W}_2\big(\bar \nu_{_{\Gamma_n}}, \nu^{\cX}_{_{\Gamma_n}} \big)
=\left\{\begin{array}{ll}
o_\eta\big(\Gamma_n^{-\frac{r}{2(1-r)}}(\log \Gamma_n)^{\frac 14 +\eta}\big)= o_\eta\big(n^{-\frac{r}{2}}(\log  n)^{\frac 14 +\eta}\big)&\mbox{ if }\; \vtstr> \frac 12 \mbox{ and }r< 1/3\\
o_\eta(\Gamma_n^{-\frac 14}(\log \Gamma_n)^{\frac 14 +\eta}\big)= o_\eta(n^{-\frac{1-r}{4}}(\log  n)^{\frac 14 +\eta}\big)&\mbox{ if }\; \vtstr > \frac12 \mbox{ and } r\ge 1/3 \\
 o_{\eta}\big( \Gamma_n ^{-\frac{r}{2(1-r)}}(\log \Gamma_n)^{\frac14 +\eta}\big)
= o_{\eta}\big(n^{-\frac r2} (\log n)^{\frac 14+\eta}\big)&\mbox{ if }\; \vtstr \le  \frac 12 \mbox{ and }r< r^{\star}\\
o_{\eta}\big( \Gamma_n^{-\frac{\a-\b}{2\a}}(\log \Gamma_n)^{\frac 14+\eta}\big)
= o_{\eta}\big(n^{(1-r)\frac{\a-\b}{2\a}}(\log n)^{\frac 14+\eta}\big)& \mbox{ if }\;  \vtstr \le \frac 12\ \mbox{ and }r\ge r^{\star}.
\end{array}\right. 
\]
\end{Proposition} 
%

\medskip
\color{black}
\proof 
We assume that $\cX_0 \in \mathscr{L}^4(\cF_0,\P)$: The proof for the case $\cX_0 \in \mathscr{L}^2(\cF_0,\P)$ only is obtained following the same arguments as in step 2.(ii) of the proof of Theorem \ref{th conv emp sid to emp mkv no rate}. 
$(a)$ First note that 
\[
 \frac1t \int_0^t |\bar \cX_{\underline s}-\cX_s|^2 \ud s \le 2\Big(\frac1t \int_0^t |\bar \cX_{s}-\cX_s|^2\ud s + \frac1t \int_0^t |\bar \cX_{s}-\bar \cX_{\underline s}|^2\ud s\Big) .
\]
We start from the upper-bounds~\eqref{eq appli pr stab} and~\eqref{eq majo conv without a rate} which reads: for every $\varrho\!\in \big(\frac{\b}{\a}, 1-\vartheta)$,
\[
\frac1t \int_0^t |\bar \cX_{s}-\cX_s|^2 \ud s \le C_{\varrho} t^{\varrho-1}\int_0^t s^{-\varrho}\mathfrak{E}_s\ud s  +o_{\eta,\varrho}(t^{-(\frac 12\wedge(1-\varrho))}(\log t)^{\frac 12 +\eta})\quad \mbox{for every $\eta>0$}
\]
where $\mathfrak{E}_t $ is given by~\eqref{eq de main disc error} and we used that $t^{-(\frac 12\wedge (1-\varrho))}= o(t^{-(\frac 12\wedge \vartheta)})$. This follows from   $\HYP{MV}_{2, K',\a',\b'}$  and $\cX_0\!\in L^4(\P)$ (see the proof of Theorem~\ref{thm:tpsdiscretsmulable}). Let us denote $ \bar{\mathfrak{E}}^{(\varrho)}_t := t^{\varrho-1}\int_0^t s^{-\varrho}\mathfrak{E}_s\ud s$. Then one has 
 
\begin{align*}
 \bar{\mathfrak{E}}^{(\varrho)}_t&= t^{\varrho-1} \int_0^t s^{-\varrho}|\bar \cX_s-\bar \cX_{\underline s}|^2 \ud s +  t^{\varrho-1} \int_0^t s^{-(1+\varrho)}\int_0^s u^{-\varrho}|\bar \cX_u-\bar \cX_{\underline u}|^2 \ud u\ud s \\
 					    & \quad + t^{\varrho-1} \int_0^t s^{-(\varrho+1)}\int_{\underline s}^su^{-\varrho}|\bar \cX_u|^2 \ud u \ud s+ t^{\varrho-1} \int_0^t s^{-(\varrho+1)}\frac{s-\underline s}{s\underline s}\int_{0}^{\underline s} u^{-\varrho}|\bar \cX_{\underline u}|^2 \ud u \ud s,
\end{align*}
with the usual convention that $\frac 10 \int_0^0 \ud u=0$. 
Since 
\[
\frac 1t \int_0^t |\bar \cX_s-\bar \cX_{\underline s}|^2\ud s\le  t^{\varrho-1} \int_0^t s^{-\varrho}|\bar \cX_s-\bar \cX_{\underline s}|^2\ud s,
\]
we obtain 
\begin{align}\label{eq as to conclude}
  \frac1t \int_0^t |\bar \cX_{\underline s}-\cX_s|^2 \ud s \le 4\bar{\mathfrak{E}}^{(\varrho)}_t +  +o_{\eta,\varrho}(t^{-(\frac 12\wedge(1-\varrho))}(\log t)^{\frac 12 +\eta})\quad \mbox{for every $\eta>0$}.
\end{align}
We now turn to the study of $\bar{\mathfrak{E}}^{(\varrho)}$.
We have only three terms to inspect since
\[
 \int_0^t s^{-(1+\varrho)}\int_0^s u^{-\varrho}|\bar \cX_u-\bar \cX_{\underline u}|^2 \ud u\ud s  \le   \tfrac{1}{\varrho} \int_0^t s^{-\varrho}|\bar \cX_s-\bar \cX_{\underline s}|^2 \ud s
\] 
owing to Fubini's Theorem.   First, as
\[
 \int_{\Gamma_1}^{+\infty} \frac{s^{-\varrho}\E\, |\bar \cX_s-\bar \cX_{\underline s}|^2 }{s^{1-(\vartheta\wedge\frac 12)-\varrho}(\log s)^{1+\epsilon}}\ud s\le  C \int_{\Gamma_1}^{+\infty} \frac{s-\underline s}{s^{1-(\vartheta\wedge\frac 12)}(\log s)^{1+\epsilon}} \ud s\le C \sum_{k\ge1} \frac{\g_{k+1}^2}{\Gamma^{1-(\vartheta\wedge\frac 12)}_{k}(\log \Gamma_k)^{1+\epsilon}}<+\infty,
\]
owing to\eqref{eq:serie-mu}, Kronecker's lemma implies since 
$(\vartheta\wedge \frac 12)+\varrho \le \vartheta +\varrho < 1$ so that $t^{1-\varrho-(\vartheta\wedge \frac12)} \uparrow +\infty$ as $t\to +\infty$, 
\[
t^{\varrho-1} \int_0^t s^{-\varrho}|\bar \cX_s-\bar \cX_{\underline s}|^2 \ud s= o_{\varrho}\big(t^{-(\vartheta\wedge \frac 12)} \big) \quad \P\mbox{-}a.s.
\]
As for the third   term in the r.h.s. of the inequality
\begin{align*}
\int_1^{+\infty} \frac{s^{-(1+\varrho)}}{s^{1-(\vartheta\wedge \frac12)-\varrho}(\log s)^{1+\epsilon}} \int_{\underline s}^su^{-\varrho}\E\, |\bar \cX_u|^2 \ud u \ud s
&\le C_{\varrho} \int_1^{+\infty} \frac{s^{-(1+\varrho)}}{s^{1-(\vartheta\wedge \frac12)}(\log s)^{1+\epsilon}} (s-\underline s)\ud s\\
&\le C_{\varrho} \int_1^{+\infty} \frac{ s-\underline s}{s^{1-(\vartheta\wedge \frac12)}(\log s)^{1+\epsilon}}  \ud s<+\infty
\end{align*}
and one concludes as above that this third  term is $o_{\varrho}(t^{-(\vartheta\wedge\frac12)})$  $\P\mbox{-}a.s$ owing to Kronecker's Lemma.
Finally, for the fourth term, one has
\begin{align*}
 \int_1^t \frac{s^{-(\varrho+1)}}{s^{1-(\vartheta\wedge \frac 12))+\varrho}(\log s)^{1+\epsilon}} \frac{s-\underline s}{s\underline s}\int_{1}^{\underline s} u^{-\varrho}\E\, |\bar \cX_{\underline u}|^2 \ud u \ud s &\le C_{\varrho}\int_1^{+\infty} s^{-2(1+\varrho)}\frac{s-\underline s}{s^{1-(\vartheta\wedge \frac12)-\varrho}(\log s)^{1+\epsilon}}\ud s\\
 &\le  C_{\varrho}\int_1^{+\infty} \frac{s-\underline s}{s^{1-(\vartheta\wedge \frac12)}(\log s)^{1+\epsilon}}\ud s<+\infty 
\end{align*}
which leads to the same estimate after applying Kronecker's Lemma. 
The proof for this step is concluded by inserting back the previous estimate into \eqref{eq as to conclude}.

\noindent $(b)$  Let us inspect the four cases.

\smallskip  
\noindent  $\blacktriangleright$ $\vtstr >\frac 12$:

\smallskip
-- If $0<r<\frac 13$, set $\vartheta= \frac{r}{1-r}\!\in (0, \frac 12)$. Then $2r+ (1-\vartheta)(1-r)= 1$ so that~\eqref{eq:serie-mu} is satisfied  and ${\cal W}^2_2\big(\bar \nu_{_{\Gamma_n}}, \nu^{\cX}_{_{\Gamma_n}} \big)= o_{\eta}\big( \Gamma_n ^{-\frac{r}{1-r}}(\log \Gamma_n)^{\frac12 +\eta}\big)=  o_{\eta}\big(\g_n |\log(1/\g_n)|^{\frac12+\eta}\big)= o_{\eta}\big(n^{-r} |\log(1/\g_n)|^{\frac12+\eta}|\big)$.

\smallskip
--  If $\frac13\le r< 1$, set $\vartheta = \frac 12$. Then $2r+ (1-\frac 12)(1-r)\ge 1$ so that~\eqref{eq:serie-mu} is satisfied and  ${\cal W}^2_2\big(\bar \nu_{_{\Gamma_n}}, \nu^{\cX}_{_{\Gamma_n}} \big)= o_{\eta}\big( \Gamma_n ^{-\frac{1}{2}}(\log \Gamma_n)^{\frac12 +\eta}\big)=  o_{\eta}\big(\g^{\frac{1-r}{2r}}_n |\log(1/\g_n)|^{\frac12+\eta}\big)= o_{\eta}\big(n^{\frac{1-r}{2}} |\log(1/\g_n)|^{\frac12+\eta}|\big)$;

\smallskip  
\noindent  $\blacktriangleright$ $\vtstr \le \frac 12$: 

\smallskip
--  If $0<r<r^{\star}= \frac{\vtstr}{1 + \vtstr}$, set 
$\vartheta= \frac{r}{1-r}<\vtstr$. Then $2r +(1-\vartheta)(1-r)> 1$ so that so that~\eqref{eq:serie-mu} is satisfied  and ${\cal W}^2_2\big(\bar \nu_{_{\Gamma_n}}, \nu^{\cX}_{_{\Gamma_n}} \big)= o_{\eta}\big( \Gamma_n ^{-\frac{r}{1-r}}(\log \Gamma_n)^{\frac12 +\eta}\big)=  o_{\eta}\big(\g_n |\log(1/\g_n)|^{1+\eta}\big)= o_{\eta}\big(n^{-r} (\log n)^{1+\eta}|\big)$.

\smallskip
 -- If $r\ge r^{\star}$, set $\bar \vartheta< \vtstr \le \frac 12 \le \vartheta = \frac{r}{1-r}$. Then $2r +(1-\bar \vartheta)(1-r)= 1 + \frac{\vartheta - \bar \vartheta}{1+\vartheta}> 1$
  so that~\eqref{eq:serie-mu} is satisfied . 
 Then ${\cal W}^2_2\big(\bar \nu_{_{\Gamma_n}}, \nu^{\cX}_{_{\Gamma_n}} \big)= o_{\eta}\big( \Gamma_n^{-\bar \vartheta}|\log \Gamma_n|^{\frac 12+\eta}\big)=  o_{\eta}\big( \gamma_n^{-(\frac 1r -1)\bar \vartheta}|\log \gamma_n|^{\frac 12+\eta}\big)= o_{\eta}\big(n^{(1-r)\bar \vartheta}(\log n)^{\frac 12+\eta}\big)$. 
\eproof
\color{black}

\bigskip
\begin{Remark}\label{re conv rate a.s. case}
  We observe from the previous upper bounds for ${\cal W}_2\big(\bar \nu_{_{\Gamma_n}}, \nu^{\cX}_{_{\Gamma_n}} \big)$, that:
  \\
 - if $\vtstr \le \frac12$, then it is optimal to set $\gamma_n = \gamma_1 n^{-r^\star}$ to obtain: 
 \[
  {\cal W}_2\big(\bar \nu_{_{\Gamma_n}}, \nu^{\cX}_{_{\Gamma_n}} \big)
  =o_\eta\left(n^{-\frac{r^\star}2 + \eta}\right)
  \]
  for any small $\eta$.
 \\
 - if $\vtstr>\frac12$, then it is optimal to set $\gamma_n = \gamma_1n^{-\frac13}$ to obtain 
 \[
  {\cal W}_2\big(\bar \nu_{_{\Gamma_n}}, \nu^{\cX}_{_{\Gamma_n}} \big)
  =o_\eta\left(n^{-\frac16}\log(n)^{\frac14+\eta}\right)
  \]
  for any small $\eta$. In this case, the convergence is limitated by the martingale term and stuck to the optimal $r^\star=\frac13$ associated to $\vtstr = \frac12$ (up to log terms).  
\end{Remark}
\color{black}


\subsection{Convergence to the stationary measure}
This section is dedicated to the proof of our main theorems announced in the introduction. 
The results, given in Section \ref{se stationnary solution MKVSDE}, Section \ref{se SID} and the current section, use the mean-reversion hypothesis formulated for the McKean-Vlasov dynamics. Our main theorems given in the introductory session make use of assumptions on $\sigma$. The next technical Lemma, whose proof is left to the reader, clarifies the relationship between these two sets of assumptions. 

\begin{Lemma}[From confluence to  mean-reversion] 
    \label{le structural assumptions3}
     Let $p\ge 1$.  Assume~\HYP{L} is in force.
     { \color{black}  
    \noindent $(a)$ Assume  and $b$ and $\s$ also satisfy  $\HYP{C}_{p, \a,\b}$. Then, the following holds:
    \begin{enumerate}[(i)]
    \color{black}
    \item for every $\bar \a\!\in (0, \a)$, there exists $\bar p>p$, $\bar K\ge 0$ and  $ \bar \b \!\in( 0,\bar \a)$, such that,  for every $p'\!\in [p,\bar p]$,    $\HYP{C}_{p', \bar \a,\bar \b}$ is satisfied and, if $\b <\a$, for every $\bar \a \!\in (\b, \a)$, one may choose $\bar \b <\bar \a$.
 \item  for every  $\bar \a\!\in (0, \a)$, there exist  $\bar p>p$,  $\bar{K},  \bar \b \ge 0,$ such that    for every $p'\in [p, \bar p]$, $\HYP{MV}_{p', \bar K, \bar \a, \bar \b}$ is satisfied. If $\b <\a$, one may also choose $\bar \b < \bar \a$. 
 \item if, furthermore, $\HYP{\sigma}$ holds then  for every $\bar p>p$ and every $ \bar \a\!\in (0, \a)$, there exist   $\bar K, \bar \b \ge 0$   such that  $\HYP{MV}_{\bar p, \bar K,  \bar \a, \bar \b}$ is satisfied and for every $ \bar \a\in (\b, \a)$, one may choose   $\bar \b< \bar \a$. 
\end{enumerate}
 
\noindent $(b)$ If  $b$ and $\s$ satisfy $\HYP{MV}_{p,K, \a, \b}$, then   for every $\bar \a\!\in (0, \a)$, there exists $\bar p>p$, $\bar K, \bar \b \ge 0$ such that, for every $p'\!\in [p,\bar p]$,     $\HYP{MV}_{p',\bar K, \bar \a,\bar \b}$  is satisfied and if, furthermore, $\b <\a$,  for every $\bar \a\in (\b, \a)$, one may  choose $\bar \b <\bar \a$. 
     }
\end{Lemma}

\subsubsection{Proof of Theorem \ref{th main intro no rate} }
We simply observe that:
\[ \cW_2(\nstr,\bar{\nu}_{\Gamma_n}) \le \cW_2(\nstr,\nu^\cX_{\Gamma_n}) + \cW_2(\nu^\cX_{\Gamma_n},\bar{\nu}_{\Gamma_n}).  \]
\noindent $\blacktriangleright$ 
The first statement is obtained by combining Theorem \ref{th conv emp sid to emp mkv no rate} (a) with Corollary \ref{co conv scheme no rate} (a).

\medskip
\noindent $\blacktriangleright$ For the second statement, we combine Theorem \ref{th conv emp sid to emp mkv no rate} (b) with  Corollary \ref{co conv scheme no rate} (b).

\subsubsection{Proof of Theorem \ref{th main results intro L2 rate} }
It follows from Minkowski inequality that 
\begin{align}\label{eq starting point proof th1.2}
\esp{\cW^2_2(\nstr,\bar{\nu}_{\Gamma_n})}^\frac12 
\le   \esp{\cW^2_2(\nstr,\nu^\cX_{\Gamma_n})}^\frac12  + \esp{\cW^2_2(\nu^\cX_{\Gamma_n},\bar{\nu}_{\Gamma_n})}^\frac12 .
\end{align}
From Theorem \ref{th l2 conv with a rate sid to nu star}, setting $\eta=(\vtstr-\vartheta)(1-r^\star)$, we obtain
\begin{align*}
  \esp{\cW^2_2(\nu^{\cX}_{\Gamma_n},\nu^\star)} \le 
  o_\eta \left({n}^{-\vtstr(1-r^\star)+\eta}\esp{|X^\star_0-\cX_0|^2}\right) +
  O_\eta \left( t^{-\zeta}(\log t)^{\kappa}\right).
\end{align*}
Then one can take the infinimum on the coupling $(X_0^\star,\cX_0)$ to obtain 
\begin{align*}
  \esp{\cW^2_2(\nu^{\cX}_{\Gamma_n},\nu^\star)} \le 
  o_\eta \left({n}^{-\vtstr(1-r^\star)+\eta} \cW_2^2(\nstr,[\cX_0])\right) +
  O_\eta \left( t^{-\zeta}(\log t)^{\kappa}\right).
\end{align*}
For (i) and (iii),The proof is concluded by inserting the previous inequality into \eqref{eq starting point proof th1.2} and invoking Remark \ref{re best L2 upper bound}. 
\\
Statement (ii) is obtained by letting $p^\star$ to be large enough so that $\zeta_{p^\star}$ is close to $\frac1{2(d+3)}$. Indeed, from Lemma \ref{le structural assumptions3}, we know that $\HYP{MV}_{\bar p, \bar K,  \bar \a, \bar \b}$ is satisfied for any $p\ge 1$ under $\HYP{\sigma}$. Then, according to Corollary \ref{co main res invariant distrib mkv}(b), one can take $p^\star$ as large as required.

\subsubsection{Proof of Theorem \ref{th main results intro L2 as} }
We simply observe that:
\[ \cW_2(\nstr,\bar{\nu}_{\Gamma_n}) \le \cW_2(\nstr,\nu^\cX_{\Gamma_n}) + \cW_2(\nu^\cX_{\Gamma_n},\bar{\nu}_{\Gamma_n}) . \]
For statement (i) and (iii) we combine Theorem \ref{th conv emp sid to emp mkv with a rate} and Remark \ref{re conv rate a.s. case} with the previous inequality. Statement (ii) follows from the same arguments  as the ones used in the proof of Theorem \ref{th main results intro L2 rate}.

\subsection{Proof of Lemma \ref{prop:cXLpbounded}}

\label{se proof of boundedness marginal scheme}

We conclude this section with the proof of the key Lemma \ref{prop:cXLpbounded}.

\medskip

\color{black}
\proof
%
We set for notational convenience
\[
\bar \cX_{k}:= \bar \cX_{\Gamma_k}, \quad b_k := b(\bar \cX_{\Gamma_k},\bar{\nu}_{\Gamma_k}) \quad \mbox{ and }\quad \s_k = \s(\bar \cX_{\Gamma_k},\bar{\nu}_{\Gamma_k}),\quad k\ge 0.
\]
Rewriting~\eqref{eq:SDEdisc2} with these notations, we have
\[
\bar \cX_{k+1} = \bar \cX_k +\g_{k+1} b_k +\sqrt{\g_k}\s_k Z_{k+1}.
\]
We also denote $\bar \nu_k= \bar \nu_{\Gamma_k}$. We will prove on our way by induction that $\bar \cX_k\!\in L^{2p}(\P)$ if $\bar \cX_0\!\in L^{2p}(\P)$.

 A second order Taylor expansion to the function $x \mapsto \cV_p(x)$, recall \eqref{de eq a lyapunov function}, at $ \bar \cX_k $  
yields 
\begin{align*}
 | \bar \cX_{k+1}|^{2p} &= |\bar \cX_k|^{2p} +2p |\bar \cX_k|^{2(p-1)}\big( \bar \cX_k\,|\,\Delta \bar \cX_{k+1} \big)+\tfrac 12 \left(\partial^2_{xx}\cV_p(\xi_{k+\frac 12})|\Delta\bar \cX_{k+1}(\Delta\bar \cX_{k+1})^{\top}\right)_F
\end{align*}
where $\xi_{k+\frac 12}$ belongs to the geometric interval $(\bar \cX_{k}, \bar \cX_{k+1})$ and $\Delta\bar \cX_{k+1}= \bar \cX_{k+1}-\bar \cX_{k}$. 
According to \eqref{eq useful majo hess}, we get
\begin{align}\label{eq:L^pbound1}
 | \bar \cX_{k+1}|^{2p} &\le |\bar \cX_k|^{2p} +p\Big(2 |\bar \cX_k|^{2(p-1)} \big(\bar \cX_k\,|\,\Delta \bar \cX_{k+1} \big)+(2p-1)|\xi_{k+\frac 12}|^{2(p-1)}|\Delta\bar \cX_{k+1}|^{2}\Big).
\end{align}
\noindent In what follows we consider the filtration ${\cal G}_k= {\cal F}^{\cX_0,W}_{\Gamma_k}$, $k\ge 0$ and  $\E_k$ will denote the conditioning w.r.t. ${\cal G}_k$. 

 \noindent {\em Case $p=1$}. Assume  $\bar \cX_\ell\!\in L^2(\P)$, $\ell=0,\ldots,k$, then $\E_k(b_k\,|\,\Delta W_{t_{k+1}})=0$ and $\E_k |\sigma_k\Delta W_{t_{k+1}}|^2= \|\sigma_k\|_{_F}^2$.   Taking the conditional expectation, denoted $\E_k$, with respect to the filtration ${\cal G}_k= {\cal F}^{X_0,W}_{\Gamma_k}$,
 \begin{align*}
\E_k | \bar \cX_{k+1}|^{2} &\le |\bar \cX_k|^{2} + \g_{k+1}\Big(2\big( \bar \cX_k\,|\,b_k \big)+\|\s_k\|_{_F}^{2}+2\g_{k+1}|b_k|^2\Big)\\
&\le |\bar \cX_k|^{2} + \g_{k+1}\Big(\bar K- \bar \a |\bar \cX_k|^2 +\bar \b {\cal W}_2^2(\bar \nu_k,\d_0)+2\g_{k+1}|b_k|^2 \Big),
\end{align*}
where we used~\eqref{eq hyp mv} in the last line.  The function $b$ being Lipschitz continuous under \HYP{L}, one has 
\begin{equation}\label{eq:boundb}
|b_k|^2 \le C_b(|\bar \cX_k|^2+{\cal W}_2^2(\bar \nu_k, \d_0)\big)
\end{equation}
for some real constant $C_b$. Let $\ve\!\in (0, \frac{\bar \a- \bar \b}{2})$ and $k_{\ve}$ such that $2C_b \g_{k+1}<\ve$ and $1- \g_{k+1}(\bar \a-\ve) >0$ for every $k\ge k_{\ve}$. Inserting, for such $k$, the above bound for $b_k$ into  the above inequality implies
 \begin{align}
\label{eq:L2bound1}\E_k | \bar \cX_{k+1}|^{2} &\le   |\bar \cX_k|^{2} +\g_{k+1}\Big(K-(\bar \a-\ve) |\bar \cX_k|^2 +(\bar \b+\ve) {\cal W}_2^2(\bar \nu_k,\d_0)\Big),\\
\label{eq:L2bound2}&\le |\bar \cX_k|^2 \big(1-\g_{k+1}(\bar \a-\ve))\big)+ \g_{k+1} (\bar \b +\ve)\Big(\frac{1}{\Gamma_n}\sum_{\ell=0}^{k-1}\g_{\ell+1}|\bar \cX_\ell|^2\Big) + \bar K\g_{k+1}.
\end{align}
Consequently taking this time regular expectation, one gets
 \begin{align*}
\E | \bar \cX_{k+1}|^{2} &\le \E\, |\bar \cX_k|^2 \big(1-\g_{k+1}(\bar \a-\ve))\big)+ \g_{k+1} (\b +\ve)\Big(\frac{1}{\Gamma_k}\sum_{\ell=0}^{k-1}\g_{\ell+1}\E\, |\bar \cX_\ell|^2\Big) +K\g_{k+1},\\
& \le  \E\, |\bar \cX_k|^2 \big(1-\g_{k+1}(\bar \a-\ve))\big)+ \g_{k+1} (\bar \b +\ve)\max_{\ell=0,\ldots,k-1}\E\, |\bar \cX_\ell|^2 +\bar K\g_{k+1}<+\infty.
\end{align*}

\noindent This proves that $\bar \cX_{k+1}\!\in L^2(\P)$.

 \noindent Now let $A\ge \max_{\ell=0,\ldots,k_{\ve}}\E\, |\bar \cX_{\ell}|^2/(1+\E\,| \bar \cX_0|^{2p})$ and set $\tau_{_A}= \min\big\{k : \E\, |\bar \cX_{k}|^2 >A(1+\E\,| \bar \cX_0|^{2p}) \big\}\!\in \{k_\ve+1,\ldots\}\cup \{\infty\}$. If $\tau_{_A}<+\infty$, then
 \begin{align*}
 A < \frac{\E\, |\bar \cX_{\tau_{_A}}|^2}{(1+\E\,| \bar \cX_0|^{2p})} \le A - A\g_{\tau_{_A}}(\bar \a- \bar \b-2\ve) +\bar K \g_{\tau_{_A}}.
 \end{align*}
 This implies that  $2\, A (\bar \b -\bar \a+2\, \ve) +\bar K>0$ which is impossible for $A>\frac{\bar K}{\bar \a-\bar \b-2\ve}$ since $\bar \a-\bar \b-2\ve>0$. For such $A$, $\tau_{_A}=+\infty$ which proves the announced boundedness result.
 
 \medskip
 \noindent {\em Case $1<p\le 3/2$}. In that case  $2(p-1)\!\in(0,1]$ and the function $u\mapsto u^{2(p-1)}$ is $2(p-1)$-H\" older so that, for every $u$, $v\ge 0$,
 \begin{equation}\label{eq:Holder2(p-1)}
(u+v)^{2(p-1)} \le u^{2(p-1)} + v^{2(p-1)}.
\end{equation}
Assume  $\bar \cX_\ell\!\in L^2(\P)$, $\ell=0,\ldots,k$. We come back to~\eqref{eq:L^pbound1}. Applying \eqref{eq:Holder2(p-1)} to $\xi_{k+\frac12}= \cX_k + \theta  \Delta \bar \cX_{k+1}$ (with$(0,1)$-valued  $\theta$) yields, once multiplied by $|\Delta \bar \cX_{k+1}|^2$, 
 \begin{align*}
 |\xi_{k+\frac12}|^{2(p-1)} |\Delta \bar \cX_{k+1}|^2&\le |\bar \cX_{k}|^{2(p-1)}| \Delta \bar \cX_{k+1}|^2 +  | \Delta \bar \cX_{k+1}|^{2p} .
  \end{align*}
Plugging this inequality  into~\eqref{eq:L^pbound1}, one has
 \begin{align*}
   | \bar \cX_{k+1}|^{2p} &\le |\bar \cX_k|^{2p} +\g_{k+1} p|\bar \cX_k|^{2(p-1)}\Big( 2\big(\bar \cX_k\,|\,\Delta \bar \cX_{k+1} \big)+(2p-1)| \Delta \bar \cX_{k+1}|^2 \Big) + p(2p-1)  | \Delta \bar \cX_{k+1}|^{2p}.
   \end{align*}
Thus conditioning by $\E_k$ yields, still using that $\E_k(b_k\,|\,\Delta W_{t_{k+1}})=0$ and $\E_k |\sigma_k\Delta W_{t_{k+1}}|^2= \gamma_{k+1}\|\sigma_k\|_{_F}^2$ entails
 \begin{align*}
\E_k | \bar \cX_{k+1}|^{2p} &\le  |\bar \cX_k|^{2p} +  \g_{k+1}p|\bar \cX_k|^{2(p-1)} \Big( 2\big( \bar \cX_k\,|\,b_k \big)+(2p-1)\|\s_k\|_{_F}^{2}+2\g_{k+1}|b_k|^2\big) +2p \E_k\,| \Delta \bar \cX_{k+1}|^{2p}.
\end{align*}
Now, note that,   there exists a positive real constant $ C= C_{b,\s,p}$ such that
\[
\E_k\, |\Delta \bar \cX_{k+1}|^{2p} \le   C \g_{k+1}^{p} |\bar \cX_k|^p.
\]
Then, there exists a positive real constant $C= C_{b,\s,p,\bar \g}$ (with $\bar \g= \sup_n \g_n$) such that 
\begin{align*}
\E_k | \bar \cX_{k+1}|^{2p} &\le |\bar \cX_k|^{2p} +  \g_{k+1}  p |\bar \cX_k|^{2(p-1)} \Big( 2\big( \bar \cX_k\,|\,b_k \big)+(2p-1)\|\s_k\|_{_F}^{2}+\g_{k+1}|b_k|^2\Big) +2p\,C \g^{p}_{k+1}|\bar \cX_k|^{2p}\\
&\le |\bar \cX_k|^{2p} \\
&+ \g_{k+1} p|\bar \cX_k|^{2(p-1)} \Big(\bar K-(\bar \a-2C \g_{k+1}^{p-1}-2C_b\g_{k+1}) |\bar \cX_k|^2 +(\bar \b+2C_b\g_{k+1}){\cal W}_2^2(\bar \nu_k,\d_0) \Big).
\end{align*}
At this stage, we mimic what we  did when $p=1$ to get~\eqref{eq:L2bound2}. Let $\ve\!\in (0, \frac{\bar \a- \bar \b}{2}) $ be small enough so that $1-p(\bar \a-\ve)>0$ and let $k_{\ve}\!\in \mathbb{N}$ be such that, for every $k\ge k_{\ve}$, 
$C \g_{k+1}^{p-1} + C_b\g_{k+1}<\ve/2$. Then 
\begin{align*}
\E_k | \bar \cX_{k+1}|^{2p} &\le |\bar \cX_k|^{2p} + p |\bar \cX_k|^{2(p-1)} \Big(\bar K-(\bar \a-\ve) |\bar \cX_k|^2 +(\bar \b+\ve) {\cal W}_2^2(\bar \nu_k,\d_0) \Big).
\end{align*}
Rearranging the terms, using the standard identity ${\cal W}_2^2(\mu,\d_0)= \int|\xi|^2\mu(\ud \xi)$  and  Young's inequality with conjugate exponents $p$ and $\frac{p}{p-1}$ yields 
\begin{align*}
\E_k | \bar \cX_{k+1}|^{2p} & \le  |\bar \cX_k|^{2p}  \big(1-p(\bar \a-\ve) \big)+p (\bar \b+\ve) \frac{1}{\Gamma_k} \sum_{\ell=1}^{k-1}\g_{\ell}|\bar \cX_{k}|^{2(p-1)}|\bar \cX_{\ell}|^{2}+ p \bar K \g^p_{k+1} |\bar \cX_k|^{2(p-1)} \\
&\le   |\bar \cX_k|^{2p}  \big(1-p(\bar \a-\ve) \big)+
p (\b+\ve) \frac{1}{\Gamma_k} \sum_{\ell=1}^{k-1}\g_{\ell}\Big(\frac{(p-1)|\bar \cX_{k}|^{2p}}{p}+\frac{|\bar \cX_{\ell}|^{2p}}{p}\Big)\\
&\qquad +  pK \g^p_{k+1}|\bar \cX_k|^{2(p-1)}.
\end{align*}
Taking expectations yields
\begin{align*}
\E |\bar \cX_{k+1}|^{2p} & \le  \E\,  |\bar \cX_k|^{2p}  \big(1-p(\bar \a-\ve) \big)\\
& \qquad +p (\bar \b+\ve) 
 \frac{1}{\Gamma_k} \sum_{\ell=1}^{k-1}\g_{\ell}\Big(\frac{(p-1)\E\, |\bar \cX_{k}|^{2p}}{p}+\frac{\E\,|\bar \cX_{\ell}|^{2p}}{p}\Big)+  pK\g_{k+1}\E\, |\bar \cX_k|^{2(p-1)}\\
 &\le   \E\,  |\bar \cX_k|^{2p}  \big(1-p(\bar \a-\ve) \big)+p (\bar \b+\ve) \max_{\ell=0,\ldots,k}\E\, |\bar \cX_{\ell}|^{2p}+ p\bar K\g_{k+1}[\E\, |\bar \cX_k|^{2p}]^{\frac{p-1}{p}},
 \end{align*}
where we used in the penultimate line   Young's inequality with conjugate coefficients $p$ and $\frac{p}{p-1}$. 

\noindent This clearly shows that $\bar \cX_{k+1}\!\in L^{2p}(\P)$. Hence $\bar \cX_n\!\in L^{2p}(\P)$ for every $n\ge 0$.

\noindent Let $A$ be such that $ A\ge \max_{\ell=0,\ldots,k_{\ve}}\E\, |\bar \cX_\ell|^{2p}/(1+\E\,| \bar \cX_0|^{2p})$. We define 
$$
\tau_{_A}= \min\big\{ k: \E\, |\bar \cX_k|^{2p}>A(1+\E\,| \bar \cX_0|^{2p})\big\}\ge k_{\ve}.
$$ 
It is clear that $A\mapsto \tau_{_A}$ is non-decreasing and goes to infinity as $A\uparrow +\infty$  since $\E\, |\bar \cX_n|^{2p}<+\infty$ for every $n\ge 0$.

\noindent Assume all $\tau_{_A}$ are  finite. Then for every $A>0$, one has
\[
A  < \frac{\E\, |\bar \cX_{\tau_{_A}}|^{2p}}{1+\E\,| \bar \cX_0|^{2p}}\le A -pA(\bar \a- \bar \b -2\ve)  +p \bar K \g_{\tau_{_A} }  \frac{A^{\frac{p-1}{p}}}{(1+\E\,| \bar \cX_0|^{2p})^{\frac 1p} }
\]
so that
\[
A<\frac{\E\, |\bar \cX_{\tau_{_A}}|^{2p}}{1+\E\,| \bar \cX_0|^{2p}}\le A -pA(\bar \a- \bar \b -2\ve)  +p \bar K \g_{\tau_{_A}} A^{\frac{p-1}{p}} 
\]
which in turn implies that, for every $A>0$
\[
A^{1/p} \le \frac{\bar K\g_{\tau_{_A}}}{\bar \a- \bar \b -2\ve}
\]
since $\bar \a- \bar \b -2\ve>0$. Letting $A\to +\infty$ yields a contradiction since $\g_n \to 0$ as $n\to +\infty$. Hence $\tau_{_A}<+\infty$ for   large enough $A$. This completes the proof for this case.

  \medskip
 \noindent {\em Case $p>3/2$}. 
 This time, we will rely on the elementary inequality which holds for every $u$, $v\ge 0$:
 \begin{equation}\label{eq:elementairep}
 (u+v)^{2(p-1)}\le u^{2(p-1)}+ \kappa_p (u^{2p-3}v+v^{2(p-1)})\quad \mbox{ with }\quad \kappa_p = (p-1)2^{2(p-1)}.
 \end{equation}
 
Assume  $\bar \cX_\ell\!\in L^2(\P)$, $\ell=0,\ldots,k$.  Applying this inequality to $\xi_{k+\frac12}$ yields, once multiplied by $|\Delta \bar \cX_{k+1}|^2$,
 \[
 |\xi_{k+\frac12}|^{2(p-1)} |\Delta \bar \cX_{k+1}|^2\le |\bar \cX_{k}|^{2(p-1)}| \Delta \bar \cX_{k+1}|^2 +\kappa_p\big( |  \bar \cX_{k}|^{2p-3} +  | \Delta \bar \cX_{k+1}|^{2p-3} \big)| \Delta \bar \cX_{k+1}|^3.
  \]
  Plugging this inequality  into~\eqref{eq:L^pbound1}, one has
 \begin{align*}
   | \bar \cX_{k+1}|^{2p} &\le |\bar \cX_k|^{2p} +p |\bar \cX_k|^{2(p-1)}\Big( 2\big(\bar \cX_k\,|\,\Delta \bar \cX_{k+1} \big)+ (2p-1)| \Delta \bar \cX_{k+1}|^2 \Big)\\
   & \qquad + p(2p-1) \kappa_p| \Delta \bar \cX_{k+1}|^3\big( |\bar \cX_k|^{2p-3}+
    |\Delta \bar \cX_{k+1}|^{2p-3}\big)
\end{align*}
so that conditioning by $\E_k$ yields, still using that $\E_k(b_k\,|\,\Delta W_{t_{k+1}})=0$ and $\E_k |\sigma_k\Delta W_{t_{k+1}}|^2= \gamma_{k+1}\|\sigma_k\|_{_F}^2$ that
 \begin{align*}
\E_k | \bar \cX_{k+1}|^{2p} &\le  |\bar \cX_k|^{2p} + \g_{k+1}p |\bar \cX_k|^{2(p-1)} \Big( 2\big( \bar \cX_k\,|\,b_k \big)+(2p-1)\|\s_k\|_{_F}^{2}+2\g_{k+1}|b_k|^2\big)\\
&\qquad +p(2p-1) \kappa_p\big(|\bar \cX_k|^{2p-3} \E_k ( |\Delta \bar \cX_{k+1}|^3)  +
    \E_k(|\Delta \bar \cX_{k+1}|^{2p})\big).
\end{align*}

One checks that, for every $q \in [3, 2p]$,  there exists a positive real constant $ \bar C= \bar C_{b,\s,q,\bar \g}$  (with $\bar \g= \sup_n \g_n$) such that
\[
\E_k\, |\Delta \bar \cX_{k+1}|^{q} \le  \bar  C \g_{k+1}^{3/2} |\bar \cX_k|^q .
\]
Then, there exists a positive real constant $C= C_{b,\s,p,\bar \g}$ such that 
\begin{align*}
\E_k | \bar \cX_{k+1}|^{2p} &\le |\bar \cX_k|^{2p} + \g_{k+1}p |\bar \cX_k|^{2(p-1)} \Big( 2\big( \bar \cX_k\,|\,b_k \big)+(2p-1)\|\s_k\|_{_F}^{2}+2\g_{k+1}|b_k|^2\Big) +p\,C \g^{3/2}_{k+1}|\bar \cX_k|^{2p}\\
&\le |\bar \cX_k|^{2p} \\
&+ \g_{k+1} p|\bar \cX_k|^{2(p-1)} \Big(\bar K-(\bar \a-C \g_{k+1}^{1/2}-2C_b\g_{k+1}) |\bar \cX_k|^2 +(\bar \b+2C_b\g_{k+1}){\cal W}_2^2(\bar \nu_k,\d_0) \Big)\\
&\qquad +p\,C \g^{3/2}_{k+1}|\bar \cX_k|^{2p}.
\end{align*}
 
\noindent At this stage, we mimic what we  did when $p=1$ in to get~\eqref{eq:L2bound2}. Let $\ve\!\in (0, \frac{\bar \a-\bar \b}{2}) $ be small enough so that $1-p( \bar\a-\ve)>0$ and let $k_{\ve}\!\in \mathbb{N}$ be such that, for every $k\ge k_{\ve}$, 
$C \g_{k+1}^{1/2} + 2C_b\g_{k+1}<\ve$. Then 
\begin{align*}
\E_k | \bar \cX_{k+1}|^{2p} &\le |\bar \cX_k|^{2p} + p |\bar \cX_k|^{2(p-1)} \Big(\bar K-(\bar \a-\ve) |\bar \cX_k|^2 +(\bar \b+\ve) {\cal W}_2^2(\bar \nu_k,\d_0) \Big)+p\,C \g^{3/2}_{k+1}|\bar \cX_k|^{2p}.
\end{align*}
\noindent At this stage one concludes like in the previous case.
\eproof
\color{black}

\small
\bibliographystyle{plain}
\bibliography{biblio}

\appendix
\section{Appendix}

\subsection{Various C\'esaro and Kronecker Lemmas}


  {
    \color{black}
    \begin{Lemma}[Continuous time C\'esaro's Lemma]\label{le Cesaro gil}
    Let $g:\R_+\to (0,+\infty)$ in $L^1_{loc}(\R_+,d\lambda)$ such that  $\int_0^{+\infty}g(t)dt =+\infty$. Then for every bounded Borel function  $U:\R_+\to \R $ such that $\lim_{t\to+\infty} U(t)= \ell\!\in \R$,  one has
    \begin{align}\label{eq  le Cesaro gil}
        \frac{\int_0^tg(s)U(s)ds}{\int_0^tg(s)ds}\ \underset{t \rightarrow +\infty }{\longrightarrow }\ell.
    \end{align}
\end{Lemma}
}
{
    \color{black}
    \proof We may assume w.l.o.g. that $\ell=0$. Let $\ve>0$ and $A_{\ve}$ be such that $|U(t)|\le \ve $ for every $t\ge A_{\ve}$. Then 
    \[
    \bigg| \frac{\int_0^tg(s)U(s)ds}{\int_0^tg(s)ds}\bigg|\le \|U\|_{\sup}\frac{\int_0^{A_{\ve}}g(s)ds}{\int_0^tg(s)ds}+ \ve 
    \] 
so that 
$$
\varlimsup_{t\to +\infty}\bigg| \frac{\int_0^tg(s)U(s)ds}{\int_0^tg(s)ds}\bigg|\le \ve\quad\mbox{for for every} \quad \ve >0.
$$ 
   \eproof
}

%
%

\begin{Lemma}[Kronecker's Lemmas]\label{le kronecker gil}
Let $ a\ge 0$ and let $g:[a,+\infty)\to (0,+\infty)$ be a continuous  function such that  $G(t)= \int_a^{t}g(s)\ud s $ is well-defined for all $t \ge a$ and $ G(t) \to +\infty$ as $t\to +\infty$.

\smallskip 
\noindent $(a)$ {\em Continuous time regular Kronecker Lemma}.  Let $u:[a,+\infty]\to \R$ be a  continuous function. 
    Assume 
    $$
 \int_a^{t} (\log G)'(s) u(s) \ud s\underset{t \rightarrow +\infty }{\longrightarrow } I_\infty. 
    $$
    Then 
    \begin{align}\label{eq re le kronecker gil}
        \frac{\int_a^t g(s)u(s)\ud s }{\int_a^t g(s)\ud s}\underset{t \rightarrow +\infty }{\longrightarrow }0 \;.
    \end{align}
    In particular if $\rho\!\in [0,1)$ and $\int_0^{t} \frac{u(s)}{s}\ud s\to \ell \in \R$ as $t\to +\infty$, then 
    \begin{align}\label{eq re le kronecker}
        \frac{1}{t^{1-\rho}}\int_0^t s^{-\rho}u(s) \ud s \underset{t \rightarrow +\infty }{\longrightarrow }0 \;.
    \end{align}
    
    \noindent $(b)$ {\em Stochastic Kronecker Lemma}. Let $W$ be an $(\cF_t)_t$-adapted standard Brownian motion  and $(u_t)_{t\ge a}$ be  an $(\cF_t)_{t\ge a}$-adapted process such that $ \int_a^{t} \big((\log G)'(s) u(s)\big)^2 \ud s<+\infty$ $\P$-$a.s.$. Assume 
    $$
 \int_a^{t} (\log G)'(s) u(s) \ud W_s \to I_\infty, \,\R  \mbox{-valued finite random variable as} \quad t\to +\infty.
    $$ 
Then 
    \begin{align}\label{eq re le kronecker stoch}
        \frac{\int_a^t g(s)u(s)\ud W_s }{\int_a^t g(s)\ud s}\underset{t \rightarrow +\infty }{\longrightarrow }0 \;.
    \end{align}
    In particular if $\rho\!\in [0,1)$ and $\int_0^{t} \frac{u(s)}{s}\ud W_s\to I_\infty$, $\R$-valued finite random variable as $t\to +\infty$, then 
    \begin{align}\label{eq re le kronecker stoch}
        \frac{1}{t^{1-\rho}}\int_0^t s^{-\rho}u(s) \ud W_s \underset{t \rightarrow +\infty }{\longrightarrow }0 \;.
    \end{align}
\end{Lemma}

{
    \color{black}
    \proof 
        $(a)$ Let $U(t) = \int_a^t (\log G)'(s) u(s) \ud s $ for $t\ge a$. We observe that 
        \begin{align*}
           G(t)U(t) &= G(a)U(a) + \int_a^t U(s)g(s) \ud s + \int_a^t G(s)U'(s) \ud s,\\
           &=  \int_a^t U(s)g(s) \ud s + \int_a^t g(s) u(s)\ud s.
        \end{align*}
        Dividing the second equality by $G(t)$ and using the above C\'esaro Lemma~\ref{le Cesaro gil} yields
        \[  
          \frac{\int_a^t g(s)u(s)\ud s }{\int_a^t g(s)\ud s} = U(t)  - \frac{\int_a^tU(s)g(s)\ud s}{\int_a^t g(s)\ud s} \underset{t \rightarrow +\infty }{\longrightarrow } I_\infty-0-I_\infty = 0.
         \]
         
         \noindent The second claim follows by considering $g(t) = t^{-\rho}$. 
         \\
 $(b)$ Setting $U(t)=\int_a^t (\log G)'(s) u(s) \ud W_s$, the proof is the same as the one for $(a)$ since the integration by part formula is formally the same      owing to the locally  finite variation property of the function $G$.  
%
    \eproof
}

\subsection{A pseudo-Gr\"onwall Lemma}

\begin{Lemma}[Pseudo-Gr\"onwall Lemma]\label{le pseudo gronwall main} Let $\rho \in (0,1)$ and $g:\R_+\rightarrow \R_+$ and let $\psi$ be a continuous function from $\R_+ \rightarrow \R_+$ satisfying
    \begin{align}
        \label{eq prop hyp gronwall}
      \forall t \ge 0,\quad   g(t) \le \frac{\rho}t \int_0^t g(s) \ud s + \psi(t), 
          \end{align}
    \begin{enumerate}[(a)]
        \item Then, for every $t\ge 0$, 
        \begin{align}
            \label{eq prop to prove gronwall}
            g(t) \le \rho t^{\rho-1}\int_0^t s^{-\rho} \psi(s) \ud s + \psi(t).
        \end{align}
        In particular, if $\psi$ is non-decreasing, then $g(t) \le \frac1{1-\rho} \psi(t)$.
        \item Moreover if $\psi(t)=\frac1t\int_0^t f(u) \ud u$, then  for every $t\ge 0$,  
        \begin{align}
            \label{eq prop gronwall bis}
            g(t) \le t^{\rho-1}\int_0^t s^{-\rho} f(s) \ud s.
        \end{align}
        \item Moreover if $\psi(t) \rightarrow 0$ as $t\rightarrow +\infty$, then so does $g$.
    \end{enumerate}

\end{Lemma}
    
{
    \color{black}
    \proof
    $(a)$ Set $G(t) := t^{-\rho}\int_0^t g(s) \ud s$, $t>0$. We compute 
    \begin{align*}
        G'(t) &= t^{-\rho}\left(g(t)-\frac{\rho}t \int_0^t g(s) \ud s\right)
        \le t^{-\rho} \psi(t) \in L^1_{loc}(\R_+,\ud \lambda). 
    \end{align*}
    Since $g$ is continuous, we observe that $\lim_{t \downarrow 0} G(t) = 0$ and thus 
    we obtain
    \begin{align*}
        G(t) \le \int_0^t s^{-\rho} \psi(s) \ud s.
    \end{align*}
    Inserting the previous inequality into \eqref{eq prop hyp gronwall} proves \eqref{eq prop to prove gronwall}.
    
   \smallskip
   \noindent  $(b)$ Direct application of Fubini's theorem.
   
     \smallskip
        \noindent  $(c)$ Since $\psi$ goes to 0 as $t$ goes to infinity, we have that forall $\epsilon > 0$, there exists $A_\epsilon$, s.t. for $t\ge A_\epsilon$, $\psi(t) \le \epsilon$. We then compute 
    \begin{align*}
        \int_0^t s^{-\rho} \psi(s) \ud s &= \int_0^{A_\epsilon} s^{-\rho} \psi(s)  \ud s 
        + \int_{A_\epsilon}^t s^{-\rho} \psi(s)  \ud s \le C_\epsilon + \epsilon \frac{ t^{1-\rho}  }{1-\rho}.
    \end{align*}
    We thus obtain that forall $\epsilon > 0$,
    \begin{align*}
        \limsup_{t \rightarrow \infty} \left( t^{\rho-1}\int_0^t s^{-\rho} \psi(s) \ud s \right)\le \epsilon.
    \end{align*}
    Letting $\epsilon$ goes to zero concludes the proof for this step since 
    \begin{align*}
   \hskip 3,75cm     0 \le  \limsup_{t \rightarrow \infty} g(t) \le  \limsup_{t \rightarrow \infty} \left( t^{\rho-1} \int_0^t s^{-\rho} \psi(s) \ud s \right).  \hskip 3,75cm   \Box
    \end{align*} 
}

\subsection{Proof of Section \ref{subse numerical illustration}}

{
    \color{black}
    \textbf{Proof of Lemma \ref{le simple majo}} One computes, using Cauchy-Schwartz inequality, that 
    \begin{align*}
        2(b(x,\mu)-b(y,\nu)\,|\, x-y\big) &\le -2a|x-y|^2 + 2L_b|x-y|\cW_2(\mu,\nu)
    \\
        &\le -2a|x-y|^2 + \epsilon L_b|x-y|^2 + \frac{L_b}\epsilon \cW_2(\mu,\nu)^2
    \end{align*}
    So that we obtain 
    \begin{align}
        2(b(x,\mu)-b(y,\nu)\,|\, x-y\big) + \|\s(\mu)-\s(\nu)\|_{_F}^2 \le -\alpha(\epsilon)|x-y|^2+\b(\epsilon) {\cal W}_2^2(\mu,\nu).
    \end{align}
    with $\alpha(\epsilon) = 2a - \epsilon L_b $ and $\beta(\epsilon)=\frac{L_b}\epsilon + L_\s^2$.
    We are interested into the quantities:
    \begin{align*}
        \vartheta(\epsilon) = 1-\frac{\beta(\epsilon)}{\alpha(\epsilon)} \text{ and } r(\epsilon) = \frac{\vartheta(\epsilon)}{1+\vartheta(\epsilon)}
        =1-\frac{1}{1+\vartheta(\epsilon)}
    \end{align*}
    that, in the end, we want to be maximum.
    We observe that maximising $\vartheta(\epsilon)$  amounts to minimise $\frac{\beta(\epsilon)}{\alpha(\epsilon)}$.  The optimal choice is given by 
    \[
    \epsilon^\star=\frac{2a}{L_b + \sqrt{L_b^2+2 a L_b L_\s^2}} .
    \]
    \eproof

   \noindent  \textbf{Proof of Lemma \ref{le fancy drift stationary}} Note that, at equilibrium, one has 
\[
    \ud X^\star_t = (- X^\star_t +\mathfrak{b} \set{\sqrt{y^\star}-\mathfrak{d}})\ud t + \sqrt{2}\mathfrak{d} \ud W_t
\]
with $y^\star = \esp{|X^\star_t|^2}$. The stationary measure are necessarily gaussian $\cN(m^\star,z^\star)$ with $z^\star=y^\star-(m^\star)^2$. The parameter $m^\star$ and $y^\star$ are fixed point of 
\begin{align}
    m = \mathfrak{b}(\sqrt{y}-\mathfrak{d}) \text{ and } y = (\mathfrak{b}(\sqrt{y}-\mathfrak{d}))^2 + \mathfrak{d}^2.
\end{align}
We compute 
\[(1-\mathfrak{b}^2)y +2\mathfrak{b}^2\mathfrak{d} \sqrt{y}-(1+\mathfrak{b}^2)\mathfrak{d}^2 = 0 \]
and verify that if $|\mathfrak{b}|\le 1$, then $y^\star=\mathfrak{d}^2$ and $m^\star=0$. If $|\mathfrak{b}|>1$, then an extra solution is given by $\sqrt{y^\star} = \mathfrak{d}\frac{1+\mathfrak{b}^2}{\mathfrak{b}^2-1} $ and $m^\star = \frac{2\mathfrak{b}\mathfrak{d}}{\mathfrak{b}^2-1}$. In all cases, $z^\star = \mathfrak{d}^2$.
\eproof
}

\end{document}